\documentclass[11pt,a4paper,psamsfonts]{article}

\usepackage{amsmath}
\usepackage{amsfonts,amssymb}
\usepackage{enumerate}
\usepackage{amsthm}
\usepackage{hyperref} 

\makeatletter\let\orig@item=\@item \def\@item[#1]{\orig@item[\rm #1]}\makeatother
\let\origcaption=\caption \renewcommand\caption[1]{\parbox{0.66\textwidth}{\origcaption{#1}}}

\theoremstyle{plain}
\newtheorem{theo}{Theorem}[section]
\newtheorem{prop}[theo]{Proposition}
\newtheorem{lemm}[theo]{Lemma}
\newtheorem{coro}[theo]{Corollary}
\newtheorem{defi}[theo]{Definition}
\theoremstyle{definition}
\newtheorem{nota}[theo]{Notation}
\newtheorem{assu}[theo]{Assumption}
\newtheorem{rema}[theo]{Remark}
\DeclareMathOperator{\cnx}{div}
\DeclareMathOperator{\RE}{Re}
\DeclareMathOperator{\IM}{Im}
\DeclareMathOperator{\Op}{Op}

\def\CZ#1{C^{0}(H^{{#1}})}
\def\defn{\mathrel{:=}}
\def\Deltax{\Delta}
\def\Deltayx{\Delta_{x,y}}

\def\dns{\lambda_\sigma}
\def\DNS{G(\sigma)}
\def\DNZ{G(0)}

\def\eps{\varepsilon}
\def\fr{\mu}

\def\k{k}
\def\la{\left\lvert}
\def\lA{\left\lVert}
\def\le{\leq}
\def\L#1{\langle{#1}\rangle}

\def\mez{\frac{1}{2}}

\def\paradif{paradifferential }
\def\Paradif{Paradifferential }
\def\paralin{paralinearization }
\def\Paralin{Paralinearization }
\def\partialx{\nabla}
\def\partialyx{\nabla_{x,y}}
\def\periode{\ell}
\def\ra{\right\rvert}
\def\rA{\right\rVert}
\def\slam{a}
\def\Slam{A}

\def\symb{\tau}
\def\symbii{\omega}
\def\sypr{p}
\def\var{\omega}

\def\xN{\mathbf{N}}
\def\xR{\mathbf{R}}
\def\xT{\mathbf{T}}
\def\xZ{\mathbf{Z}}
\def\X2{d}

\numberwithin{equation}{section}

\begin{document}

\title{\Paralin of the Dirichlet to Neumann operator,
and regularity of three-dimensional water waves}
\author{\textsc{Thomas Alazard}
\thanks{CNRS \& Univ Paris-Sud 11. Laboratoire de Math\'ematiques d'Orsay; Orsay F-91405.
Support by the french Agence Nationale de la Recherche, project EDP Dispersives, reference ANR-07-BLAN-0250, is acknowledged.}
\and \textsc{Guy M\'etivier}
\thanks{Univ Bordeaux 1. Institut de Math\'ematiques de Bordeaux; Talence Cedex, F-33405. Support by the french
Institut Universitaire de France is acknowledged.}}

\date{}
\maketitle

\begin{abstract}
This paper is concerned with {\em a priori\/}
$C^\infty$ regularity for three-dimensional doubly periodic travelling
gravity waves whose
fundamental domain is a symmetric diamond.
The existence of such waves was a long
standing open problem solved recently by Iooss and Plotnikov.
The main difficulty is that, unlike conventional free boundary problems,
the reduced boundary system is not elliptic for three-dimensional pure gravity
waves, which leads to small divisors problems.
Our main result asserts that
sufficiently smooth diamond waves which satisfy a diophantine condition
are automatically~$C^\infty$. In particular, we prove that the solutions
defined by Iooss and Plotnikov
are $C^\infty$. Two notable
technical aspects are that (i) no smallness condition is required and (ii)
we obtain an exact paralinearization formula
for the Dirichlet to Neumann operator.
\end{abstract}

\tableofcontents
\vspace{5mm}

\section{Introduction}
The question is to prove the {\em a priori}
regularity of known travelling waves solutions to the water waves equations.
We here start an analysis
of this problem for diamond waves, which are
three-dimensional doubly periodic travelling gravity waves whose
fundamental domain is a symmetric diamond. The existence of such waves
was established by Iooss and Plotnikov in a recent beautiful memoir (\cite{IP}).

After some standard changes of unknowns which are recalled below in \S \ref{sec:equations}, for a wave
travelling in the direction $Ox_1$, we are led
to a system of two scalar equations which reads
\begin{equation}\label{reduceda}
\left\{
\begin{aligned}
&\DNS\psi -\partial_{x_1}\sigma
=0,\\
&\fr \sigma+\partial_{x_1} \psi+ \frac{1}{2}\la\partialx \psi\ra^2  -\frac{1}{2}
\frac{\bigl(\partialx  \sigma\cdot\partialx \psi +\partial_{x_1}\sigma\bigr)^2}{1+|\partialx  \sigma|^2}
= 0,
\end{aligned}
\right.
\end{equation}
where the unknowns are $\sigma,\psi\colon \xR^2\rightarrow \xR$, $\fr$ is a given positive constant and
$\DNS$ is the Dirichlet to Neumann operator, which is defined by
$$
\DNS \psi (x) =
\sqrt{1+|\partialx\sigma|^2}\,
\partial _n \phi\arrowvert_{y=\sigma(x)}
=(\partial_y \phi)(x,\sigma(x))-\partialx \sigma (x)\cdot (\partialx \phi)(x,\sigma(x)),
$$
where $\phi=\phi(x,y)$ is the solution of the Laplace equation
\begin{equation}\label{defiphi1}
\Delta_{x,y} \phi =0 \quad \text{in}\quad
\Omega\defn \{\,(x,y)\in\xR^{2}\times\xR \,\arrowvert\, y<\sigma (x)\,\},
\end{equation}
with boundary conditions
\begin{equation}\label{defiphi2}
\phi (x,\sigma(x))=\psi(x),\quad
\nabla_{x,y}\phi (x,y)\rightarrow 0 \text{ as } y\rightarrow -\infty.
\end{equation}
Diamond waves are the simplest solutions of \eqref{reduceda} one can think of. 
These 3D waves come from the nonlinear interaction of two simple oblique waves 
with the same amplitude. Henceforth, by definition, 
Diamond waves are solutions $(\sigma,\psi)$ of System~\eqref{reduceda} such that: (i)
$\sigma,\psi$ are doubly-periodic with period $2\pi$ in $x_1$ and period $2\pi\periode$ in $x_2$ for some
fixed $\ell>0$ and (ii)
$\sigma$ is even in $x_1$ and even in $x_2$; $\psi $ is odd in $x_1$ and even in $x_2$
(cf Definition~\ref{defi:diamond}).

It was proved by H. Lewy~\cite{Lewy} in the fifties that, in the two-dimensional case, if the free boundary
is a $C^{1}$ curve, then it is a $C^{\omega}$ curve (see also the independent papers of
Gerber \cite{Gerber51,Gerber55,GerberCRAS}).
Craig and Matei obtained an analogous result
for three-dimensional (i.e.\ for a 2D surface)
capillary gravity waves in~\cite{CM,CM2}.
For the study of {\em pure gravity waves} the main difficulty is that
System~\eqref{reduceda} is {\em not elliptic}.
Indeed, it is well known that $\DNZ=\la D_x\ra$ (cf \S\ref{basicexample}).
This implies that the determinant of the symbol of the
linearized system at the
trivial solution $(\sigma,\psi)=(0,0)$ is
$$
\fr | \xi | -  \xi_1^2,
$$
so that the characteristic variety $\{\xi\in \xR^2\,:\, \fr | \xi | -  \xi_1^2=0\}$ is unbounded.

This observation contains the key dichotomy between two-dimensional waves and
three-dimensional waves. Also, it explains why the problem is much more intricate for
pure gravity waves (cf \S\ref{secst} where we
prove {\em a priori} regularity for capillary waves
by using the ellipticity given by surface tension).
More importantly, it suggests that the main technical issue is that
small divisors enter into the analysis of
three-dimensional waves, as observed by Plotnikov in~\cite{Plot80} and Craig and Nicholls in~\cite{CN}.

In~\cite{IP}, Iooss and Plotnikov give a bound for the
inverse of the symbol of the linearized system
at a non trivial point under a diophantine condition,
which is the key ingredient to prove the existence of non trivial solutions to \eqref{reduceda} by means of a Nash-Moser scheme.
Our main result, which is Theorem~\ref{theo:smooth}, asserts that
sufficiently smooth diamond waves which satisfy a refined variant of their diophantine condition
are automatically~$C^\infty$. We shall prove that
there are three functions $\nu,\kappa_0,\kappa_1$ defined on the set of $H^{12}$ diamond waves such that,
if for some $0\le \delta<1$ there holds
\begin{equation*}
\la \k_2  - \left( \nu(\fr,\sigma,\psi) \k_1^2 +\kappa_0(\fr,\sigma,\psi)
+\frac{\kappa_1(\fr,\sigma,\psi)}{\k_1^{2}}\right)\ra\geq
\frac{1}{\k_1^{2+\delta}},
\end{equation*}
for all but finitely many $(\k_1,\k_2)\in\xN^2$, then
$(\sigma,\psi)\in C^{\infty}$.
Two interesting features of this result are that, firstly
{\em no smallness condition is required}, and secondly
this diophantine condition is weaker than the one
which ensures that the solutions of Iooss and Plotnikov exist.

The main corollary of this theorem states 
that the solutions of Iooss and Plotnikov are $C^\infty$.
Namely, consider the family of solutions whose existence was established in~\cite{IP}. These 
diamond waves are of the form
\begin{equation}\label{scw}
\begin{aligned}
\sigma^\eps(x)&=\eps \sigma_{1}(x) + \eps^{2}\sigma_2(x)+
\eps^{3}\sigma_3(x)
+O(\eps^4),\\
\psi^\eps(x)&=\eps \psi_1(x)+\eps^{2}\psi_2(x)+\eps^{3}\psi_3(x)+O(\eps^4),\\
\fr^\eps &= \fr_{c}+\eps^{2}\fr_{1}+O(\eps^{4}),
\end{aligned}
\end{equation}
where $\eps \in [0,\eps_0]$ is a small parameter and
$$
\fr_c\defn \frac{\periode}{\sqrt{1+\periode^2}},~
\sigma_{1}(x)\defn - \frac{1}{\fr_c} \cos x_1 \cos \left(\frac{x_2}{\periode}\right),
~ \psi_{1}(x)\defn \sin x_1 \cos \left(\frac{x_2}{\periode}\right),
$$
so that $(\sigma_{1},\psi_{1})\in C^{\infty}(\xT^2)$
solves the linearized system around the trivial solution $(0,0)$. 
Then it follows from the small divisors analysis in \cite{IP} and 
Theorem~\ref{theo:smooth} below that $(\sigma^\eps,\psi^\eps) \in C^{\infty}$.

The main novelty is to perform a full paralinearization of System~\eqref{reduceda}.
A notable technical aspect is that we obtain exact identities with remainders having
optimal regularity.
This approach depends on a careful study of the Dirichlet to Neumann operator,
which is inspired by the important paper of Lannes~\cite{LannesJAMS}.
The corresponding result about the \paralin of the Dirichlet to Neumann operator
is stated in Theorem~\ref{theo:paraDN}. This strategy has a number of consequences.
For instance, we shall see that this approach
simplifies the analysis of the diophantine condition (see Remark~\ref{numbertheory} in \S\ref{SIP}).
Also, one might in a future work use Theorem~\ref{theo:paraDN} to prove 
the existence of the solutions without the Nash--Moser iteration scheme.
These observations might be useful in a wider context.
Indeed, it is easy to prove a
variant of Theorem~\ref{theo:paraDN} for time-dependent free boundaries.
With regards to the analysis of the Cauchy problem for the water waves, this tool reduces the proof of some
difficult nonlinear estimates to easy symbolic calculus questions for symbols.

\section{Main results}

\subsection{The equations}\label{sec:equations}

We denote the spatial variables
by $(x, y) =(x_1, x_{2},y)\in \xR^{2}\times\xR$ and use the notations
\begin{align*}
\partialx=(\partial_{x_1},\partial_{x_{2}}),\quad
\Deltax=\partial_{x_1}^2+\partial_{x_2}^2,\quad
\partialyx=(\partialx,\partial_y),\quad
\Deltayx=\partial_y^2+\Deltax.
\end{align*}
We consider a three-dimensional gravity wave travelling with velocity $c$ on the free surface of an infinitely deep fluid. 
Namely, we consider a solution of the 
three-dimensional incompressible Euler equations for an irrotational flow 
in a domain of the form
$$
\Omega=\{\,(x,y)\in\times\xR^{2}\times\xR \,\arrowvert\, y<\sigma (x)\,\},
$$
whose boundary is a free surface, which means that $\sigma$
is an unknown (think of an interface between air and water).
The fact that we consider an incompressible, irrotational flow
implies that the velocity field is the gradient of a potential 
which is an harmonic function. 
The equations are then given by two boundary conditions: 
a kinematic condition 
which states that the free surface moves with the fluid, 
and a dynamic condition that expresses a balance of forces across the free 
surface. 
The classical system reads
\begin{equation}\label{System}
\left\{
\begin{aligned}
&\partial_{y}^2\phi+\Deltax \phi=0 &&\text{in } \Omega,\\
&\partial_y \phi -\partialx\sigma\cdot\partialx\phi-c\cdot\partialx\sigma=0&&\text{on }
\partial\Omega,\\
&g \sigma +\frac{1}{2}\la \partialx\phi\ra^2 +\frac{1}{2}( \partial_y\phi)^2
+c\cdot \partialx \phi
=0&&\text{on }\partial\Omega,\\
&(\partialx\phi,\partial_y \phi) \rightarrow (0,0)\quad&&\text{as } y\rightarrow -\infty,
\end{aligned}
\right.
\end{equation}
where the unknowns are $\phi\colon \Omega\rightarrow \xR$ 
and $\sigma\colon \xR^2\rightarrow \xR$, 
$c\in\xR^2$ is the wave speed and $g>0$ is the acceleration of gravity.

A popular form of the water waves equations is obtained by working
with the trace of $\phi$ at the free boundary. 
Define $\psi\colon \xR^2\rightarrow \xR$ by
$$
\psi(x)\defn\phi(x,\sigma(x)).
$$
The idea of introducing $\psi$ goes back to Zakharov.
It allows us to reduce the problem to the analysis of a system of two equations on $\sigma$ and $\psi$
which are defined on $\xR^2$. The most direct computations
show that $(\sigma,\psi)$ solves
\begin{equation*}
\left\{
\begin{aligned}
&\DNS\psi -c\cdot\partialx \sigma
=0,\\
&g \sigma+c\cdot\partialx
\psi+ \frac{1}{2}\la\partialx \psi\ra^2  -\frac{1}{2}
\frac{\bigl(\partialx  \sigma\cdot\partialx \psi +c\cdot\partialx\sigma\bigr)^2}{1+|\partialx  \sigma|^2}
= 0.
\end{aligned}
\right.
\end{equation*}
Up to rotating the axes and replacing $g$ by $\fr\defn g/\la c\ra^2$
one may assume that
$$
c=(1,0),
$$
thereby obtaining System \eqref{reduceda}.

\begin{rema}
Many variations are possible. In~\S\ref{secst} we study capillary gravity waves.
Also, we consider in~\S\ref{Secgn} the case with source terms.
\end{rema}

\subsection{Regularity of three-dimensional diamond waves}
Now we specialize to the case of diamond patterns. Namely we consider solutions which are
periodic in both horizontal directions, of the form
\begin{gather*}
\sigma(x)=
\sigma(x_1+2\pi ,x_2)=\sigma\left(x_1,x_2+2\pi\periode\right),\\
\psi(x)
=\psi(x_1+2\pi ,x_2)=\psi\left(x_1,x_2+2\pi\periode\right),
\end{gather*}
and which are symmetric with respect to the direction of propagation $Ox_1$.
\begin{defi}\label{defi:diamond}
i) Hereafter, we fix $\ell>0$ and denote by $\xT^2$ the $2$-torus
$$
\xT^2 =  
( \xR /2\pi\xZ)\times (\xR/ 2\pi\periode\xZ).
$$
Bi-periodic functions on $\xR^2$ are identified with functions on $\xT^2$, so that
the Sobolev spaces of bi-periodic functions are denoted by $H^{s}(\xT^2)$ ($s\in\xR$).

ii) Given $\fr>0$ and $s>3$, the set $D_{\fr}^s(\xT^{2})$ consists of the
solutions $(\sigma,\psi)$ of System~\eqref{reduceda} which belong to $H^{s}(\xT^2)$ and which satisfy,
for all $x\in\xR^2$,
\begin{gather*}
\sigma(x)=\sigma(-x_1,x_2)=\sigma(x_1,-x_2),\\[0.5ex]
\psi(x)=-\psi(-x_1,x_2)=\psi(x_1,-x_2),
\end{gather*}
and
\begin{equation}\label{cond:x1}
1+(\partial_{x_1}\phi)(x,\sigma(x))\neq 0,
\end{equation}
where $\phi$ denotes the harmonic extension of $\psi$ defined by
\eqref{defiphi1}--\eqref{defiphi2}.

iii) The set $D^{s}(\xT^2)$
of $H^s$ diamond waves is the set of all triple $\omega=(\fr,\sigma,\psi)$ such that
$(\sigma,\psi) \in  D_{\fr}^{s}(\xT^{2})$.
\end{defi}
\begin{rema}
A first remark about these spaces is that they are not empty; at least
since 2D waves are obviously 3D waves (independent of $x_2$) and since we know
that 2D symmetric waves exist, as proved in the twenties by Levi-Civita~\cite{Civita}, Nekrasov~\cite{Nekrasov} and Struik~\cite{Struik}.
The existence of {\em really three-dimensional pure gravity waves} 
was a well known problem in the theory of surface waves. It has been
solved by Iooss and Plotnikov in~\cite{IP}.
We refer to~\cite{IP,BDM,CN,Groves}
for references and an historical survey of the
background of this problem.
\end{rema}
\begin{rema}\label{remaOx1}
Two observations are in order about~\eqref{cond:x1}, which is not an usual assumption.
We first note that~\eqref{cond:x1} 
is a natural assumption which ensures that, in the moving frame where the waves look steady, 
the first component of the velocity evaluated at the free surface does not vanish 
(cf the proof of Lemma~\ref{propBx1}, which
is the only step in which we use \eqref{cond:x1}).
On the other hand, observe that \eqref{cond:x1} is automatically
satisfied for small amplitude waves such that $\phi=O(\eps)$ in $C^1$.
\end{rema}

For all $s\geq 23$, Iooss and Plotnikov prove the existence of
$H^{s}$-diamond waves having the above form~\eqref{scw} for $\eps\in\mathcal{E}$ where
$\mathcal{E}=\mathcal{E}(s,\periode)$ has
asymptotically a full measure when $\eps$ tends to $0$
(we refer to Theorem~\ref{theo:IP} below for a precise statement).
The set $\mathcal{E}$ is the set of parameters $\eps \in [0,\eps_0]$
(with $\eps_0$ small enough) such that a diophantine condition is satisfied. The following theorem states 
that solutions satisfying a refined diophantine condition are automatically 
$C^\infty$. We postpone to the next paragraph for a statement which asserts that this condition is not empty. 
As already mentioned, a nice technical feature is that
no smallness condition is required in the following statement.

\begin{theo}\label{theo:smooth}
There exist three real-valued functions
$\nu,\kappa_0,\kappa_1$ defined on
$ D^{12}(\xT^{2})$ such that, for all $\var=(\fr,\sigma,\psi)\in D^{12}(\xT^{2})$:

i) if there exist $\delta \in [0,1[$ and $N\in \xN^*$
such that
\begin{equation}\label{diophantine}
\la \k_2 - \left(\nu(\var) \k_1^2 +
\kappa_0(\var)+\frac{\kappa_1(\var)}{\k_1^{2}}
\right)\ra \geq
\frac{1}{\k_1^{2+\delta}},
\end{equation}
for all $(\k_1,\k_2)\in \xN^2$ with $ \k_1 \geq N$,
then $(\sigma,\psi)\in C^{\infty}(\xT^{2})$.

ii) $\nu(\var)\geq 0$ and there holds the estimate
\begin{multline*}
\la \nu(\var)-\frac{1}{\fr}\ra
+\bigg\lvert\kappa_0(\var)-\kappa_{0}(\mu,0,0)\bigg\rvert + \bigg\lvert
\kappa_1(\var)
-\kappa_{1}(\mu,0,0)\bigg\rvert \\
\le C\left(\lA (\sigma,\psi)\rA_{H^{12}}+\mu+\frac{1}{\mu}\right)\lA (\sigma,\psi)\rA_{H^{12}}^2,
\end{multline*}
for some non-decreasing function $C$ independent of $(\fr,\sigma,\psi)$.
\end{theo}

\begin{rema}
i) To define the coefficients $\nu(\var),\kappa_0(\var),\kappa_{1}(\var)$ we shall
use the principal, sub-principal and sub-sub-principal symbols of the Dirichlet to Neumann operator.
This explains the reason why we need to know that $(\sigma,\psi)$ belongs at least to $H^{12}$
in order to define these coefficients.

ii) The important thing to note about the estimate is that it is second order in
$\lA (\sigma,\psi)\rA_{H^{12}}$.
This plays a crucial role to prove that small amplitude solutions exist
(see the discussion preceding Proposition~\ref{theo:uniform}).
\end{rema}


\subsection{The small divisor condition for small amplitude waves}

The properties of an ocean surface wave are easily obtained
assuming the wave has an infinitely small amplitude (linear Airy theory).
To find nonlinear waves of small amplitude, one seeks
solutions which are small perturbations of small amplitude solutions
of the linearized system at the trivial solution $(0,0)$.
To do this, a basic strategy which goes back to Stokes
is to expand the waves in a power series of the amplitude $\eps$.
In \cite{IP}, the authors use a third order nonlinear theory to find 
3D-diamond waves (this means that they
consider solutions of the form~\eqref{scw2}).
We now state the main part of their
results (see \cite{IP} for further comments).

\begin{theo}[from~\cite{IP}]\label{theo:IP}
Let $\periode >0$ and $s\geq 23$, and set $\fr_c= \frac{\periode}{\sqrt{1+\periode^2}}$.
There is a set $A\subset [0,1]$
of full measure such that, if $\fr_c \in A$ then
there exists a set $\mathcal{E}=\mathcal{E}(s,\fr_c)$ satisfying
\begin{equation}\label{pde}
\lim_{r \rightarrow 0 }\frac{2}{r^2}\int_{\mathcal{E}\cap [0,r]} t \, dt= 1,
\end{equation}
such that there exists a family of diamond waves
$(\fr^\eps,\sigma^\eps,\psi^\eps)\in D^{s}(\xT^2)$ with
$\eps\in \mathcal{E}$,
of the special form
\begin{equation}\label{scw2}
\begin{aligned}
\sigma^\eps(x)&=\eps \sigma_{1}(x) + \eps^{2}\sigma_2(x)+\eps^3 \sigma_3(x)
+\eps ^4 \Sigma^\eps(x), \\
\psi^\eps(x)&=\eps \psi_1(x)+\eps^{2}\psi_2(x)+\eps^3\psi_3(x)+\eps^4 \Psi^\eps(x), \\
\fr^\eps &= \fr_{c}+\eps^{2}\fr_{1}+O(\eps^{4}),
\end{aligned}
\end{equation}
where $\sigma_1,\sigma_2,\sigma_3,\psi_1,\psi_2,\psi_3\in H^{\infty}(\xT^2)$ with
$$
\sigma_{1}(x)= - \frac{1}{\fr_c} \cos x_1 \cos  \Big(\frac{x_2}{\periode} \Big),
\quad \psi_{1}(x)= \sin x_1 \cos \Big(\frac{x_2}{\periode} \Big),
$$
the remainders $\Sigma^\eps,\Psi^\eps$ are uniformly bounded in $H^{s}(\xT^2)$
and
$$
 \fr_1= \frac{1}{4\fr_c^3} -\frac{1}{2 \fr_c^2} -\frac{3}{4 \fr_c} + 2 +\frac{\fr_c}{2} -
\frac{9}{4(2-\fr_c)}.
$$
\end{theo}
To prove this result, the main difficulty is to give a bound for the
inverse of the symbol of the linearized system
at a non trivial point. Due to the occurence of small divisors, this is proved in \cite{IP} under a diophantine condition. 
Now, it follows from the small divisors analysis by Iooss and Plotnikov in \cite{IP} that, for all $\eps\in\mathcal{E}$, 
\begin{equation*}
\la \k_2 - \left(\nu(\fr^\eps,\sigma^\eps,\psi^\eps) \k_1^2 +
\kappa_0(\fr^\eps,\sigma^\eps,\psi^\eps)\right)\ra \geq
\frac{c}{\k_1^{2}},
\end{equation*}
for some positive constant $c$ and all $k=(k_1,k_2)$ such that $k\neq 0, k\neq \pm (1,1)$, $k\neq \pm (-1,1)$. 
As a result, Theorem~\ref{theo:smooth} implies that, for all $\eps\in \mathcal{E}$,
$$
(\sigma^\eps,\psi^\eps)\in C^{\infty}(\xT^2).
$$
The main question left open here is to prove that, in fact,
$(\sigma^\eps,\psi^\eps)$ is analytic or at least
have some maximal Gevrey regularity. 

\smallbreak
To prove Theorem~\ref{theo:IP}, the first main step in \cite{IP}
is to define approximate solutions. 
Then, 
Iooss and Plotnikov used a Nash--Moser iterative scheme
to prove that there exist exact solutions near these approximate solutions. 
Recall that the Nash method allows
to solve functional equations of the form 
$\Phi(u)=\Phi(u_0)+f$ 
in situations where there are
loss of derivatives so that one cannot apply the usual implicit function Theorem.
It is known that the solutions thus obtained are smooth provided that $f$ is smooth
(cf Theorem 2.2.2 in~\cite{HorGeodesy}). This remark raises a question: Why
the solutions constructed by Iooss and Plotnikov are not automatically smooth?
This follows from the fact that the problem depends on the parameter
$\eps$ and hence one is led to consider functional equations of the form 
$\Phi(u,\eps)=\Phi(u_0,\eps)+f$. In this context, the estimates established in~\cite{IP} allow to
prove that, for any $\ell\in \xN$, one can define solutions $(\sigma,\psi)\in C^{\ell}(\xT^2)$
for $\eps \in \mathcal{E}\cap[0,\eps_{0}]$,
for some positive constant $\eps_{0}$ depending on $\ell$.

The previous discussion raises a second question. Indeed, to prove that the solutions exist
one has to establish uniform estimates in the following sense:
one has to prove that some diophantine condition
is satisfied for all $\k$ such that $\k_1$ is greater
than a {\em fixed } integer independent of $\eps$. In \cite{IP}, the authors
establish such a result by using an ergodic argument.
We shall explain how to perform this analysis by means of our refined diophantine condition. 
This step depends in a crucial way on the fact that the estimate of
$\nu (\mu,\sigma,\psi)-\nu(\mu,0,0)$, 
$\kappa_{0}(\mu,\sigma,\psi)-\kappa_{0}(\mu,0,0)$ and 
$\kappa_{1}(\mu,\sigma,\psi)-\kappa_{1}(\mu,0,0)$ are of second order in the amplitude. 
Namely, we make the following assumption.

\begin{assu}\label{assunukappa0}
Let $\nu=\nu(\eps)$, $\kappa_0=\kappa_{0}(\eps)$ and
$\kappa_1=\kappa_1(\eps)$ be three real-valued functions defined on $[0,1]$.
In the following proposition it is assumed that
\begin{equation}\label{nukappa0}
\begin{aligned}
&\nu(\eps)= \underline\nu+ \underline\nu' \eps^2 +\eps\varphi_1(\eps^2),\\
&\kappa_0(\eps)= \underline{\kappa_0}+ \varphi_2(\eps^2) ,\\
&\kappa_1(\eps)= \underline{\kappa_1}+\varphi_3(\eps^2),
\end{aligned}
\end{equation} 
for some constants $ \underline\nu, \underline\nu',\underline{\kappa_0},\underline{\kappa_1}'$ with
$$
 \underline\nu'\neq 0,
$$
and three Lipschitz functions $\varphi_j \colon [0,1]\rightarrow \xR$ satisfying $\varphi_j(0)=0$.
\end{assu}
\begin{rema}In \cite{IP}, the authors prove that the assumption $\underline{\nu'}\neq 0$ is satisfied for 
$\nu(\eps)=\nu(\fr^\eps,\sigma^\eps,\psi^\eps)$ where $(\fr^\eps,\sigma^\eps,\psi^\eps)$
are the solutions of Theorem~\ref{theo:IP}. Assumption~\ref{assunukappa0} is satisfied by these solutions.
\end{rema}
\begin{prop}\label{theo:uniform}
Let $\delta$ and $\delta'$ be such that
$$
1>\delta>\delta'>0.
$$
Assume in addition to Assumption~\ref{assunukappa0} that
there exists $n\geq 2$ such that
\begin{equation}\label{hd}
\la \k_2 -\underline\nu \k_1^2 -\underline{\kappa_0}\ra \geq \frac{1}{\k_1^{1+\delta'}},
\end{equation}
for all $\k \in \xN^2$ with $\k_1\geq n$. 
Then there exist $K>0$, $r_0>0$, $N_{0}\in\xN$ and a set $\mathcal{A}\subset [0,1]$ satisfying
\begin{equation}\label{pda}
\forall r \in [0,r_0],\quad \frac{1}{r}\la  \mathcal{A} \cap [0,r] \ra \ge 1- K r^{\frac{\delta-\delta'}{3+\delta'}},
\end{equation}
such that, if $\eps^2\in\mathcal{A}$ and 
$\k_1\geq N_{0}$ then
\begin{equation}\label{df2}
\la \k_2 -\nu(\eps)\k_1^2 -\kappa_0(\eps)-\frac{\kappa_1(\eps)}{\k_1^2}\ra
\geq \frac{1}{\k_1^{2+\delta}},
\end{equation}
for all $k_2\in \xN$.
\end{prop}
\begin{rema}
$(i)$ It follows from the classical argument introduced by Borel in~\cite{Borel} that there exists a null
set $\mathcal{N}\subset [0,1]$ such that, for all $(\underline\nu,\underline{\kappa_0})\in ([0,1]\setminus \mathcal{N})\times [0,1]$, the inequality
\eqref{hd} is satisfied for all $(k_1,k_2)$ with $k_1$ sufficiently large.

$(ii)$ If $\mathcal{A}$ satisfies \eqref{pda} then the set $\mathcal{E}=\{\eps\in [0,1]\,:\, \eps^2\in \mathcal{A}\}$ satisfies \eqref{pde}. 
The size of set of those parameters $\eps$ such that the diophantine condition \eqref{df2} is satisfied is bigger than the size of the 
set $\mathcal{E}$ given by Theorem~\ref{theo:IP}. 
\end{rema}

Proposition~\ref{theo:uniform} is proved in Section~\ref{SIP}.

\subsection{\Paralin of the Dirichlet to Neumann operator}\label{SparaDN}
To prove Theorem~\ref{theo:smooth}, we use the strategy of Iooss and Plotnikov~\cite{IP}.
The main novelty is that we paralinearize the water waves system.
This approach depends on a careful study of the Dirichlet to Neumann operator,
which is inspired by a paper of Lannes~\cite{LannesJAMS}.

\smallbreak
Since this analysis has several applications (for instance to the study of the Cauchy problem),
we consider the general multi-dimensional case
and we do not assume that the functions have some symmetries.
We consider here a domain $\Omega$ of the form
\begin{equation*}
\Omega = \{\,(x,y)\in\xT^{d}\times\xR \,\arrowvert\,  y<\sigma (x)\,\},
\end{equation*}
where $\xT^d$ is any $d$-dimensional torus with $d\geq 1$.
Recall that, by definition, the Dirichlet to Neumann operator is the operator $\DNS$ given by
\begin{equation*}
\begin{aligned}
\DNS\psi=\sqrt{1+|\nabla\sigma|^2}\partial_n \varphi \arrowvert_{ y=\sigma(x)},
\end{aligned}
\end{equation*}
where $n$ is the exterior normal and $\varphi$ is given by
\begin{equation}\label{defi:phi2}
\Delta_{x,y} \varphi =0, \quad \varphi\arrowvert_{y=\sigma(x)} =\psi,\quad
\nabla_{x,y}\varphi \rightarrow 0 \text{ as } y\rightarrow -\infty.
\end{equation}
To clarify notations, the Dirichlet to Neumann operator is defined by
\begin{equation}\label{defi:dn}
(\DNS\psi)(x)=
(\partial_y \varphi)(x,\sigma(x))-\partialx\sigma(x)\cdot(\partialx\varphi)(x,\sigma(x)).
\end{equation}
Thus defined, $\DNS$ differs from the usual definition of
the Dirichlet to Neumann operator
because of the scaling factor $\sqrt{1+|\nabla\sigma|^2}$; yet, as in \cite{LannesJAMS,IP}
we use this terminology for the sake of simplicity.

\smallbreak

It is known since Calder\'on that, if $\sigma$ is a given $C^\infty$ function, then
the Dirichlet to Neumann operator $\DNS$ is a classical pseudo-differential
operator, elliptic of order $1$ (see~\cite{ABB,SU,Taylor1,Treves}). We have
$$
\DNS\psi= \Op(\dns)\psi, 
$$
where the symbol $\dns$ has the asymptotic expansion
\begin{equation}\label{ae}
\dns(x,\xi)\sim\dns^1(x,\xi)+\dns^0(x,\xi)+\dns^{-1}(x,\xi)+\cdots
\end{equation}
where $\dns^{k}$ are homogeneous of degree $k$ in $\xi$, and
the principal symbol $\dns^1$ is elliptic of order $1$, given by
\begin{equation}\label{dmu1}
\dns^1(x,\xi)=\sqrt{(1+\la\nabla\sigma(x)\ra^2)\la\xi\ra^2-(\nabla \sigma(x)\cdot\xi)^2}.
\end{equation}
Moreover, the symbols $\dns^{0},\dns^{-1},\ldots$ are defined by induction so that one can easily check that
$\dns^{k}$ involves only derivatives of $\sigma$ of order $\le \la k\ra+2$ (see~\cite{ABB}).

There are also various results when $\sigma\not\in C^{\infty}$.
Expressing $G(\sigma)$ as a singular integral operator, it was proved
by Craig, Schanz and C. Sulem~\cite{CSS} that
\begin{equation}\label{CMi}
\sigma\in C^{k+1},~\psi\in H^{k+1} \text{ with }k\in \xN
\Rightarrow \DNS\psi\in H^{k}.
\end{equation}
Moreover, when $\sigma$ is a given function with
limited smoothness,
it is known that $\DNS$ is a pseudo-differential operator
with symbol of limited regularity\footnote{We do not explain here the way we define pseudo-differential operators
with symbols of limited smoothness since this problem will be fixed by using
paradifferential operators, and since all that matters in \eqref{opfr}
is the regularity of the remainder term $r(\sigma,\psi)$.} (see \cite{Taylor3,ES}).
In this direction,
for $\sigma\in H ^{s+1}(\xT^2)$ with $s$ large enough,
it follows from the analysis by Lannes (\cite{LannesJAMS}) and a small additional
work that
\begin{equation}\label{opfr}
\DNS\psi = \Op(\dns^1)\psi+r(\sigma,\psi),
\end{equation}
where the remainder $r(\sigma,\psi)$ is such that
$$
\psi\in H^{s}(\xT^d) \Rightarrow r(\sigma,\psi) \in H^{s}(\xT^d).
$$
For the analysis of the water waves,
the think of great interest here is that this gives a result for $\DNS \psi$
when $\sigma$ and $\psi$ have exactly the same regularity. Indeed, \eqref{opfr}
implies that, if
$\sigma\in H^{s+1}(\xT^d)$ and $\psi\in H^{s+1}(\xT^d)$ for some $s$ large enough,
then $\DNS\psi\in H^{s}(\xT^d)$. This result was first established by 
Craig and Nicholls in \cite{CN} and Wu in \cite{WuInvent, WuJAMS} by different methods. 
We refer to \cite{LannesJAMS} for comments on the estimates
associated to these regularity results as well as to \cite{ASL} for the rather different case
where one considers various dimensional parameters.

\smallbreak

A fundamental difference with these results is that we shall 
determine the full structure of $\DNS$ by performing a full \paralin of $\DNS \psi$ with respect to $\psi$
and $\sigma$. A notable technical aspect is that we obtain exact
identities where the remainders have optimal regularity. We shall establish a formula of the form
\begin{equation*}
\DNS\psi = \Op(\dns)\psi+B(\sigma,\psi)+R(\sigma,\psi),
\end{equation*}
where $B(\sigma,\psi)$ shall be given explicitly and $R(\sigma,\psi)\sim 0$ in the following sense: 
$R(\sigma,\psi)$ is twice more regular than $\sigma$ and $\psi$.

\smallbreak

Before we state our result, two observations are in order.

Firstly, observe that we can extend the definition of
$\dns$ for $\sigma\not\in C^\infty$
in the following obvious manner: we consider in the asymptotic expansion
\eqref{ae} only the terms which are meaningful. This means that, for
$\sigma\in C^{k+2}\setminus C^{k+3}$ with $k\in \xN$, we set
\begin{equation}\label{aef}
\dns(x,\xi)=\dns^1(x,\xi)+\dns^0(x,\xi)+\cdots+\dns^{-k}(x,\xi).
\end{equation}
We associate operators to these symbols by means of the paradifferential quantization
(we recall the definition of \paradif operators in~\S\ref{sectionpara}).

\smallbreak
Secondly, recall that a classical idea in free boundary problems is to use a
change of variables to reduce the problem to
a fixed domain. This suggests to map the graph domain $\Omega$ to
a half space via the correspondence
$$
(x,y)\mapsto (x,z) \quad\text{where}\quad z=y-\sigma(x).
$$
This change of variables
takes $\Delta_{x,y}$ to a
strictly elliptic operator and $\partial_{n}$ to a vector field which is transverse to the boundary $\{z=0\}$.
Namely, introduce
$v\colon \xT^d\times ]-\infty,0] \rightarrow \xR$ defined by
\begin{equation*}
v(x, z)  = \varphi (x, z + \sigma (x) ),
\end{equation*}
so that $v$
satisfies
$$
v\arrowvert_{z=0}= \varphi\arrowvert _{y=\sigma(x)}=\psi,
$$
and
\begin{equation}\label{dnint1a}
(1+|\partialx  \sigma|^2 ) \partial_z ^2 v+
\Deltax  v -2 \partialx  \sigma \cdot\partialx \partial_z v - \partial_z v \Delta\sigma=0,
\end{equation}
in the fixed domain $\xT^d\times ]-\infty,0[$. Then,
\begin{equation}\label{dnam.6a}
\DNS\psi=(1+|\partialx  \sigma|^2)\partial_z v
-\partialx \sigma\cdot\partialx  v \Big\arrowvert_{z=0}.
\end{equation}
Since $v$ solves the strictly elliptic equation \eqref{dnint1a} with the Dirichlet boundary condition
$v\arrowvert_{z=0}=\psi$, there is a clear
link between the regularity of $\psi$ and the regularity of $v$.
We formulate this link in Remark~\ref{rema:reg} below.
However, to state our result, the assumptions are better formulated in terms of $\sigma$ and~$v$.
Indeed, this enables us to state a result which remains valid for
the case of finite depth. The trick is that,
even if $v$ is defined for $(x,z)\in \xT^d\times ]-\infty,0]$,
we shall make an assumption on $v\arrowvert _{\xT^d \times [-1,0]}$ only
(we can replace $-1$ by any negative constant).
Below, we denote by
$C^0 ([-1,0];H^r(\xT^d))$ the space of functions which are
continuous in $z\in [-1,0]$ with values in $H^r(\xT^d)$.

\begin{theo}\label{theo:paraDN}
Let $d\geq 1$ and  $s\ge 3+d/2$ be such that $s-d/2\not\in \xN$.
If
\begin{equation}\label{assu:para}
\sigma\in H^s(\xT^{d}),\quad
v\in C^{0}([-1,0];H^s(\xT^d)),\quad \partial_{z} v \in C^{0}([-1,0];H^{s-1}(\xT^d)),
\end{equation}
then
\begin{equation}\label{paralinearize}
\DNS\psi= T_{\dns} \big( \psi-T_{\mathfrak{b}} \sigma\big) - T_{V}\cdot\partialx\sigma -T_{\cnx V} \sigma
+R(\sigma,\psi),
\end{equation}
where $T_{a}$ denotes
the paradifferential operator with symbol $a$
(cf \S\ref{sectionpara}),
the function $\mathfrak{b}=\mathfrak{b}(x)$ and the vector field $V=V(x)$ belong to
$H^{s-1}(\xT^{d})$, the symbol $\dns\in \Sigma^1_{s-1-d/2}(\xT^d)$ (see Definition~\ref{defi:Sigmamrho})
is given by $\eqref{aef}$ applied with $k=s-2-d/2$,
and $R(\sigma,\psi)$ is twice more regular than the unknowns:
\begin{equation}\label{2R}
\forall\eps >0,\quad R(\sigma,\psi)\in H^{2s-2-\frac{d}{2}-\eps}(\xT^{d}).
\end{equation}
Explicitly, $\mathfrak{b}$ and $V$ are given by
$$
\mathfrak{b}=  \frac{\partialx \sigma \cdot\partialx \psi+ \DNS\psi}{1+|\partialx  \sigma|^2},
\qquad
V=\partialx \psi -\mathfrak{b} \partialx\sigma.
$$
\end{theo}

There are a few further points that should be added to Theorem~\ref{theo:paraDN}.
\begin{rema}\label{rema:reg}
The first point to be made is a clarification of how one passes from an assumption on
$(\sigma,v)$ to an assumption on $(\sigma,\psi)$. As in \cite{ASL}, it follows 
from standard elliptic theory that
$$
\sigma\in H^{k+\frac{1}{2}}(\xT^d),~ \psi\in H^{k}(\xT^d)\Rightarrow
v\in H^{k+\frac{1}{2}}([-1,0]\times \xT^d),
$$
so that $v\in C^{0}([-1,0];H^k(\xT^d))$ and $\partial_{z} v \in C^{0}([-1,0];H^{k-1}(\xT^d))$. 
As a result, we can replace \eqref{assu:para} by the assumption that
$\sigma\in H^{s+\frac{1}{2}}(\xT^d)$ and $\psi\in H^{s}(\xT^d)$.
\end{rema}

\begin{rema}
Theorem~\ref{theo:paraDN} still holds true for non periodic functions.
\end{rema}
\begin{rema}\label{rema:fd}
The case with which we are chiefly concerned is that of an infinitely deep
fluid. However, it is worth remarking that Theorem~\ref{theo:paraDN} remains valid
in the case of finite depth where one considers a domain $\Omega$ of the form
\begin{equation*}
\Omega\defn \{\,(x,y)\in\xT^{d}\times\xR \,\arrowvert\, b(x)< y<\sigma (x)\,\},
\end{equation*}
with the assumption that $b$ is a given smooth function such that $b+2 \le \sigma$,
and define $\DNS\psi$ by \eqref{defi:dn} where $\varphi$ is given by
$$
\Delta_{x,y} \varphi =0,\quad \varphi\arrowvert_{y=\sigma(x)} =\psi,
\quad \partial_n \varphi\arrowvert_{y=b(x)}=0.
$$
\end{rema}
\begin{rema}\label{rH}
Since the scheme of the proof of Theorem~\ref{theo:paraDN} is reasonably 
simple, the reader should be able
to obtain further results in other scales of Banach spaces without too much work.
We here mention an analogous result in H\"older spaces $C^{s}(\xR^d)$
which will be used in \S\ref{secst}. If
$$
\sigma\in C^s(\xR^{d}),\quad 
v\in C^{0}([-1,0];C^s(\xR^d)),\quad \partial_{z} v \in C^{0}([-1,0];C^{s-1}(\xR^d)),
$$
for some $s\in [3,+\infty]$, then we have \eqref{paralinearize} with
$$
\mathfrak{b}\in C^{s-1}(\xR^d), \quad
V \in C^{s-1}(\xR^{d}), \quad
\dns \in \Sigma^1_{s-1}(\xR^d),
$$
and $R(\sigma,\psi)\in C^{2s-2-\eps}(\xR^{d})$ for all $\eps>0$.
\end{rema}

\begin{rema}We can give other expressions of the coefficients. We have
\begin{gather*}
\mathfrak{b}(x)= (\partial_y \varphi)(x,\sigma(x))=(\partial_{z} v)(x,0),
\\
V(x) = (\partialx \varphi)(x,\sigma(x))=(\partialx v)(x,0)-(\partial_{z} v)(x,0) \partialx \sigma(x),
\end{gather*}
where $\varphi$ is as defined in~\eqref{defi:phi2}. This clearly shows that
$\mathfrak{b},V\in H^{s-1}(\xT^d)$.
\end{rema}

As mentioned earlier, Theorem~\ref{theo:paraDN} has a number of consequences.
For instance, this permits us to reduce estimates for commutators with the
Dirichlet to Neumann operator  to symbolic calculus questions for symbols.
Similarly, we shall use Theorem~\ref{theo:paraDN}
to compute the effect of changes of variables by means of the paracomposition operators of Alinhac.
As shown by H\"ormander in~\cite{HorNash},
another possible application is to prove the existence of the solutions
by using elementary nonlinear functional analysis instead of using the Nash--Moser iteration scheme.

The proof of Theorem~\ref{theo:paraDN} is given in \S\ref{Secpara}.
The heart of the entire argument is a sharp \paralin of the interior equation performed
in Proposition~\ref{prop:key}. To do this, following Alinhac~\cite{Ali}, the idea is to work with the good unknown
$$
\psi-T_{\mathfrak{b}}\sigma.
$$
At first we may not expect to have to take this unknown into
account, but it comes up on its own when we compute the
linearized equations (cf~\S\ref{sec:3}). For the study of the linearized equations,
this change of unknowns amounts to introduce $\delta \psi -\mathfrak{b}\delta \sigma$.
The fact that this leads to a key cancelation was
first observed by Lannes in \cite{LannesJAMS}.

\subsection{An example}\label{basicexample}

We conclude this section by discussing a classical example which
is Example~$3$ in \cite{KN} (see~\cite{IP} for an analogous discussion).
Consider
$$
\phi=0\quad\text{and} \quad\sigma=\sigma(x_2).
$$
Then, for any $\sigma\in C^{1}$, this defines a solution of \eqref{System} with $g=0$, and no
further smoothness of the free boundary can be inferred.
Therefore, if $g=0$ (i.e.\ $\fr=0$) then there is no {\em a priori} regularity.

\smallbreak

In addition, the key dichotomy $d=1$ or $d=2$ is well illustrated by this
example. Indeed, consider the linearized system at the
trivial solution $(\sigma,\phi)=(0,0)$.
We are led to analyse the following system (cf \S\ref{sec:3}):
\begin{equation*}
\left\{
\begin{aligned}
&\Delta_{z,x}v=0 &\qquad &\text{in}\quad  z < 0 ,\\
& \partial_z v -\partial_{x_1}  \sigma  = 0
&&\text{on}\quad  z =0 ,\\
& \fr \sigma  +  \partial_{x_1} v  =0&&\text{on}  \quad   z = 0,\\
& \nabla_{z,x} v \rightarrow 0 && \text{as } z\rightarrow -\infty.
\end{aligned}
\right.
\end{equation*}
For $\sigma=0$, it is straightforward to compute the Dirichlet to
Neumann operator $\DNZ$. Indeed, we have to consider the solutions of
$(\la\xi\ra^2 -\partial_z^2) V(z)=0$,
which are bounded when $z<0$. It is clear that $V$ must be proportional to
$e^{z\la\xi\ra}$, so that $\partial_{z}V=\la\xi\ra V$.
Reduced to the boundary, the system thus becomes
\begin{equation*}
\left\{
\begin{aligned}
& |D_x| v   -  \partial_{x_1}  \sigma  = 0
&&\text{on}\quad  z =0 ,\\
& \fr \sigma  + \partial_{x_1} v  =0&&\text{on}  \quad   z = 0.
\end{aligned}
\right.
\end{equation*}
The symbol of this system is
\begin{equation}\label{pre24}
\begin{pmatrix}
| \xi |  & - i \xi_1 \\
i  \xi_1 & \fr
\end{pmatrix},
\end{equation}
whose determinant is
\begin{equation}\label{24}
\fr | \xi | -  \xi_1^2.
\end{equation}
If $d=1$ (or if $\fr<0$), this is a (quasi-)homogeneous elliptic symbol.
Yet, if $d=2$ (and $\fr>0$), the symbol~\eqref{24} is not elliptic.
It vanishes when $\fr|\xi |=\xi_1^2$, that is when $ | \xi_1 |   \ll   |\xi_2 | $.
The singularities are linked to the set $\{ \fr|\xi_2|=\xi_1^2\}$.
We thus have a Schr\"odinger equation on the boundary which may propagate
singularities for rational values of the parameter $\fr$.
This explains why, to prove regularity, some diophantine criterion is necessary.

To conclude, let us explain why surface tension simplifies the analysis.
Had we worked instead with capillary waves, the corresponding
symbol \eqref{pre24} would have read
$$
\begin{pmatrix}
| \xi |  & - i \xi_1 \\
i \xi_1 & \fr + \la\xi\ra^2
\end{pmatrix}.
$$
The simplification presents itself: this is an elliptic matrix-valued symbol
for all $\fr\in\xR$ and all $d\geq 1$.

\section{Linearization}\label{sec:3}

Although it is not essential for the rest of the paper, it helps
if we begin by examining the linearized equations.
Our goal is twofold. First we want to prepare for the \paralin of the equations.
And second we want
to explain some technical but important points related to changes of variables.

\smallbreak
We consider the system
\begin{equation*}
\left\{
\begin{aligned}
&\partial_y^2\phi+\Deltax \phi=0 &&\text{in}\quad \{y<\sigma(x)\},\\
&\partial_y \phi -\partialx\sigma\cdot\partialx\phi-c\cdot\partialx\sigma=0&&\text{on}\quad \{y=\sigma(x)\},\\
&\fr \sigma +\frac{1}{2}\la \partialx\phi\ra^2 +\frac{1}{2}( \partial_y\phi)^2
+c\cdot \partialx \phi
=0&&\text{on}\quad \{y=\sigma(x)\},\\
&\partialyx \phi \rightarrow 0
&&\text{as}\quad y\rightarrow -\infty,\\
\end{aligned}
\right.
\end{equation*}
where $\fr>0$ and $c\in\xR^2$.
We shall perform the linearization of this system.
These computations are well known. In particular it is known that the Dirichlet to Neumann operator  $\DNS$
is an analytic function of $\sigma$ (\cite{CoMe,NR1}). Moreover,
the shape derivative of $\DNS$ was computed by Lannes \cite{LannesJAMS}
(see also \cite{BLS,IP}).
Here we explain some key cancelations differently, by means of
the good unknown of Alinhac~\cite{Ali}.

\subsection{Change of variables}\label{sec:3.1}

One basic approach
toward the analysis of solutions of a boundary value problem is to flatten the boundary.
To do so, most directly, one can use the following
change of variables, involving the unknown $\sigma$,
\begin{equation}\label{simple}
z=y-\sigma(x),
\end{equation}
which means we introduce $v$ given by
$$
v(x,z)=\phi(x,z+\sigma(x)).
$$
This reduces the problem to the domain $\{ -\infty<z<0\}$.

\smallbreak
The first elementary step is to compute the equation satisfied
by the new unknown $v$ in $\{z<0\}$ as well as the boundary conditions
on $\{z=0\}$. We easily find the following result.

\begin{lemm}\label{lemm:cov1}
If $\phi$ and $\sigma$ are $C^2$, then $v(x, z)=\phi(x,z+\sigma(x))$ satisfies
\begin{align}
&(1+|\partialx  \sigma|^2 ) \partial_z ^2 v+
\Deltax  v -2 \partialx  \sigma \cdot\partialx \partial_z v
- \partial_z v \Delta\sigma=0&&\text{in}\quad z<0, \label{am.6a}\\[1ex]
&(1+|\partialx  \sigma|^2)\partial_z v
-\partialx \sigma\cdot (\partialx  v+c)=0&&\text{on}\quad z=0,\label{am.6b}\\[1ex]
&\fr \sigma+c\cdot\partialx v+ \frac{1}{2}\la\partialx v\ra^2  -\frac{1}{2}
\frac{\bigl(\partialx  \sigma\cdot(\partialx v+c)\bigr)^2}{1+|\partialx  \sigma|^2}
= 0&&\text{on}\quad z=0.\label{am.6c}
\end{align}
\end{lemm}
\begin{rema}
It might be tempting to use a general change of variables
of the form
$y=\rho(x,z)$ (as in \cite{CM,CM2,KLM,LannesJAMS}).
However, these changes of variables do not modify
the behavior of the functions on $z=0$ and hence they
do not modify the Dirichlet to Neumann operator (see the discussion in \cite{Uhlmann}).
Therefore, the fact that we use
the most simple change of variables
one can think of is an interesting feature of our approach.
\end{rema}
\begin{rema}
By following the strategy used in~\cite{IP},
a key point below is to use a change of variables
in the {\em tangential\/} variables, of the form
$x'=\chi(x)$.
In \cite{IP}, this change of variables is performed before the linearization.
Our approach goes the opposite
direction. We shall paralinearize first and then
compute the effect of this change of variables by means of paracomposition operators.
This has the advantage of simplifying the computations.
\end{rema}

\subsection{Linearized interior equation}
Introduce the operator
\begin{equation}\label{am.5}
L\defn (1+|\partialx  \sigma|^2 ) \partial_z ^2 +
\Deltax   -2 \partialx  \sigma \cdot\partialx \partial_z,
\end{equation}
and set
\begin{equation*}
\mathcal{E}(v,\sigma)\defn  L v - \Delta\sigma \partial_z v,
\end{equation*}
so that the interior equation~\eqref{am.6a} reads $\mathcal{E}(v,\sigma)=0$.
Denote by $\mathcal{E}'_v$ and $\mathcal{E}'_\sigma$,
the linearization of~$\mathcal{E}$ with respect
to~$v$ and~$\sigma$ respectively, which are given~by
\begin{align*}
\mathcal{E}'_{v}(v,\sigma)\dot{v} &\defn \lim_{\eps\rightarrow 0}
\frac{1}{\eps} \Bigl(\mathcal{E}(v+\eps\dot v,\sigma)-\mathcal{E}(v,\sigma)\Bigr),\\
\mathcal{E}'_{\sigma}(v,\sigma)\dot{\sigma}&\defn \lim_{\eps\rightarrow 0}
\frac{1}{\eps}\Bigl(\mathcal{E}(v,\sigma+\eps\dot{\sigma})-\mathcal{E}(v,\sigma)\Bigr).
\end{align*}
To linearize the equation $\mathcal{E}(v,\sigma)=0$, we use a standard remark in the comparison between
partially and fully linearized equations for systems obtained by the
change of variables $z=y-\sigma(x)$.
\begin{lemm}
There holds
\begin{equation}\label{am.7}
\mathcal{E}'_{v}(v,\sigma)\dot{v}
+\mathcal{E}'_{\sigma}(v,\sigma)\dot{\sigma}
=\mathcal{E}'_v(v,\sigma)\Bigl(\dot v- (\partial_{z}v)\dot\sigma \Bigr).
\end{equation}
\end{lemm}
\begin{proof}
See \cite{Ali} or \cite{MeK}.
\end{proof}
The identity~\eqref{am.7} was pointed out by S. Alinhac (\cite{Ali}) along with the role
of what he called ``the good unknown'' $\dot{u}$ defined by
$$
\dot{u}=\dot{v}-(\partial_{z}v)\dot{\sigma}.
$$
Since $\mathcal{E}(v,\sigma)$ is linear with respect to $v$, we have
$$
\mathcal{E}'_{v}(v,\sigma)\dot{v}=\mathcal{E}(\dot{v},\sigma)=L \dot{v}
-\Delta\sigma\partial_z \dot{v},
$$
from which we obtain the following formula for the linearized interior equation.

\begin{prop}
There holds
\begin{equation*}
(1+|\partialx  \sigma|^2 ) \partial_z ^2\dot{u} +
\Deltax \dot{u}  -2 \partialx  \sigma \cdot\partialx \partial_z\dot{u} - \Delta\sigma\partial_z \dot{u}=0,
\end{equation*}
where $\dot{u}\defn\dot{v}-(\partial_{z}v )\dot{\sigma}$.
\end{prop}

We conclude this part by making two remarks concerning the good unknown of Alinhac.

\begin{rema}
The good unknown $\dot{u}=\dot{v}-(\partial_{z}v)\dot{\sigma}$ was
introduced by Lannes \cite{LannesJAMS} in the analysis of
the linearized equations of the Cauchy problem for the water waves.
The computations of Lannes play a key role
in \cite{IP}. We have explained differently the reason
why $\dot u$ simplifies the computations by means of the
general identity~\eqref{am.7} (compare with the proof of Prop.\ 4.2 in~\cite{LannesJAMS}).
We also refer to a very recent paper by Trakhinin~(\cite{Trakhinin}) where
the author also uses the good unknown of Alinhac to study the Cauchy problem.
\end{rema}
\begin{rema}
A geometrical way to understand the role of the good unknown $\dot{v}-\partial_{z} v \dot{\sigma}$
is to note that the vector field $D_x\defn \partialx -\partialx \sigma \partial_{z}$
commutes with the interior equation \eqref{am.6a} for $v$: we have
\begin{equation*}
\Bigl(L-\Delta\sigma\partial_z \Bigr)D_x v=0.
\end{equation*}
The previous result can be checked directly. Alternatively, it follows from the identity
\begin{equation*}
\Bigl(L-\Delta\sigma\partial_z \Bigr)D_x v
=(D_x^2 + \partial_z^2 )D_x v,
\end{equation*}
and the fact that $D_x$ commutes with $\partial_z$. This explains why $\dot{u}$ is the natural unknown
whenever one solves a free boundary problem by straightening
the free boundary.
\end{rema}

\subsection{Linearized boundary conditions}
It turns out that the good unknown $\dot{u}$ is also useful to compute
the linearized boundary conditions. Indeed,
by differentiating the first boundary condition~\eqref{am.6b},
and replacing $\dot{v}$ by
$\dot{u}+(\partial_z v)\dot{\sigma}$ we obtain
\begin{multline*}
(1+|\partialx \sigma|^2)\partial_{z}\dot{u} -\partialx \sigma\cdot \partialx \dot{u}-
(c+\partialx  v - \partial_{z}v \partialx \sigma)\cdot\partialx \dot{\sigma}\\
+\dot{\sigma}\Big((1+|\partialx \sigma|^2)\partial_z^2v
-\partialx \sigma\cdot\partialx  \partial_z v\Big)=0.
\end{multline*}
The interior equation~\eqref{am.6a} for $v$ implies that
$$
(1+|\partialx \sigma|^2)\partial_z^2v-\partialx \sigma\cdot\partialx  \partial_z v
= - \cnx\Big(\partialx  v - \partial_{z}v \partialx \sigma\Big).
$$
which in turn implies that
\begin{equation*}
(1+|\partialx \sigma|^2)\partial_{z}\dot{u} -\partialx \sigma\cdot \partialx \dot{u}-
\cnx\Big(\big(c+\partialx  v - \partial_{z}v\partialx \sigma\big)\dot{\sigma}\Big)
=0.
\end{equation*}
With regards to the second boundary condition, we easily find that
$$
\mathfrak{a}\dot{\sigma}+ \big(c+\partialx  v - \partial_{z}v \partialx \sigma\big)\cdot\partialx \dot{u} =0,
$$
with $\mathfrak{a}\defn \fr+\big(c+\partialx  v - \partial_{z}v \partialx \sigma\big)\cdot\partialx \partial_{z}v$.

\smallbreak
Hence, we have the following proposition.
\begin{prop}\label{prop:34}
On $\{z=0\}$, the linearized boundary conditions are
\begin{equation}\label{bound:syst}
\left\{
\begin{aligned}
&N\dot{u}-\cnx(V \dot{\sigma})=0,\\[0.5ex]
&\mathfrak{a}\dot{\sigma}+ (V\cdot\partialx )\dot{u} =0,
\end{aligned}
\right.
\end{equation}
where $N$ is the Neumann operator
\begin{equation}\label{am.5N}
N=(1+\la\partialx \sigma\ra^2)\partial_{z}-\partialx \sigma\cdot\partialx ,
\end{equation}
and
$$
V=c+\partialx v -\partial_z v \partialx  \sigma ,\quad
\mathfrak{a}= \fr+V\cdot\partialx \partial_{z}v.
$$
\end{prop}
\begin{rema}
On $\{z=0\}$,
directly from the definition, we compute
$$
V=c+(\nabla \phi)(x,\sigma(x)).
$$
With regards to the coefficient $\mathfrak{a}$,
we have (cf~Lemma~\ref{lemm:hp})
$$
\mathfrak{a}=-(\partial_{y}P)(x,\sigma(x)).
$$
\end{rema}

\section{\Paralin of the Dirichlet to Neumann operator}\label{Secpara}

In this section we prove Theorem~\ref{theo:paraDN}.

\subsection{\Paradif calculus}\label{sectionpara}
We start with some basic reminders and a few more technical
issues about \paradif operators.

\subsubsection{Notations}

We denote by $\widehat{u}$ or $\mathcal{F}u$ the Fourier transform acting on temperate
distributions $u\in\mathcal{S}'(\xR^d)$, and in particular on periodic distributions.
The spectrum of $u$ is the support of $\mathcal{F}u$.
Fourier multipliers are defined by the formula
$$
p(D_x)u = \mathcal{F}^{-1}\left( p \mathcal{F}u\right),
$$
provided that the multiplication by $p$ is defined at least from
$\mathcal{S}(\xR^d)$ to $\mathcal{S}'(\xR^d)$;
$p(D_x)$ is the operator associated to the symbol $p(\xi)$.

According to the usual definition, for $\rho\in ]0,+\infty[\setminus \xN$, we denote by 
$C^\rho$ the space of bounded functions whose derivatives of order $[\rho]$  are uniformly H\"older continuous with 
exponent $\rho- [\rho]$. 

\subsubsection{\Paradif operators}\label{defi:symb}
The \paradif calculus was introduced by J.-M. Bony~\cite{Bony}
(see also~\cite{Hormander,MePise,Meyer,Taylor2}).
It is a quantization of symbols $a(x,\xi)$,
of degree $m$ in $\xi$ and limited regularity in $x$,
to which are associated operators denoted by $T_a$, of order $\le m$.

We consider symbols in the following classes.
\begin{defi}\label{defi:Gammamrho}
Given $\rho\geq 0$ and $m\in\xR$, $\Gamma_{\rho}^{m}(\xT^d)$ denotes the space of
locally bounded functions $a(x,\xi)$
on $\xT^d\times(\xR^d\setminus 0)$,
which are $C^\infty$ with respect to $\xi$ for $\xi\neq 0$ and
such that, for all $\alpha\in\xN^d$ and all $\xi\neq 0$, the function
$x\mapsto \partial_\xi^\alpha a(x,\xi)$ belongs to $C^{\rho}(\xT^d)$ and there exists a constant
$C_\alpha$ such that,
\begin{equation}\label{para:10}
\forall \la \xi\ra \ge \frac{1}{2},
\quad \lA \partial_\xi^\alpha a(\cdot,\xi)\rA_{C^{\rho}}\le C_\alpha
(1+\la\xi\ra)^{m-\la\alpha\ra}.
\end{equation}
\end{defi}
\begin{rema}\label{entier}
The analysis remains valid if we replace $C^\rho$ by
$W^{\rho,\infty}$ for $\rho\in \xN$.
\end{rema}

Note that
we consider symbols $a(x,\xi)$ that need not be smooth for $\xi=0$
(for instance $a(x,\xi)=\la\xi\ra^{m}$ with $m\in\xR^*$).
The main motivation for considering such symbols comes
from the principal symbol
of the Dirichlet to Neumann operator. As already mentioned, it is known that this symbol is given by
$$
\dns^1(x,\xi)\defn\sqrt{(1+\la\nabla\sigma(x)\ra^2)\la\xi\ra^2 - (\nabla\sigma(x)\cdot\xi)^2}.
$$
If $\sigma\in C^s(\xT^d)$ then this symbol belongs to $\Gamma^1_{s-1}(\xT^d)$. Of course, this
symbol is not $C^\infty$ with respect to $\xi\in\xR^d$.

\smallbreak
The consideration of the symbol $\dns$ also suggests that we
shall be led to consider pluri-homogeneous symbols.

\begin{defi}\label{defi:Sigmamrho}
Let $\rho\ge 1$, $m\in\xR$.
The classes $\Sigma^m_{\rho}(\xT^d)$ are
defined as the spaces of symbols $a$ such that
$$
a(x,\xi)=\sum_{0\le j <\rho} a_{m-j}(x,\xi),
$$
where $a_{m-j}\in \Gamma^{m-j}_{\rho-j}(\xT^d)$ is homogeneous of degree $m-j$ in $\xi$,
$C^\infty$ in $\xi$ for $\xi\neq 0$ and with regularity $C^{\rho-j}$ in $x$.
We call $a_m$ the principal symbol of $a$.
\end{defi}

The definition of \paradif operators needs two arbitrary but fixed
cutoff functions $\chi$ and $\psi$.
Introduce $\chi=\chi(\theta,\eta)$ such that $\chi$ is  a $C^\infty$
function on $\xR^d\times \xR^d\setminus 0$, homogeneous
of degree $0$ and satisfying, for $0<\eps_1<\eps_2$ small enough,
\begin{alignat*}{3}
\chi(\theta,\eta)&=1 \quad &&\text{if}\quad &&\la\theta\ra\le \eps_1\la \eta\ra,\\
\chi(\theta,\eta)&=0 \quad &&\text{if}\quad &&\la\theta\ra\geq \eps_2\la\eta\ra.
\end{alignat*}
We also introduce a $C^\infty$ function $\psi$ such that $0\le\psi\le 1$,
$$
\psi(\eta)=0\quad \text{for } \la\eta\ra\le 1,\qquad
\psi(\eta)=1\quad \text{for }\la\eta\ra\geq 2.
$$
Given a symbol $a(x,\xi)$, we then define
the \paradif operator $T_a$ by
\begin{equation}\label{defi:para}
\widehat{T_a u}(\xi)=(2\pi)^{-d}\int \chi(\xi-\eta,\eta)\widehat{a}(\xi-\eta,\eta)\psi(\eta)\widehat{u}(\eta)
\, d\eta,
\end{equation}
where
$\widehat{a}(\theta,\xi)=\int e^{-ix\cdot\theta}a(x,\xi)\, dx$
is the Fourier transform of $a$ with respect to the first variable.
We call attention to the fact that this notation is
not quite standard since $u$ and $a$ are periodic in $x$.
To clarify notations, fix $\xT^d= \xR^d/ L$ for some lattice $L$.
Then we can write \eqref{defi:para} as
\begin{equation*}
\widehat{T_a u}(\xi)=(2\pi)^{-d}\sum_{\eta \in L^*}
\chi(\xi-\eta,\eta)\widehat{a}(\xi-\eta,\eta)\psi(\eta)\widehat{u}(\eta).
\end{equation*}

Also, we call attention to the fact that, if $Q(D_x)$ is a Fourier multiplier with symbol $q(\xi)$,
then we do not have $Q(D_x)=T_q$,
because of the cut-off function~$\psi$. However, this is obviously almost true since
we have $Q(D_x)=T_q+R$ where
$R$ maps $H^t$ to $H^\infty$ for all $t\in\xR$.

\bigbreak

Recall the following definition, which is used continually in the sequel.

\begin{defi}Let $m\in\xR$.
An operator $T$ is said of order $\le m$ if, for all $s\in\xR$,
it is bounded from $H^{s+m}$ to $H^{s}$.
\end{defi}

\begin{theo}\label{lemM1.2b}
Let $m\in \xR$. If $a \in \Gamma^m_0(\xT^d)$, then $T_a$ is of order $\le m$.
\end{theo}

We refer to~\eqref{esti:quant1} below
for operator norms estimates.

We next recall the main feature of symbolic calculus, which is
a symbolic calculus lemma for composition of \paradif operators.
The basic property, which will be freely used in the sequel, is the following
$$
a \in \Gamma^m_1(\xT^d),
b  \in \Gamma^{m'}_1(\xT^d)
\Rightarrow
T_a T_b - T_{a b} \text{  is of order }\le m+m'-1.
$$
More generally, there is an asymptotic formula for the composition of two such operators,
whose main term is the pointwise product of their symbols.

\begin{theo}\label{theo:sc}
Let $m,m'\in\xR$. Consider $a\in \Gamma_{\rho}^{m}(\xT^d)$ and $b\in \Gamma^{m'}_{\rho}(\xT^d)$ where
$\rho\in ]0,+\infty[$, and set
$$
a\sharp b (x,\xi)= \sum_{\la\alpha\ra<\rho}
\frac{1}{i^\alpha\alpha!} \partial_{\xi}^\alpha a (x,\xi)
 \partial_x^\alpha b(x,\xi) \in \sum_{j<\rho}\Gamma^{m+m'-j}_{\rho-j}(\xT^d).
$$
Then, the operator $T_a T_b -T_{a\sharp b}$ is of order $\le m+m'-\rho$.
\end{theo}

Proofs of these two theorems can be found in the references cited above.
Clearly, the fact that we consider symbols which are periodic
in $x$ does not change the analysis.
Also, as noted in~\cite{MePise}, the fact that we consider symbols
which are not smooth at the origin $\xi=0$ is not a problem.
Here, since we added the extra function $\psi$ in the definition~\eqref{defi:para},
following the original definition in~\cite{Bony}, the argument is elementary:
if $a\in \Gamma^{m}_{\rho}(\xT^d)$, then
$\psi(\xi)a(x,\xi)$ belongs to the usual class of symbols.

\subsubsection{Paraproducts}

If $a=a(x)$ is a function of $x$ only, the \paradif operator $T_a$ is a called a paraproduct.
For easy reference, we recall a few results about paraproducts.

We already know from Theorem~\ref{lemM1.2b} that, if $\beta>d/2$
and $b\in H^{\beta}(\xT^d)\subset C^{0}(\xT^d)$, then
$T_b$ is of order $\le 0$ (note that this holds true if we
only assume that $b\in L^\infty$).
An interesting point is that one can extend the analysis
to the case where $b\in H^{\beta}(\xT^d)$ with $\beta<d/2$.

\begin{lemm}\label{lemM1.3}For all $\alpha\in \xR$ and all $\beta <d/2$,
\begin{equation*}
a  \in H^{\alpha}(\xT^d),~ b\in H^{\beta} (\xT^d)
\quad \Rightarrow \quad T_{b}a \in H^{\alpha + \beta-\frac{d}{2}}(\xT^d).
\end{equation*}
\end{lemm}

We also have the two following key lemmas about paralinearization.

\begin{lemm}\label{lemParalin}
For $a \in H^{\alpha}(\xT^d)$  with $\alpha>d/2$ and $F \in C^\infty$,
\begin{equation}\label{lemM1.1}
F( a) - T_{F'(a)}a \in H^{2\alpha-\frac{d}{2}}(\xT^d).
\end{equation}
For all $\alpha,\beta\in\xR$ such that $\alpha+\beta>0$,
\begin{equation}\label{lemM1.2bis}
a \in H^{\alpha}(\xT^d) ,  b \in H^{\beta}(\xT^d)
~  \Rightarrow ~ a b - T_a b - T_b a \in H^{\alpha + \beta-\frac{d}{2}}(\xT^d).
\end{equation}
\end{lemm}

There is also one straightforward consequence of Theorem~\ref{theo:sc} that
will be used below.
\begin{lemm}\label{lemm:sc}
Assume that $t> d/2$ is such that $t-d/2\not\in \xN$.
If $a\in H^{t}$ and $b\in H^{t}$,
then
$$
T_a T_b - T_{a b} \text{ is of order }\le - \left(t-\frac{d}{2}\right).
$$
\end{lemm}

\subsubsection{Maximal elliptic regularity}
In this paragraph,
we are concerned with scalar elliptic evolution equations of the form
$$
\partial_z u + T_a u =T_b u+f \qquad (z\in [-1,0], x\in\xT^d),
$$
where
$b\in \Gamma^0_0(\xT^d)$ and
$a\in \Gamma_2^1(\xT^d)$ is a first-order elliptic symbol with positive real part and with regularity
$C^2$ in $x$. 

With regards to further applications, 
we make the somewhat unconventional choice to take the Cauchy datum on $z=-1$. 
Recall that we denote by $C^0 ([-1,0];H^m(\xT^d))$ the space of continuous functions
in $z\in [-1,0]$ with values in $H^m(\xT^d)$. 
We prove that, if $f\in C^{0}([-1,0];H^{s}(\xT^d))$, then 
$u (0)\in H^{s+1-\eps}(\xT^d)$ for any $\eps>0$ (where 
$u(0)(x)=u\arrowvert_{z=0}=u(x,0)$).
This corresponds the usual gain of $1/2$ derivative for the Poisson kernel.
This result is not new.
Yet, for lack of a reference, we include a detailed analysis.

\begin{prop}\label{prop:max}
Let $r\in [0,1[$, $a(x, \xi)\in \Gamma_{1+r}^1(\xT^d)$ and $b(x,\xi)\in \Gamma^0_{0}(\xT^d)$. Assume
that there exists $c>0$ such that
$$
\forall (x,\xi)\in \xT^d \times \xR^{d},\quad \RE a (x,\xi)\geq c \la\xi\ra.
$$
If $u\in C^{1}([-1,0];H^{-\infty}(\xT^d))$ solves the elliptic evolution equation
$$
\partial_z u + T_a u =T_b u +f,
$$
with
$f\in C^{0}([-1,0];H^{s}(\xT^d))$ for some $s\in\xR$, then
$$
u(0) \in H^{s+r}(\xT^d).
$$
\end{prop}
\begin{proof}
The following proof gives the stronger conclusion that $u$ is continuous in $z\in ]-1,0]$ with values in 
$H^{s+r}(\xT^d)$. Therefore, by an elementary induction argument, we can assume without loss of generality
that $b=0$ and $u\in C^{0}([-1,0];H^{s}(\xT^d))$. In addition one can assume that 
$u(t,x,z)=0$ for $z\le -1/2$.

\smallbreak

Introduce the symbol
\begin{align*}
e(z;x,\xi)& \defn e_{0}(z;x,\xi)+e_{-1}(z;x,\xi)\\
&=
\exp\left( z a (x,\xi) \right)+\exp\left( z a (x,\xi)\right)\frac{z^2}{2i} \partial_{\xi}a(x,\xi)\cdot\partial_{x}a(x,\xi),
\end{align*}
so that $e(0;x,\xi)=1$ and 
\begin{equation}\label{eqe}
\partial_{z}e=e_{0}a+e_{-1} a+\frac{1}{i}\partial_{\xi}e_{0}\cdot\partial_{x}a.
\end{equation}
According
to our assumption that $\RE a \geq c\la\xi\ra$, we have the simple estimates
$$
(z\la \xi\ra)^\ell \exp\left( z a (x,\xi)\right) \le C_{\ell}.
$$
Therefore  
$$e_{0}\in C^{0}([-1,0];\Gamma^{0}_{1+r}(\xT^d)),\quad e_{-1}\in C^{0}([-1,0];\Gamma^{-1}_{r}(\xT^d)).$$ 
According to \eqref{eqe} and Theorem~\ref{theo:sc}, then, 
$T_{\partial_{z}e}-T_{e}T_{a}$ is of order $\le -r$.
Write
$$
\partial_{z}\left(T_{e} u \right)
=T_{e} f + F,
$$
with $F\in C^{0}([-1,0];H^{s+r}(\xT^d))$ and integrate on $[-1,0]$ to obtain
\begin{equation}\label{eq:duhamel}
T_{1}u(0)=
\int_{-1}^{0} F(y)\, dy+
\int_{-1}^{0} (T_{e} f)(y)\, dy.
\end{equation}
Since $F\in C^{0}([-1,0];H^{s+r}(\xT^d))$, the first term in the right-hand side belongs to $H^{s+r}(\xT^d)$. 
Moreover $u(0)-T_{1}u(0)\in H^{+\infty}(\xT^d)$
 and hence it remains only to study the second term in the right-hand side of \eqref{eq:duhamel}. Set
$$
\widetilde{u}(0)\defn \int_{-1}^{0} (T_{e} f)(y)\, dy.
$$
To prove that $\widetilde{u}(0)$ belongs to $H^{s+r}(\xT^d)$,
the key observation is that, since $\RE a \geq c\la\xi\ra$,
the family
$$
\left\{\,  (\la y\ra \la \xi\ra)^{r} e(y;x,\xi) \mid -1 \le y \le 0 \,\right\}
$$
is bounded in~$\Gamma^{r}_{0}(\xT^d)$.
Once this is granted, we use the following result 
(see \cite{MePise}) about operator norms estimates.
Given $s\in\xR$ and $m\in\xR$,
there is a constant $C$ such that,
for all $\symb\in \Gamma_0^m(\xT^d)$
and all $v\in H^{s+m}(\xT^d)$,
\begin{equation}\label{esti:quant1}
\lA T_\symb v\rA_{H^{s}}\le C \sup_{\la\alpha\ra\le \frac{d}{2}+1}\sup_{\la\xi\ra \ge 1/2}
\lA (1+\la\xi\ra)^{\la\alpha\ra-m}\partial_\xi^\alpha \symb(\cdot,\xi)\rA_{C^{0}(\xT^d)}
\lA v\rA_{H^{s+m}}.
\end{equation}
This estimate implies that there is a constant $K$ such that,
for all $-1\le y\le 0$ and all $v\in H^s(\xT^d)$,
$$
\lA (\la y\ra \la D_x\ra)^{r} (T_{e} v)\rA_{H^{s}}\le K
\lA v\rA_{H^{s}}.
$$
By applying this result we obtain that
there is a constant $K$ such that, for all $y\in [-1,0[$,
$$
\lA (T_{e} f)(y)\rA_{H^{s+r}}
\le \frac{K}{\la y\ra^{r}}\lA f(y)\rA_{H^{s}}.
$$
Since $\la y\ra^{-r} \in L^{1}(]-1,0[)$, this implies that 
$\widetilde{u}(0)\in H^{s+r}(\xT^d)$. 
This completes the proof.
\end{proof}

\subsection{\Paralin of the interior equation}
With these preliminaries established, we start the proof of Theorem~\ref{theo:paraDN}.
From now on we fix $s\geq 3+d/2$ such that $s-d/2\not\in \xN$, $\sigma\in H^{s}(\xT^d)$ and $\psi\in H^{s}(\xT^d)$.
As already explained,
we use the change of variables $z  = y-\sigma(x)$ to reduce the problem to the fixed domain
$$
\{(x,z)\in\xT^d\times\xR\,:\, z<0\}.
$$
That is, we set
\begin{equation*}
v(x, z)  = \varphi (x, z + \sigma (x) ),
\end{equation*}
which satisfies
\begin{equation}\label{dnint1}
(1+|\partialx  \sigma|^2 ) \partial_z ^2 v+
\Deltax  v -2 \partialx  \sigma \cdot\partialx \partial_z v - \partial_z v \Delta\sigma=0
\quad\text{in }\{z<0\},
\end{equation}
and  the following boundary condition
\begin{equation}\label{dnam.6}
(1+|\partialx  \sigma|^2)\partial_z v
-\partialx \sigma\cdot\partialx  v =
\DNS\psi  \quad\text{on } \{z=0\}.
\end{equation}
Henceforth we denote simply by $C^0 (H^r)$ the space of continuous functions
in $z\in [-1,0]$ with values in $H^r(\xT^d)$.
By assumption, we have
\begin{equation}\label{reg:sigmav}
v\in C^{0}(H^{s}),\quad
\partial_z  v\in C^{0}(H^{s-1}).
\end{equation}

There is one observation
that will be useful below.

\begin{lemm}
For $k=2,3$,
\begin{equation}\label{reg:sigmav2}
\partial_z ^k v\in C^{0}(H^{s-k}).
\end{equation}
\end{lemm}
\begin{proof}
This follows directly from the equation \eqref{dnint1}, the assumption
\eqref{reg:sigmav} and
the classical rule of product in Sobolev spaces which we recall here.
For $t_1,t_2\in \xR$, the product maps $H^{t_1}(\xT^d)\times H^{t_2}(\xT^d)$ to
$H^{t}(\xT^d)$ whenever
\begin{equation*}
t_1+t_2\geq 0,\quad t\le \min \{t_1,t_2\} \quad\text{and}\quad t\le t_1+t_2-d/2,
\end{equation*}
with the third inequality strict if $t_1$ or $t_2$ or $-t$ is equal to $d/2$ (cf \cite{Hormander}).
\end{proof}

We use the tangential  \paradif calculus, that is the \paradif
quantization $T_a$  of symbols $a(z, x, \xi)$ depending on the phase space variables
$(x, \xi) \in T^*\xT^{d}$ and possibly on the parameter $z \in [-1, 0]$.
Based on the discussion earlier,
to paralinearize the interior equation \eqref{dnint1}, it is natural
to introduce what we call the \emph{good unknown}
\begin{equation}
\label{m1.10}
u\defn v - T_{\partial_z v} \, \sigma.
\end{equation}
(A word of caution:  $\psi-T_{\mathfrak{b}}\sigma$, which is 
what we called the good unknown in the introduction, corresponds to the trace on $\{z=0\}$ of $u$.)

The following result is the key technical point.
\begin{prop}\label{prop:key}
The good unknown $u = v - T_{\partial_z v} \, \sigma$
satisfies the \paradif equation
\begin{equation}\label{2m.11}
T_{(1 + | \partialx \sigma|^2)}\partial_z^2 u  - 2T_{\partialx \sigma}\cdot \partialx \partial_z u
+\Deltax u   - T_{\Deltax\sigma}\partial_z u=f_0,
\end{equation}
where
$$
f_0\in \CZ{2s-3-\frac{d}{2}}.
$$
\end{prop}
\begin{proof}
Introduce the notations
$$
E\defn (1 + | \partialx \sigma|^2)\partial_z^2   - 2 \partialx \sigma\cdot
\partialx \partial_z    + \Deltax    - \Deltax\sigma \partial_z
$$
and
\begin{equation*}
P \defn  T_{(1 + | \partialx \sigma|^2)}\partial_z^2   - 2T_{\partialx \sigma}\cdot \partialx \partial_z
+\Deltax    - T_{\Deltax\sigma}\partial_z.
\end{equation*}
We begin by proving that $v$ satisfies
\begin{equation}\label{m1.9}
E v  - P v  -T_{\partial_{z}^2 v} \la \partialx \sigma\ra^2 +
2 T_{\partialx \partial_z v} \partialx  \sigma + T_{\partial_z v} \Deltax \sigma
\in \CZ{2s-3-\frac{d}{2}}.
\end{equation}
This follows from the paralinearization lemma~\ref{lemParalin}, which implies that
\begin{align*}
&\partialx\sigma\cdot\partialx\partial_z v
-T_{\partialx\sigma}\cdot\partialx \partial_z v
-T_{\partialx \partial_z v} \cdot\partialx \sigma  \in \CZ{2s-3-\frac{d}{2}} ,\\
& | \partialx  \sigma|^2   \partial_z^2  v   -
T _{ | \partialx  \sigma|^2  }    \partial_z^2 v
-T_{ \partial^2_z v}   | \partialx  \sigma|^2   \in  \CZ{2s-3-\frac{d}{2}},\\
&\partial_z v \Delta _x \sigma
-  T _{ \partial_z v}\Deltax\sigma -T_{\Deltax\sigma}\partial_z v\in  \CZ{2s-3-\frac{d}{2}}.
\end{align*}
We next substitute $v=u+T_{\partial_z v}\sigma$ in \eqref{m1.9}.
Directly from the definition of~$u$, we obtain
\begin{align*}
&\partial_z^2 u= \partial_z^2 v - T_{\partial_z^3 v}\sigma,\\
&\partialx \partial_z u = \partialx \partial_z v -  T_{\partial_z^2 v}\partialx \sigma
-T_{\partialx \partial_z^2 v}\sigma,\\
&\Deltax  u = \Deltax  v - T_{\partial_z v}\Deltax\sigma
+ 2 T_{\partialx \partial_z v}\cdot\partialx \sigma
-T_{\Deltax \partial_z v}\sigma.
\end{align*}
Since
$$
(1 + | \partialx \sigma|^2)\partial_z^2v   - 2 \partialx \sigma\cdot
\partialx \partial_z v   + \Deltax v   - \Deltax\sigma \partial_z v=0,
$$
by using Lemma~\ref{lemm:sc} and \eqref{reg:sigmav2} we obtain the key cancelation
\begin{equation}\label{keypoint}
T_{(1 + | \partialx \sigma|^2)}T_{\partial_z^3 v}\sigma   - 2T_{\partialx \sigma}
T_{\partialx \partial_z^2 v}\sigma
+T_{\Delta\partial_z v}\sigma    - T_{\partial_z^2 v\Delta\sigma}\sigma
\in \CZ{2s-3-\frac{d}{2}}.
\end{equation}
Then,
\begin{equation*}
Pu - P v + T_{\partial_z v}\Delta\sigma
-\left( 2 T_{\partialx \sigma} \cdot T_{\partial_z^2 v}\partialx \sigma
-2 T_{\partialx \partial_z v}\partialx \sigma\right) \in \CZ{2s-3-\frac{d}{2}},
\end{equation*}
so that
\begin{equation*}
E v  - P u  +\left(2 T_{\partialx \sigma} \cdot T_{\partial_z^2 v}\partialx \sigma-
T_{\partial_{z}^2 v} \la \partialx \sigma\ra^2 \right) \in \CZ{2s-3-\frac{d}{2}},
\end{equation*}
The symbolic calculus implies that
\begin{equation*}
2 T_{\partial_z^2 v} T_{\partialx \sigma} \cdot \partialx \sigma-
T_{\partial_{z}^2 v} \la \partialx \sigma\ra^2\in \CZ{2s-2-\frac{d}{2}}.
\end{equation*}
Which concludes the proof.
\end{proof}

\subsection{Reduction to the boundary}
As already mentioned, it is known that, if $\sigma$
is a $C^\infty$ given function, then
the Dirichlet to Neumann operator  $\DNS$ is a classical pseudo-differential operator.
The proof of this result is based on elliptic factorization. We here perform this elliptic
factorization for the equation for the good unknown. We next apply this lemma to determine
the normal derivatives of $u$ at the boundary in terms of tangential derivatives.

\smallbreak

We have just proved that
\begin{equation}\label{2m.11bis}
T_{(1 + | \partialx \sigma|^2)}\partial_z^2u   -
2T_{\partialx \sigma}\cdot\partialx \partial_z u
+\Deltax u   - T_{\Deltax\sigma}\partial_z u = f_0
\in \CZ{2s-3-\frac{d}{2}}.
\end{equation}
Set
\begin{equation*}
b = \frac{1}{ 1 + | \partialx \sigma |^2}.
\end{equation*}
Since $b\in H^{s-1}(\xT^d)$, by applying Lemma~\ref{lemm:sc}, we find that one can equivalently rewrite
equation~\eqref{2m.11bis}
as
\begin{equation}
\label{m1.11bis}
\partial_z^2 u - 2T_{b\partialx \sigma} \cdot\partialx \partial_z u + T_b \Deltax   u
- T_{b\Deltax\sigma}\partial_z u = f_1\in \CZ{2s-3-\frac{d}{2}}.
\end{equation}
Following the strategy in~\cite{Taylor1},
we shall perform a full decoupling into a forward and a backward elliptic evolution equations.
Recall that the classes
$\Sigma^m_{\rho}(\xT^d)$ have been defined in~\S\ref{defi:symb}.
\begin{lemm}\label{lemm:total}
There exist two symbols $\slam,\Slam\in \Sigma^{1}_{s-1-d/2}(\xT^d)$
such that,
\begin{equation}\label{2m2.b}
\begin{aligned}
( \partial_z - T_ \slam ) (\partial_z - T_\Slam)u =f\in \CZ{2s-3-\frac{d}{2}}.
\end{aligned}
\end{equation}
\end{lemm}
\begin{proof}
We seek $\slam$ and $\Slam$ in the form
\begin{equation}\label{defilL}
\slam(x,\xi)=\sum_{0\le j <t}\slam_{1-j}(x,\xi), \qquad
\Slam(x,\xi)=\sum_{0\le j <t}\Slam_{1-j}(x,\xi),
\end{equation}
where
$$
t\defn s-3-d/2,
$$
and
$$
\slam_{m}, \Slam_{m}\in \Gamma^{m}_{t+1+m} \qquad (-t < m\le 1).
$$
We want to solve the system
\begin{equation}\label{cascade1}
\begin{aligned}
\slam\sharp\Slam\defn\sum \slam_k \sharp \Slam_\ell&=-b \la \xi\ra^2+ r(x,\xi),\\
\slam+\Slam=\sum \slam_k+\Slam_k &=2  b  (i \partialx \sigma\cdot\xi) + b \Delta\sigma,
\end{aligned}
\end{equation}
for some admissible remainder $r\in \Gamma^{-t}_{0}(\xT^d)$.
Note that the notation $\sharp$, as given in Theorem~\ref{theo:sc}
depends on the regularity of the symbols. To clarify notations, we explicitly
set
$$
\slam \sharp \Slam \defn
\underset{\la\alpha\ra< t+\min\{k,\ell\}}{\sum\sum\sum} \,\frac{1}{i^\alpha \alpha!}
\partial_\xi^\alpha \slam_k \partial_x^\alpha \Slam_\ell.
$$

Assume that we have defined $\slam$ and $\Slam$ such that
\eqref{cascade1} is satisfied, and let us then prove the desired result~\eqref{2m2.b}.
For $r\in [1,+\infty)$, use the notation
$$
a\sharp_{r} b (x,\xi)= \sum_{\la\alpha\ra < r} \frac{1}{i^\alpha\alpha!} \partial_{\xi}^\alpha a (x,\xi)
 \partial_x^\alpha b(x,\xi).
$$
Then, Theorem~\ref{theo:sc} implies that
\begin{alignat*}{3}
&T_{\slam_1}T_{\Slam_1}-T_{\slam_1 \sharp_{s-1} \Slam_1 } \quad &&\text{is of order }
&&\le 1+1-(s-1)-\frac{d}{2}=-t,\\
&T_{\slam_1}T_{\Slam_0}-T_{\slam_1 \sharp_{s-2} \Slam_0 } \quad &&\text{is of order }
&&\le 1+0-(s-2)-\frac{d}{2}=-t,\\
&T_{\slam_0}T_{\Slam_1}-T_{\slam_0  \sharp_{s-2}\Slam_1 } \quad &&\text{is of order }
&&\le 0+1-(s-2)-\frac{d}{2}=-t,
\end{alignat*}
and, for $-t\le k,\ell \le 0$,
\begin{equation*}
T_{\slam_k}T_{\Slam_\ell}-T_{\slam_k \sharp_{s-2+\min\{ k,\ell\}} \Slam_\ell }
\text{ is of order } \le k+\ell-(s-2+\min\{ k,\ell\})-\frac{d}{2}\le -t-1.
\end{equation*}
Consequently,
$T_\slam T_\Slam - T_{\slam\sharp\Slam}$ is of order
$\le -t$. The first equation in \eqref{cascade1} then implies that
$$
T_\slam T_\Slam u - b \Delta u \in C^{0}(H^{s+t}),
$$
while the second equation directly gives
$$
\partial_z T_\Slam + T_\slam\partial_z u
- \big( 2 T_{b\partialx \sigma}\cdot \partialx\partial_z u - T_{b\Delta\sigma}\partial_z u\big)
\in C^{0}(H^{s+t}).
$$
We thus obtain the desired result \eqref{2m2.b} from \eqref{m1.11bis}.

\smallbreak
To define the symbols $\slam_m,\Slam_m$,
we remark first that $\slam \sharp \Slam  = \symb+\symbii$  with $\symbii\in \Gamma^{3-s}_{0}(\xT^d)$ and
\begin{equation}\label{defitau}
\symb\defn \underset{ \la\alpha\ra< s-3+k+\ell}{\sum\sum\sum} \,\frac{1}{i^\alpha \alpha!}
\partial_\xi^\alpha \slam_k \partial_x^\alpha \Slam_\ell.
\end{equation}
We then write $\tau=\sum \tau_m$ where $\tau_m$ is of order $m$.
Together with the second equation in~\eqref{cascade1}, this
yields a cascade of equations that allows to determine $\slam_m$ and $\Slam_m$ by induction.

\smallbreak
Namely, we determine $\slam$ and $\Slam$ as follows.
We first solve the principal system:
\begin{align*}
\slam_{1}\Slam_{1}&=-b \la \xi\ra^2,\\
\slam_{1}+\Slam_{1}&=2  i b \partialx \sigma\cdot\xi,
\end{align*}
by setting
\begin{align*}
\slam_1 (x, \xi) &=  i  b\partialx \sigma \cdot \xi
-   \sqrt{ b   | \xi |^2 -  (b\partialx \sigma \cdot \xi)^2},\\
\Slam_1(x, \xi)  &=  i  b\partialx \sigma \cdot \xi
+    \sqrt{ b   | \xi |^2 -( b\partialx \sigma \cdot \xi)^2}.
\end{align*}
Note that $b  | \xi |^2 - (b\partialx\sigma \cdot \xi) ^2  \geq  b^2 | \xi |^2$ so that
the symbols $\slam_1,\Slam_1$ are well defined and belong to $\Gamma^1_{s-1-d/2}(\xT^d)$.

We next solve the sub-principal system
\begin{align*}
&\slam_{0}\Slam_{1}+\slam_1\Slam_0+
\frac{1}{i}\partial_\xi \slam_1 \partial_x \Slam_1=0,\\
&\slam_{0}+\Slam_{0}= b\Delta\sigma.
\end{align*}
It is found that
\begin{align*}
\slam_0=\frac{i\partial_\xi \slam_1 \cdot\partial_x\Slam_1
-b\Delta\sigma \slam_1}{\Slam_1-\slam_1} ,\quad
\Slam_0=\frac{i\partial_\xi \slam_1 \cdot \partial_x \Slam_1
- b \Delta\sigma \Slam_1}{\slam_1-\Slam_1} .
\end{align*}
Once the principal and sub-principal symbols have been defined, one can
define the other symbols by induction.
By induction, for $-t+1\le m\le 0$, suppose that $\slam_1,\ldots, \slam_{m}$ and
$\Slam_1,\ldots, \Slam_{m}$ have been determined.
Then define $\slam_{m-1}$ and $\Slam_{m-1}$ by
$$
\Slam_{m-1}=-\slam_{m-1},
$$
and
$$
\slam_{m-1}=
\frac{1}{\slam_1 -\Slam_1}\sum\sum\sum\frac{1}{i^\alpha \alpha!}
\partial_\xi^\alpha \slam_k \partial_x^\alpha \Slam_\ell
$$
where the sum is over all triple $(k,\ell,\alpha)\in\xZ\times\xZ\times\xN^d$ such that
$$
m\le k\le 1,\quad m\le \ell\le 1, \quad \la\alpha\ra = k+\ell -m.
$$
By definition, one has $\slam_{m}, \Slam_{m}\in \Gamma^{m}_{t+1+m}$
for $-t+1<m\le 0$. Also, we obtain that $\tau=-b\la\xi\ra^2$ and
$\slam+\Slam=2  b  (i \partialx \sigma\cdot\xi) + b \Delta\sigma$.

This completes the proof.
\end{proof}

We now complete the reduction to the boundary.
As a consequence of the precise parametrix exhibited in Lemma~\ref{lemm:total}, we
describe the boundary value of $\partial_z u$ up to an error in
$H^{2s-2-\frac{d}{2}-0}(\xT^d)$.

\begin{coro}\label{prop:total}
Let $\eps>0$. On the boundary $\{ z = 0 \}$, there holds
\begin{equation}
\label{2m2.5}
(\partial_z u  -    T_\Slam u)\arrowvert_{z=0}   \in H^{2s-2-\frac{d}{2}-\eps}(\xT^d),
\end{equation}
where $\Slam$ is given by Lemma~$\ref{lemm:total}$.
\end{coro}
\begin{proof}
Introduce $ w \defn (\partial_z - T_\Slam)u$ and write $\slam=\slam_1+\widetilde{\slam}$
where $\slam_1\in \Gamma^1_2(\xT^d)$ is the principal symbol of $\slam$ and $\widetilde{\slam}\in \Gamma^0_{0}(\xT^d)$. 
Then $w$ satisfies
\begin{equation*}
\partial_z w  - T_{\slam_1} w = T_{\widetilde{\slam}}w+f .
\end{equation*}
Since $f\in C^{0}(H^{2s-3-\frac{d}{2}})$, and since $\RE \slam_1 \le -K \la \xi\ra$, 
Proposition~\ref{prop:max} applied with $r=1-\eps$ implies that
$$
(\partial_z u  -    T_\Slam u)\arrowvert_{z=0}  =w(0)\in H^{2s-2-\frac{d}{2}-\eps}(\xT^d).
$$
\end{proof}

\subsection{\Paralin of the Neumann boundary condition}

We now conclude the proof of Theorem~\ref{theo:paraDN}. Recall that, by definition,
\begin{equation*}
\DNS\psi=
(1 + | \partialx \sigma |^2) \partial_z v  - \partialx \sigma \cdot\partialx v
\arrowvert_{z=0}.
\end{equation*}
As before, on $ \{ z = 0 \}$ we find that
\begin{equation*}
T_{ (1 + | \partialx \sigma |^2)}  \partial_z v  + 2
 T_{ \partial_z v \partialx \sigma   } \cdot\partialx \sigma
 - T_{ \partialx \sigma}  \cdot\partialx v - T_ {\partialx v} \cdot\partialx \sigma \in
 H^{2s -2-\frac{d}{2}}(\xT^d).
\end{equation*}
We next translate this equation to the good unknown $u=v-T_{\partial_z v}\sigma $.
It is found that
\begin{equation*}
T_{(1+\la \partialx  \sigma\ra^2)} \partial_z u - T_{\partialx  \sigma}\cdot\partialx  u
-T_{\partialx  v-\partial_z v\partialx \sigma}\cdot\partialx  \sigma +
 T_{\alpha}\sigma\in H^{2s-2-\frac{d}{2}}(\xT^d),
\end{equation*}
with
$$
\alpha = (1+\la \partialx  \sigma\ra^2)\partial_z^2 v - \partialx  \sigma\cdot \partialx  \partial_z v.
$$
The interior equation for $v$ implies that
$$
\alpha = -\cnx  (\partialx  v-\partial_z v \partialx \sigma),
$$
so that
\begin{equation}\label{IVc1}
T_{(1+\la \partialx  \sigma\ra^2)} \partial_z u - T_{\partialx  \sigma}\cdot\partialx  u
-T_{\partialx  v-\partial_z v\partialx \sigma}\cdot\partialx  \sigma -T_{\cnx(\partialx  v-\partial_z v\partialx \sigma)}  \sigma\in H^{2s-2-\frac{d}{2}}.
\end{equation}
Furthermore, Corollary~\ref{prop:total} implies that
\begin{equation}\label{IVc2}
T_{ (1 + | \partialx \sigma |^2)}  \partial_z u - T_{\partialx \sigma}\cdot \partialx  u
=T_{\dns}  u  +R,
\end{equation}
with $R\in H^{2s-2-\frac{d}{2}-\eps}(\xT^d)$ and
$$
\dns=(1 + | \partialx \sigma |^2)\Slam -i \partialx \sigma\cdot  \xi.
$$
In particular, $\dns\in \Sigma^1_{s-1-d/2}(\xT^d)$ is a
complex-valued elliptic symbol of degree  $1$,
with principal symbol
$$
\dns^{1}(x,\xi)=
\sqrt{(1+|\partialx \sigma(x)|^2)| \xi |^2
-\left(\partialx \sigma (x)\cdot \xi\right)^2}.
$$
By combining \eqref{IVc1} and \eqref{IVc2}, we conclude the proof of
Theorem~\ref{theo:paraDN}.

\section{Regularity of diamond waves}

In this section, we prove Theorem~\ref{theo:smooth},
which is better formulated as follows.

\begin{theo}\label{toprove}
There exist three real-valued functions
$\nu,\kappa_0,\kappa_1$ defined on
$D^{12}(\xT^2)$ such that, for
all $\var=(\fr,\sigma,\psi)\in D^{s}(\xT^2)$ with $s\ge 12$,

\begin{enumerate}[i)]
\item
if there exist $\delta \in [0,1[$ and $N\in \xN^*$
such that
\begin{equation*}
\la \k_2 - \left(\nu(\var) \k_1^2 +
\kappa_0(\var)+\frac{\kappa_1(\var)}{\k_1^{2}}
\right)\ra \geq
\frac{1}{\k_1^{2+\delta}},
\end{equation*}
for all $(\k_1,\k_2)\in \xN^2$ with $ \k_1 \geq N$, then
$$
(\sigma,\psi)\in H^{s+\frac{1-\delta}{2}}(\xT^2),
$$
and hence $(\sigma,\psi)\in C^{\infty}(\xT^2)$ by an immediate induction argument.

\item $\nu(\var)\geq 0$ and there holds the estimate
\begin{multline*}
\la \nu(\var)-\frac{1}{\fr}\ra
+\bigg\lvert\kappa_0(\var)-\kappa_{0}(\mu,0,0)\bigg\rvert + \bigg\lvert
\kappa_1(\var)
-\kappa_{1}(\mu,0,0)\bigg\rvert \\
\le C\left(\lA (\sigma,\psi)\rA_{H^{12}}+\mu+\frac{1}{\mu}\right)\lA (\sigma,\psi)\rA_{H^{12}}^2,
\end{multline*}
for some non-decreasing function $C$ independent of $(\fr,\sigma,\psi)$.
\end{enumerate}
\end{theo}

We shall define explicitly the coefficients $\nu,\kappa_0,\kappa_1$.
The proof of the estimate is left to the reader.

\subsection{\Paralin of the full system}

From now on, we fix $\periode>0$ and $s\geq 12$ and consider a given diamond wave
$(\fr,\sigma,\psi)\in D^s(\xT^2)$. Recall that the system reads
\begin{equation}\label{FullSyst}
\left\{
\begin{aligned}
&\DNS\psi -c\cdot\partialx\sigma=0,\\
&\fr \sigma+c\cdot\partialx \psi+ \frac{1}{2}\la\partialx \psi\ra^2  -\frac{1}{2}
\frac{\bigl(\partialx  \sigma\cdot(\partialx \psi+c)\bigr)^2}{1+|\partialx  \sigma|^2}
= 0,
\end{aligned}
\right.
\end{equation}
where $c=(1,0)$ so that $c\cdot\partialx=\partial_{x_1}$. In analogy with the previous section, we introduce the coefficient
$$
\mathfrak{b}\defn\frac{\partialx \sigma \cdot(c+\partialx \psi)}{1+|\partialx  \sigma|^2},
$$
and what we call the good unknown
$$
u\defn \psi-T_{\mathfrak{b}}\sigma.
$$
A word of caution: this corresponds to the trace on $\{z=0\}$ of what we called 
the good unknown (also denoted by $u$) in the previous section.

The first main step is to paralinearize System~\eqref{FullSyst}.

\begin{prop}\label{prop:csystem}
The good unknown $u=\psi-T_{\mathfrak{b}}\sigma$ and $\sigma$ satisfy
\begin{align}
&T_{\dns } u - T_V \cdot\partialx \sigma -T_{\cnx V} \sigma  =f_1\in H^{2s-5}(\xT^2),\label{parasystem1}
\\
&T_\mathfrak{a} \sigma +  T_V \cdot\partialx u =f_2\in H^{2s-3}(\xT^2),\label{parasystem2}
\end{align}
where the symbol $\dns=\dns(x,\xi)\in \Sigma^{1}_{s-2}(\xT^2)$ is as given by Theorem~$\ref{theo:paraDN}$.
The coefficient $\mathfrak{a}=\mathfrak{a}(x)\in\xR$ and the vector field $V=V(x)\in\xR^2$ are given by
$$
V\defn c+\partialx \psi -\mathfrak{b} \partialx\sigma,
\qquad
\mathfrak{a}\defn \fr+V\cdot\partialx \mathfrak{b}.
$$
\end{prop}
\begin{rema}
The Sobolev embedding gives
$\dns\in \Sigma^{1}_{s-2}(\xT^2)$ if and only if $s\not\in \xN$. For $s\in\xN$ we only have $\dns\in \Sigma^{1}_{s-2-\eps}(\xT^2)$
for all $\eps>0$. Since this changes nothing in the
following analysis, we allow ourself to write abusively
$\dns\in \Sigma^{1}_{s-2}(\xT^2)$ for all $s\ge 12$.
\end{rema}
\begin{proof}
The main part of the result, which is \eqref{parasystem1}, follows directly
from Theorem~\ref{theo:paraDN} and the regularity result in Remark~\ref{rema:reg}.
The proof of \eqref{parasystem2} is easier. Indeed, note that for
$$
F (a, b) = \mez   \frac{ (a \cdot b)^2}{1 + | a |^2},
$$
there holds
$$
\partial_b F =   \frac{   (a \cdot b)}{1 + | a |^2}  a ,
\quad
\partial_ a F  =  \frac{  (a \cdot b) }{1 + | a |^2}  \Big( b - \frac{  (a \cdot b) }{1 + | a |^2} a \Big).
$$
Using these identities for $a = \partialx \sigma$ and $b = c+\partialx \psi$, the \paralin lemma (i.e. Lemma~\ref{lemParalin})
implies that
\begin{equation*}
  \fr \sigma +   T_{V } \cdot\partialx \psi
- T _{ V \mathfrak{b}} \cdot \partialx \sigma  \in H^{2s - 3}(\xT^2).
\end{equation*}
Since $s-1>d/2$, the product rule in Sobolev spaces successively implies that
$$
\mathfrak{b}=\frac{\partialx\sigma\cdot(c+\partialx \psi)}{1+\la\partialx\sigma\ra^2}\in H^{s-1}(\xT^2),
\quad V=c+\partialx\psi-\mathfrak{b}\partialx\sigma \in H^{s-1}(\xT^2).
$$
Then we use Lemma~\ref{lemm:sc} (applied with $t=s-1$), which implies the following:
\begin{equation*}
T _{ V \mathfrak{b}}\cdot  \partialx \sigma -   T _{ V} \cdot\partialx (T_{ \mathfrak{b}}    \sigma)
+ T_{V} \cdot T_{\partialx\mathfrak{b}} \sigma =
\big(T _{ V \mathfrak{b}}  -   T _{ V}  T_{ \mathfrak{b}} \big) \cdot\partialx   \sigma
 \in H^{2s -3}(\xT^2).
\end{equation*}
Similarly, Lemma~\ref{lemm:sc} applied with $t=s-2$ implies that 
$$
T_{V} \cdot T_{\partialx\mathfrak{b}} \sigma-  T_{V\cdot\partialx\mathfrak{b}} \sigma \in H^{2s -3}(\xT^2).
$$ 
As a corollary, with $\mathfrak{a}  = \fr +  V  \cdot \partialx \mathfrak{b}$,
there holds
\begin{equation*}
T_\mathfrak{a} \sigma+T_{V}\cdot\partialx u    \in H^{2 s - 3}(\xT^2).
\end{equation*}
This completes the proof.
\end{proof}

\subsection{The Taylor sign condition}
Let $\phi$ be the harmonic extension of $\psi$ as defined in \S\ref{sec:equations}, so that
\begin{equation}\label{Systembis}
\left\{
\begin{aligned}
&\partial_{y}^2\phi+\Deltax \phi=0 &&\text{in } \Omega,\\
&\partial_y \phi -\partialx\sigma\cdot\partialx\phi-c\cdot\partialx\sigma=0&&\text{on }
\partial\Omega,\\
&\fr \sigma +\frac{1}{2}\la \partialx\phi\ra^2 +\frac{1}{2}( \partial_y\phi)^2
+c\cdot \partialx \phi
=0&&\text{on }\partial\Omega,\\
&(\partialx\phi,\partial_y \phi) \rightarrow (0,0)\quad&&\text{as } y\rightarrow -\infty,
\end{aligned}
\right.
\end{equation}
where $\Omega = \{\,(x,y)\in\xR^{2}\times\xR \,\arrowvert\, y<\sigma (x)\,\}$. Define the pressure by
$$
P(x,y)\defn
-\fr y-\frac{1}{2}| \partialx  \phi (x,y)|^2-\frac{1}{2}(\partial_y\phi (x,y))^2
-c\cdot\nabla\phi (x,y).
$$
The following result gives the coefficient $\mathfrak{a}$ (which appeared in Proposition~\ref{prop:csystem}) 
in terms of the pressure $P$.
\begin{lemm}\label{lemm:hp} There holds
$\mathfrak{a}(x)=-(\partial_{y}P)(x,\sigma(x))$.
\end{lemm}
\begin{proof}
We have
$$
\mathfrak{a}(x)=\fr+ (c+ (\partialx  \phi) (x,\sigma(x)) )\cdot\partialx \bigl( (\partial_y \phi) (x,\sigma(x))\bigr).
$$
This yields
$$
\mathfrak{a}(x)=\fr+\Bigl( \partial_{y}\Bigl(\frac{1}{2}|\partialx \phi|^2
+c\cdot\partialx \phi\Bigr)
+(c+\partialx \phi)\cdot\partialx \sigma  \partial_{y}^2\phi\Bigr)(x,\sigma(x)).
$$
The Neumann condition
$\partial_y \phi-\partialx \sigma\cdot\partialx \phi-c\cdot\partialx\sigma =0$
implies that
\begin{align*}
\mathfrak{a}(x)=\partial_y \Big( \fr y +\frac{1}{2}|\partialx \phi|^2
+\frac{1}{2}(\partial_{y}\phi)^2+c\cdot\partialx\phi\Bigr)(x,\sigma(x))=-(\partial_y P)(x,\sigma(x)),
\end{align*}
which concludes the proof.
\end{proof}

The Taylor sign condition is the physical assumption that
the normal derivative
of the pressure in the flow at the free surface is negative.
The equation
$$
\partial_y \phi-\partialx\sigma\cdot\partialx\phi-c\cdot\partialx\phi=0,
$$
implies that $\partial_n P = \partial_{y} P$ at the free surface. Therefore, the Taylor sign condition
reads
\begin{equation}\label{tsc}
\forall x\in \xT^2,\quad (\partial_{y} P)(x,\sigma(x))<0.
\end{equation}
It is easily proved that this property is satisfied under a smallness assumption (see~\cite{IP}).
Indeed,
if $\lA (\sigma,\psi)\rA_{C^2}=O(\eps)$, then
$$
\lA (\partial_{y} P)(x,\sigma(x)) + \fr\rA_{L^\infty}=O(\eps).
$$
Our main observation is that diamond waves
satisfy \eqref{tsc} automatically: No smallness assumption is required to prove \eqref{tsc}.
This a consequence of the following general proposition, which is a variation of one of Wu's key
results (\cite{WuInvent,WuJAMS}). Since the result is not restricted
to diamond waves,
the following result is stated in a little more generality than is needed.

\begin{prop}\label{theo:WU}
Let $\fr>0$ and $c\in\xR^2$.
If $(\sigma,\phi)$ is a $C^2$ solution of~\eqref{Systembis}
which is doubly periodic in the horizontal variables $x_1$ and $x_2$,
then the Taylor sign condition
is satisfied: $(\partial_y P )(x,\sigma(x))<0$ for all $x\in\xT^2$. As a result it follows from the previous lemma that
$$
\mathfrak{a}(x)>0,
$$
for all $x\in \xT^2$.
\end{prop}
\begin{rema}
$(i)$ In the case of finite depth, it was shown by Lannes
(\cite{LannesJAMS}) that the Taylor sign condition is satisfied
under an assumption on the second fundamental form of the bottom surface
(cf Proposition $4.5$ in \cite{LannesJAMS}). 

$(ii)$ Clearly the previous result is false for $\mu=0$. Indeed, if $\mu=0$ then $(\sigma,\phi)=(0,0)$ solves~\eqref{Systembis}.
\end{rema}
\begin{proof}[Proof (from~\cite{WuJAMS,LannesJAMS})] 
We have $P=0$ on the free surface $\{y=\sigma(x)\}$. On the other hand,
since $\mu>0$ and since $\nabla_{x,y}\phi\rightarrow 0$ when $y$ tends to $-\infty$, there exists
$h>0$ such that
$$
P(x,y)\ge 1 \quad\text{for}\quad y\le -h.
$$
Define
$$
\Omega_h = \{ (x,y)\in \xR^2\times\xR \,:\, -h\le y\le \sigma(x)\}.
$$
Since $P$ is bi-periodic in $x$, $P$ reaches its minimum on $\Omega_h$.
The key observation is that the equation $\Deltayx\phi=0$ implies that
$$
\Deltayx P = - | \nabla_{y,x}\nabla_{y,x} \phi |^2 \leq 0,
$$
and hence $P$ is a super-harmonic function.
In particular $P$ reaches its minimum on $\partial \Omega_h$ and at such a point we have
$\partial_n P <0$. We conclude the proof by means of the following three ingredients: (i) $P$ reaches its
minimum on the free surface since $P\arrowvert_{y=-h}\ge 1 > 0=P\arrowvert_{y=\sigma(x)}$;
(ii) $P=0$ on the free surface so that $P$ reaches its minimum
at any point of the free surface, hence $\partial_{n}P<0$ on $\{y=\sigma(x)\}$; and (iii)
$\partial_n P= \partial_y P$ on the free surface .
\end{proof}
\begin{rema}
According to the Acknowledgments in \cite{WuJAMS}, the idea of expressing 
 $\mathfrak{a}$ as $-(\partial_{y}P)(x,\sigma(x))$ and using the maximum principle to prove that 
 $-(\partial_{y}P)(x,\sigma(x))<0$ is credited to 
R. Caflisch and Th. Hou.
\end{rema}

After this short d\'etour, we return to the main line of our development. 
By using the fact that $\mathfrak{a}$ does not vanish, one can form a
second order equation from~\eqref{parasystem1}--\eqref{parasystem2}.
\begin{lemm}\label{lem55}
Set
\begin{equation}\label{defimV}
\mathcal{V}(x,\xi) \defn - \mathfrak{a}(x)^{-1}(V(x)\cdot \xi)^2 +i \cnx \left(\mathfrak{a}(x)^{-1}(V(x)\cdot \xi) V(x)\right).
\end{equation}
Then,
\begin{equation}\label{oneeq}
T_{\dns+\mathcal{V}}\, u \in H^{2s-5}(\xT^2).
\end{equation}
\end{lemm}
\begin{rema}
The fact that $\mathfrak{a}$ is positive implies that
the symbol $\dns+\mathcal{V}$ may vanish or be arbitrarily small.
If $\mathfrak{a}$ were negative, the analysis would have been
much easier (cf Section~\ref{Secgn}).
\end{rema}
\begin{proof}As already seen we have 
$\mathfrak{b}\in H^{s-1}(\xT^2)$ and $V\in H^{s-1}(\xT^2)$. Therefore the product rule in Sobolev space implies that
$$
\mathfrak{a}=\fr+V \cdot\partialx\mathfrak{b} \in H^{s-2}(\xT^2).
$$
Then, by applying Theorem~\ref{theo:sc} with $\rho=s-3$  
we obtain that,
$$
T_{\mathcal{V}}u- \left( T_{V}\cdot\partialx + T_{\cnx V}\right)
T_{\mathfrak{a}^{-1}} T_{V} \cdot\partialx u  \in H^{s+(s -5)}(\xT^2).
$$
On the other hand, since $\mathfrak{a},\mathfrak{a}^{-1}\in H^{s-2}(\xT^2)$, Lemma~\ref{lemm:sc} implies that
$$
T_{\mathfrak{a}^{-1}}T_{\mathfrak{a}} -I\quad\text{is of order }\le -(s-3),
$$
and hence 
\begin{equation}\label{sansa}
\sigma -( - T_{\mathfrak{a}^{-1}} T_{V}\cdot\partialx u)\in H^{2s-4}(\xT^2).
\end{equation}
The desired result \eqref{oneeq} is an immediate consequence
of~\eqref{parasystem1}--\eqref{parasystem2}.
\end{proof}

\subsection{Notations}

The following notations are used continually in this section.
\begin{nota}\label{Nota}
i)
The set $C(o,e)$ is the set of function
$f=f(x_1,x_2)$ which are odd in $x_1$ and even in $x_2$. Similarly we define the sets
$C(o,o)$, $C(e,o)$ and $C(e,e)$.

ii) The set $\Gamma(o,e)$ is the set of symbols $a=a(x_1,x_2,\xi_1,\xi_2)$ 
such that
\begin{equation*}
\left\{
\begin{aligned}
&a(-x_1,x_2,-\xi_1,\xi_2)=-a(x_1,x_2,\xi_1,\xi_2),\\
&a(x_1,-x_2,\xi_1,-\xi_2)= a(x_1,x_2,\xi_1,\xi_2).
\end{aligned}
\right.
\end{equation*}
Similarly we define the sets
$\Gamma(o,o)$, $\Gamma(e,o)$ and $\Gamma(e,e)$.
\end{nota}
\begin{rema}\label{remaoe}
$(i)$ With these notations, by assumption we have $\sigma\in C(e,e)$ and $\psi\in C(o,e)$ so that
$$
\mathfrak{b}=\frac{\partial_{x_1} \sigma (1+\partial_{x_1} \psi)+\partial_{x_2}\sigma \partial_{x_2}\psi}{1+|\partialx  \sigma|^2}\in C(o,e),
$$
and hence $u=\psi-T_{\mathfrak{b}}\sigma\in C(o,e)$. Also, we check that the vector field $V=(V_1,V_2)=
c+\partialx \psi-\mathfrak{b}\partialx \sigma$ is such that 
$V_1\in C(e,e)$ and $V_2\in C(o,o)$.

$(ii)$ If $v \in C(o,e)$ and $a\in \Gamma(e,e)$ then $T_a v\in C(o,e)$ (provided that $T_a v$ is well defined). 
Clearly, the same property is true for the three other classes of symmetric functions.
\end{rema}

To simplify the presentation,
we will often only check only one half of
the symmetry properties which are claimed in the various statements below.
We will only check the symmetries with respect to the axis $\{x_1=0\}$.
To do this, it will be convenient to use the following notation
(as in \cite{IP}).
\begin{nota}\label{notaetoile}
By notation, given $z=(z_1,z_2)\in \xR^2$,
$$
z^\star=(-z_1,z_2).
$$
\end{nota}
\subsection{Change of variables}
We have just proved that $T_{\dns+\mathcal{V}}\, u \in H^{2s-5}(\xT^2)$.
We now study the sum of the principal symbols of $\dns$ and $\mathcal{V}$. Introduce
$$
p(x,\xi)
=\sqrt{(1+\la\partialx\sigma(x)\ra^2)\la\xi\ra^2-(\partialx\sigma(x)\cdot\xi)^2}
-\mathfrak{a}(x)^{-1}(V(x)\cdot\xi)^2.
$$
By following the analysis in \cite{IP},
we shall prove that there exists a change of variables $\xR^2\ni x\mapsto \chi(x)\in \xR^2$ such that
$\sypr\bigl(x,{}^t\chi'(x) \xi\bigr)$ has a simple expression.

\smallbreak

Since we need to consider change of variables $ x\mapsto \chi(x)$  such that
$u\circ \chi$ is doubly periodic whenever $u$ is doubly periodic, we introduce the following definition.

\begin{defi}\label{def:Crdiffeo}
Let $\chi\colon \xR^2\rightarrow \xR^2$ be a continuously differentiable diffeomorphism.
    For $r>1$, we say that
$\chi$ is a $C^{r}(\xT^2)$-diffeomorphism if there exists $\tilde{\chi}\in C^{r}(\xT^2)$ such that
$$
\forall x\in \xR^2, \quad \chi(x)=x+\tilde{\chi}(x).
$$
(Recall that bi-periodic functions on $\xR^2$ are identified with functions on $\xT^2$.)
\end{defi}

In this paragraph we show
\begin{prop}
\label{prop:r1}
There exist a $C^{s-4}(\xT^2)$-diffeomorphism~$\chi$, a constant $\nu\geq 0$, a positive function
$\gamma\in C^{s-4}(\xT^2)$  and a symbol $\alpha\in \Gamma^{0}_{s-6}(\xT^2)$ homogeneous of degree $0$ in $\xi$
such that, for all $(x,\xi)\in \xT^2\times\xR^2$,
\begin{equation*}
\sypr\bigl(x,{}^t\chi'(x) \xi\bigr)=\gamma(x)\bigl(\la \xi\ra - \nu \xi_1^2\bigr)
+i \alpha(x,\xi)\xi_1,
\end{equation*}
and such that the following properties hold:
\begin{enumerate}[i)]
\item $\chi=(\chi^1,\chi^2)$ where $\chi^1 \in C(o,e)$ 
and $\chi^2\in C(e,o)$; 
\item $\alpha \in \Gamma (o,e)$, $\gamma\in C(e,e)$.
\end{enumerate}
\end{prop}

Proposition~\ref{prop:r1} will be deduced from the following
lemma.
\begin{lemm}\label{propBx1}
There exists a $C^{s-3}(\xT^2)$-diffeomorphism $\chi_1$ of the form
$$
\chi_1 (x_1,x_2) = \begin{pmatrix} x_1 \\ \X2(x_1,x_2)\end{pmatrix},
$$
such that $\X2$ solves the transport equation
$$
V(x)\cdot\partialx \X2 (x)=0,
$$
with initial data $\X2(0,x_2)=x_2$ on $x_1=0$, and such that, for all $x\in \xR^2$,
\begin{align}
\X2(x_1,x_2)&=\X2(-x_1,x_2)=-\X2(x_1,-x_2), \quad (d\in C(e,o))\label{propertiesd1}\\
\X2(x_1,x_2)&=\X2(x_1+2\pi,x_2)=\X2\left(x_1,x_2+2\pi\periode\right)
-2\pi\periode\label{propertiesd2}.
\end{align}
\end{lemm}
\begin{rema}\label{remaBx1}
The result is obvious for the trivial solution $(\sigma,\psi)=(0,0)$ with $\X2(x)=x_2$.
It can be easily inferred from the analysis in Appendix~C in \cite{IP} that
this result is also satisfied in a neighborhood of $(0,0)$. The new point here
is that we prove the result under the only assumption that
$(\sigma,\phi)$ satisfies condition~\eqref{cond:x1} in Definition~\ref{defi:diamond}.
\end{rema}
\begin{proof}
Assumption~\eqref{cond:x1} implies that $V_1(x)\neq 0$ for all $x\in\xT^2$.
We first write that, if $\X2$ satisfies $V\cdot\partialx \X2=0$ with initial data $\X2(0,x_2)=x_2$ on $x_1=0$, then
$w=\partial_{x_2} \X2$ solves the Cauchy problem
\begin{equation}\label{CPw}
\partial_{x_1}w + \frac{V_2}{V_1}\partial_{x_2} w + w \partial_{x_2}\left(\frac{V_2}{V_1}\right) =0,\qquad
w(0,x_2)=1.
\end{equation}
To study this Cauchy problem, we work in the Sobolev spaces
of $2\pi\periode$-periodic functions $H^{s}(T)$ where~$T$ is the circle $T\defn \xR/(2\pi\periode\xZ)$. Since
$$
V_2/V_1, \partial_{x_2}(V_2/V_1)\in H^{s-2}(\xT^2)\subset L^{\infty}(\xR;H^{s-3}(T)),
$$
and since $s>4$, standard results for hyperbolic equations imply that \eqref{CPw} has a
unique solution $w\in C^{0}(\xR;H^{s-3}(T))$. We define~$\X2$ by
\begin{equation}\label{defid}
\X2(x_1,x_2)\defn\int_{0}^{x_2} w(x_1,t)\, dt.
\end{equation}
We then obtain a solution of $V\cdot\partialx \X2=0$, and we easily checked that
$$
\X2(x)-x_2\in \bigcap _{0\le  j <s-2}C^{j}(\xR;H^{s-2-j}(T)).
$$
The Sobolev embedding thus implies that $\X2(x)-x_2\in C^{s-3}(\xR^2)$.

\smallbreak
We next prove that $\X2$ satisfies \eqref{propertiesd1}--\eqref{propertiesd2}.
Firstly, by uniqueness for the Cauchy problem~\eqref{CPw}, using that $V_1,V_2$ are $2\pi\ell$-periodic in $x_2$ and using 
the symmetry properties of $V_1$, $V_2$ (cf Remark~\ref{remaoe}) we easily obtain
\begin{equation*}
w(x_1,x_2)=w(-x_1,x_2)=w(x_1,-x_2)
=w\left(x_1,x_2+2\pi\periode\right).
\end{equation*}
To prove that $w$ is periodic in $x_1$, following \cite{IP}, we use in a essential way
the fact that $w$ is even in $x_1$ to write
\begin{equation*}
w(-\pi,x_2)=w(+\pi,x_2),
\end{equation*}
which yields, by uniqueness for the Cauchy problem,
\begin{equation*}
w(x_1-\pi,x_2)=w(x_1+\pi,x_2).
\end{equation*}
This proves that $w$ is $2\pi$-periodic in $x_1$. Next, directly from the definition~\eqref{defid}, we obtain
that $d$ is $2\pi$-periodic in $x_1$ and that $d$ satisfies
\eqref{propertiesd1}. Moreover, this yields
$$
\X2\left(x_1,x_2+2\pi\periode\right)-\X2(x_1,x_2)=
\int_{0}^{2\pi\periode} w(x_1,x_2)\, dx_2.
$$
Differentiating the right-hand side with respect to $x_1$, and using the identity
$\partial_{x_1}w = -\partial_{x_2}(V_2 w / V_1)$, we obtain
$$
\int_{0}^{2\pi\periode} w(x_1,x_2)\, dx_2=\int_{0}^{2\pi\periode} w(0,x_2)\, dx_2=2\pi\periode,
$$
which completes the proof of \eqref{propertiesd2}.

\smallbreak

We next prove that
\begin{equation}\label{claimwneq0}
\forall x\in\xT^2, \quad w(x)=\partial_{x_2} \X2 (x)\neq 0.
\end{equation}
Suppose for contradiction that there exists $x\in [0,2\pi)\times [0,2\pi\periode)$
such that $w(x)=0$. Set
$$
\alpha = \inf \{ x_1\in [0,2\pi)\,:\, \exists x_2\in [0,2\pi\periode] \text{ s.t. } w(x_1,x_2)=0\}.
$$
By continuity, there exists $y$ such that $w(\alpha,y)=0$. Since $w(0,x_2)=1$, we have $\alpha>0$.
For $0\le x_1<\alpha$, we compute that $1/w$ satisfies
$$
\left(\partial_{x_1} + \frac{V_2}{V_1}\partial_{x_2}  -  \partial_{x_2}\left(\frac{V_2}{V_1}\right) \right)\frac{1}{w}=0,\qquad
\frac{1}{w}(0,x_2)=1.
$$
Let $0<\delta<1$. By Sobolev embedding, there exists a constant $K$ such that
$$
\sup_{(x_1,x_2)\in [0,\delta\alpha]\times [0,2\pi\periode]} \la \frac{1}{w}(x_1,x_2)\ra \le K
\sup_{x_1\in [0,\delta\alpha]} \lA \frac{1}{w}(x_1,\cdot)\rA_{H^{1}(T)}.
$$
Therefore, classical energy estimates imply that
$$
\sup_{(x_1,x_2)\in [0,\delta\alpha]\times [0,2\pi\periode]} \la \frac{1}{w}(x)\ra \le
(K+KC)e^{4 C \alpha}
$$
with
$$
C\defn \sup_{x_1\in [0,\delta\alpha]} \lA \frac{V_2}{V_1}(x_1,\cdot) \rA_{C^{2}(T)}.
$$
Therefore, if $0\le x_1<\alpha$, then
$$
w(x)> (K+KC)^{-1}e^{-4C\alpha}.
$$
This gets us to $w(\alpha,\cdot)>0$, whence the contradiction which proves~\eqref{claimwneq0}.
Consequently, $\det \chi'(x)\neq 0$ for all $x\in \xR^2$.
The above argument also establishes that there exists a constant $c\ge 1$ such that,
for all $x\in\xR^2$,
$$
\frac{x_2}{c} \le \X2(x) \le c x_2.
$$
This implies that
$\chi_1$ is a diffeomorphism of $\xR^2$, which completes the proof.
\end{proof}

The end of the proof of Proposition~\ref{prop:r1} follows from Section~3 in \cite{IP}. 

Firstly, directly from the identity $V\cdot\partialx \X2=0$, we obtain
$$
V(x)\cdot ({}^t\chi_1'(x)\xi)=V_1(x)\xi_1,
$$
and hence
$$
\sypr(x,{}^t\chi_1'(x)\xi)= \dns^{1}(x,{}^t\chi_1'(x)\xi)-\mathfrak{a}(x)^{-1}(V_1(x)\xi_1)^2.
$$
Our next task is to rewrite $\dns^{1}(x,{}^t\chi_1'(x)\xi)$ in an appropriate form.

\begin{lemm}\label{lemm:prep1}
There holds
\begin{equation}\label{desireddns1}
\dns^{1}(x,{}^t\chi_1'(x)\xi)=\dns^{1}\left(x,{}^t\chi_1'(x)\begin{pmatrix} 0 \\ 1\end{pmatrix}\right)\la\xi\ra
+ i r(x,\xi)\xi_1,
\end{equation}
where 
$r\in \Gamma^{0}_{s-4}(\xT^2)$ is homogeneous of degree $0$ in $\xi$ and such that
$$
r\in \Gamma (o,e).
$$
\end{lemm}
\begin{proof}
This follows from the fact that $\dns^{1}(x,{}^t\chi_1'(x)\xi)$ is the square root of an homogeneous polynomial function of degree $2$ in $\xi$. 
Indeed, write $\dns^{1}(x,{}^t\chi_1'(x)\xi)^2$ under the form 
\begin{equation}\label{poly}
\dns^{1}(x,{}^t\chi_1'(x)\xi)^2 = P_1(x)\xi_1^2+P_{12}(x)\xi_1\xi_2+P_2(x)\xi_2^2. 
\end{equation}
Then
$$
\dns^{1}(x,{}^t\chi_1'(x)\xi)^2 -P_{2}(x)\la\xi\ra^2= \left(P_1(x)-P_2(x)\right)\xi_1^2+P_{12}(x)\xi_1\xi_2. 
$$
Directly from \eqref{poly} applied with $\xi=e_2\defn \left(\begin{smallmatrix} 0 \\ 1\end{smallmatrix}\right)$, we get
\begin{equation}\label{defiP2}
P_2(x)=\dns^{1}\left(x,{}^t\chi_1'(x)e_2\right)^2.
\end{equation}
Since $\dns\ge 0$, using the elementary 
identity $\sqrt{a}-\sqrt{b}=(a-b)/(\sqrt{a}+\sqrt{b})$, we obtain 
the desired result \eqref{desireddns1} with
$$
r(x,\xi) =-i\frac{\left(P_1(x)-P_2(x)\right)\xi_1+P_{12}(x)\xi_2}{\dns^{1}(x,{}^t\chi_1'(x)\xi)+\sqrt{P_2(x)}\la\xi\ra}.
$$

Since $d\in C^{s-3}(\xR^2)$, we verify that $r$ belongs to $\Gamma^{0}_{s-4}(\xT^2)$. 
It remains to show that $r(x,\xi)\in \Gamma(o,e)$. To fix idea we study the symmetry with respect to the change of variables
$(x,\xi)\mapsto (x^\star,\xi^\star)$ (see Notation~\ref{notaetoile}). We want to prove that 
$r(x^\star,\xi^\star)=-r(x,\xi)$. To see this, $r$ is better written as
$$
r(x,\xi)=-i\frac{1}{\xi_1}\frac{(\dns^{1}(x,{}^t\chi_1'(x)\xi))^2 -P_2(x)\la \xi\ra^2}{\dns^{1}(x,{}^t\chi_1'(x)\xi)+\sqrt{P_2(x)}\la\xi\ra}. 
$$
Note that \eqref{propertiesd1} implies that  
${}^{t}\chi'_1(x^\star)\xi^\star=\bigl[{}^{t}\chi'_1(x)\xi\bigr]^{\star}$ which in turn implies that 
$\dns^{1}(x^\star,{}^t\chi_1'(x^\star)\xi^\star)=\dns^{1}(x,{}^t\chi_1'(x)\xi)$. Moreover, directly from the definition of $P_2$ we see that 
$P_2(x^\star)=P_2(x)$. This yields the desired result $r(x^\star,\xi^\star)=-r(x,\xi)$.
\end{proof}

Introduce next the symbol $p_1$ given by
\begin{equation*}
\sypr_1(\chi_1(x),\xi)\defn \sypr(x,{}^t \chi_1'(x)\xi).
\end{equation*}

\begin{lemm}
There exists a constant $\nu\geq 0$, a $C^{s-4}(\xT^2)$-diffeomorphism $\chi_2$,
a positive function $M\in C^{s-4}(\xT^2)$
and a symbol
$\widetilde{\alpha}\in \Gamma^{0}_{s-6}(\xT^2)$ homogeneous of degree $0$ in $\xi$
such that, for all $(x,\xi)\in \xT^2\times\xR^2$,
\begin{equation*}
\sypr_{1}\bigl(x,{}^t\chi_{2}'(x) \xi\bigr)=
M(x)\left(\la \xi\ra - \nu \xi_1^2\right)+i \widetilde{\alpha}(x,\xi)\xi_1,
\end{equation*}
and such that the following properties hold: $M\in C(e,e)$ and $\widetilde{\alpha}\in \Gamma (o,e)$, and 
$\chi_2$ is of the form
\begin{equation}\label{chi2}
\chi_2(x_1,x_2)=\begin{pmatrix} x_1 + \widetilde{d}(x_1,x_2)\\ x_2+\widetilde{e}(x_2)\end{pmatrix},
\end{equation}
where $\widetilde{d}\in C^{s-4}(\xT^2)$ is odd in $x_1$ and even in $x_2$,
and $\widetilde{e}\in C^{s-3}(\xR/2\pi\periode\xZ)$ is odd in $x_2$.
\end{lemm}
\begin{proof}If $\chi_2$ is of the form~\eqref{chi2}, then
$$
^t\chi_2'(x)\begin{pmatrix}\xi_1 \\ \xi_2\end{pmatrix}=
\begin{pmatrix} (1 + \partial_{x_1}\widetilde{d}(x))\xi_1  \\
\partial_{x_2}\widetilde{d} (x)\xi_1 + (1+\partial_{x_2}\widetilde{e}(x_2))\xi_2\end{pmatrix}.
$$
By transforming the symbols as in the proof of Lemma~\ref{lemm:prep1},
we find that it is sufficient to find $\nu\geq 0$,
$\widetilde{d}=\widetilde{d}(x_1,x_2)$ and $\widetilde{e}=\widetilde{e}(x_2)$ such that
$$
(1 + \partial_{x_1}\widetilde{d}(x_1,x_2))^{2}=\nu \Gamma(x)(1+\partial_{x_2}\widetilde{e}(x_2)),
\quad \text{ where } \Gamma(x)\defn \frac{P_2 \mathfrak{a}}{V_1^{2}}(\chi_1^{-1}(x)).
$$
Therefore, we set
$$
\widetilde{d}(x_1,x_2)=\int_{0}^{x_1}\sqrt{\nu  \Gamma(t,x_2) (1+\partial_{x_2}\widetilde{e}(x_2))}\, dt -x_{1}.
$$
Then, $\widetilde{d}$ is $2\pi$-periodic in $x_1$ if and only if,
$$
\forall x_2\in [0,2\pi\periode],\qquad
\sqrt{ \nu (1+\partial_{x_2}\widetilde{e}(x_2))}\int_{0}^{2\pi} \sqrt{\Gamma (x_1,x_2)}\, dx_1-2\pi = 0.
$$
This yields an equation for $\widetilde{e}=\widetilde{e}(x_2)$ which has a
$(2\pi\periode)$-periodic
solution if and only if $\nu$ is given by
$$
\nu\defn \frac{2\pi}{\periode}\int_{0}^{2\pi\periode}\left(\int_{0}^{2\pi}
\sqrt{\Gamma(x_1,x_2)} \, dx_1\right)^{-2}\, dx_2.
$$

Now, notice that there exists a positive constant $C$ such that 
$P_2\ge C$ (by definition of $P_2$, cf \eqref{defiP2}), $\mathfrak{a}\ge C$ (by Proposition~\ref{theo:WU}), 
$V_1^2\ge C$ (by the assumption~\eqref{cond:x1}). Hence, with $\nu$, $\widetilde{d}$ and $\widetilde{e}$ as previously determined, 
we easily check that $\chi_2$ is a $C^{s-4}(\xT^2)$-diffeomorphism.
\end{proof}

To complete the proof of Proposition~\ref{prop:r1}, set
$$
\chi(x)=\chi_2(\chi_1(x)),
$$
to obtain
\begin{align*}
\sypr\bigl(x,{}^t\chi'(x) \xi\bigr)&=\sypr\bigl(x,{}^t\chi_{1}'(x){}^t \chi_2'(\chi_1(x)) \xi\bigr)\\
&=\sypr_{1}\bigl(\chi_1(x),{}^t\chi_{2}'(\chi_1(x)) \xi\bigr)\\
&=M(\chi_1(x))\left(\la \xi\ra - \nu \xi_1^2\right)+i \widetilde{\alpha}(\chi_1(x),\xi)\xi_1,
\end{align*}
so that we obtain the desired result with
$$
\gamma(x)= M(\chi_1(x)),\quad \alpha(x,\xi)=\widetilde{\alpha}(\chi_1(x),\xi).
$$

\subsection{Paracomposition}\label{paracomposition}
To compute the effect of the
change of variables $x\mapsto \chi(x)$,
we shall use Alinhac's paracomposition operators.
We refer to~\cite{Alipara,Taylor3} for general theory about Alinhac's paracomposition operators.
We here briefly state the basic definitions and results
for periodic functions. Roughly speaking, these results assert that, given $r>1$,
one can associate to a $C^{r}(\xT^2)$-diffeomorphism $\chi$
an operator $\chi^*$ of order $\le 0$ such that, on the one hand,
$$
\forall \alpha\in\xR,\quad u\in H^{\alpha}(\xT^2)\Leftrightarrow \chi^*u\in H^\alpha(\xT^2),
$$
and on the other hand, there is a symbolic calculus to compute the commutator of $\chi^*$ to a \paradif operator.

Let $\phi\colon \xR\rightarrow \xR$ be a smooth even function with
$\phi(t)=1$ for $\la t\ra\le 1.1$ and $\phi(t)=0$ for $\la t\ra\geq 1.9$.
For $k\in \xZ$, we introduce the symbol
$$
\phi_k(\xi)=\phi \left(2^{-k}(1+\la\xi\ra^2)^{1/2}\right),
$$
and then the operator
$$
\widehat{\Delta_{k} f}(\xi)\defn \left(\phi_{k}(\xi) -\phi_{k-1}(\xi)\right)\widehat{f}(\xi).
$$
For all temperate distribution $f\in\mathcal{S}'(\xR^d)$, the spectrum of $\Delta_k f$ satisfies
${\rm spec}\, \Delta_k f \subset \{ \xi\,:\, 2^{k-1} <\L{\xi}< 2^{k+1}\}$.
Hence $\Delta_k f=0$ when $k<0$. Thus, one has
the Littlewood--Paley decomposition:
$$
f=\sum_{k\geq 0}\Delta_k f.
$$
\begin{defi}
Let
$\chi$ be a $C^{r}(\xT^2)$ diffeomorphism with $r>1$. By definition
$$
\chi^* f = \sum_{j\in\xN\,} \sum_{\,\la k -j\ra \le N} \Delta_k \left( (\Delta_j f) \circ \chi )\right),
$$
where $N$ is large enough (depending on $\lA \tilde{\chi}\rA_{C^1}$ only, where $\tilde{\chi}(x)=\chi(x)-x$).
\end{defi}

Two of the principal facts about paracomposition operators
are the following theorems, whose proofs follow from~\cite{Alipara} by adapting the analysis to the case
of $C^{r}(\xT^2)$-diffeomorphisms. The first basic result is that
$\chi^*$ is an operator of order $\le 0$ which can be inverted in the following sense.

\begin{theo}\label{theo:paracomp1}
Let $\chi$ be a $C^{r}$-diffeomorphism with $r>1$. 
For all $\alpha\geq 0$, $f\in H^{\alpha}(\xT^2)$ if and only if $\chi^*f\in H^{\alpha}(\xT^2)$. 
Moreover, $\chi^* (\chi^{-1})^* - I $ is of order $\le -(r-1)$.
\end{theo}

This theorem reduces the study of the regularity of $u$ to the study of the regularity of
$\chi^* u$. To study the regularity of $\chi^* u$ we need to compute the equation satisfied by $\chi^* u$.
To do this, we shall use a symbolic calculus theorem which
allows to compute the equation satisfied by~$\chi^* u$ in terms of the equation satisfied by $u$
(in analogy with the \paradif calculus). For what follows, it is convenient to work with $(\chi^{-1})^*$.

\begin{theo}\label{theo:paracomp2}
Let $m\in\xR$, $r>1$, $\rho>0$ and set $\sigma\defn \inf\{\rho,r-1\}$. 
Consider a $C^{r}(\xT^2)$-diffeomorphism $\chi$ 
and a symbol $a\in \Sigma^{m}_{\rho}(\xT^2)$,
then there exists $a^*\in \Sigma^{m}_{\sigma}(\xT^2)$ such that
$$
(\chi^{-1})^* T_{a}  - T_{a^*}(\chi^{-1})^* \quad\text{is order }\le m- \sigma.
$$
Moreover, one can give an explicit formula for $a^*$:
If one decomposes $a$ as a sum of homogeneous symbols, then
\begin{equation}\label{defi*}
a^*(\chi(x),\eta)=\sum  
\frac{1}{i^\alpha \alpha !}
\partial_{\xi}^{\alpha} a_{m-k}(x,{}^t\chi'(x) \eta) \partial_{y}^\alpha (e^{i\Psi_{x}(y)\cdot\eta})\arrowvert
_{y=x},
\end{equation}
where the sum is taken over all $\alpha\in \xN^2$ such that the summand is well defined, 
$\chi'(x)$ is the differential of $\chi$, $t$ denotes transpose and
\begin{equation}\label{defiPsi}
\Psi_{x}(y)=\chi(y)-\chi(x)-\chi'(x)(y-x).
\end{equation}
\end{theo}

\subsection{The first reduction}
We here apply the previous results to perform a first reduction to a
case where the ``principal'' part of the equation has constant coefficient.

\begin{prop}\label{C514}
Let $\chi$ be as given by Proposition~$\ref{prop:r1}$. Then
$\bigl(\chi^{-1}\bigr)^{*} u$ satisfies an equation of the form
\begin{equation}\label{eq:inter2}
\Bigl(\la D_x \ra  +\nu\partial_{x_1}^2 + T_A \partial_{x_1}  + T_B \Bigr) \bigl(\chi^{-1}\bigr)^{*} u =f\in H^{s+2}(\xT^2),
\end{equation}
where $A\in \Gamma^{0}_{s-6}(\xT^2)$ and $B\in \Sigma^0_{s-6}(\xT^2)$ are such that:
\begin{enumerate}[i)]
\item $A$ is homogeneous of degree $0$ in $\xi$; $B=B_0+B_{-1}$ where $B_{\ell}$ is homogeneous of degree $\ell$ in $\xi$;
\item $A\in \Gamma (o,e)$ and $B_{\ell}\in \Gamma(e,e)$ for $\ell=0,1$.
\end{enumerate}
\end{prop}
\begin{rema}
$(i)$ For what follows it is enough 
to have remainders in $H^{s+2}(\xT^2)$.
From now on, to simplify the presentation we do not try to give results with remainders
having a regularity higher than what is needed.

$(ii)$ $A\in \Gamma (o,e)$ and $B_{\ell}\in \Gamma(e,e)$ is equivalent to the fact that 
the symbol of $\la D_x \ra  +\nu\partial_{x_1}^2 + T_A \partial_{x_1}  + T_B$ is invariant by the symmetries 
\begin{align*}
&(x_1,x_2,\xi_1,\xi_2)\mapsto (-x_1,x_2,-\xi_1,\xi_2),\\
&(x_1,x_2,\xi_1,\xi_2)\mapsto (x_1,-x_2,\xi_1,-\xi_2).
\end{align*}
\end{rema}
\begin{proof}
We begin by applying the results in \S\ref{paracomposition} to compute the equation satisfies by $\bigl(\chi^{-1}\bigr)^{*} u$.
Recall that, by notation, $\mathcal{V}$ is as given by \eqref{defimV} and
$$
p(x,\xi)
=\sqrt{(1+\la\partialx\sigma\ra^2)\la\xi\ra^2-(\partialx\sigma\cdot\xi)^2}
-\mathfrak{a}(x)^{-1}(V(x)\cdot\xi)^2.
$$
We define $\dns^*$ by \eqref{defi*} applied with $m=1$ and $\rho=s-1$.
Similarly, we define $\mathcal{V}^*$ and $p^*$ by \eqref{defi*} applied with $m=2$ and $\rho=s-4$.
To prove Proposition~\ref{C514}, the key step is to compare the principal symbol of $p^*$
with $\dns^*+\mathcal{V}^*$.

\begin{lemm}\label{prop:subsub}
There exist $a\in \Gamma^{0}_{s-6}(\xT^2)$ and $b=b_{0}+b_{-1}\in \Sigma^0_{s-6}(\xT^2)$ and a remainder 
$r\in \Sigma^{-2}_{s-8}(\xT^2)$ such that
$$
\dns^*(\chi(x),\xi)+\mathcal{V}^*(\chi(x),\xi)=p(x,{}^t\chi'(x)\xi)+i a(x,\xi)\xi_1
+b(x,\xi)+r(x,\xi),
$$
and such that $a\in \Gamma(o,e)$, $b\in \Gamma(e,e)$. 
\end{lemm}

\begin{proof} The proof is elementary.
We first study $\dns^*(\chi(x),\xi)$.
Since $\dns$ is a symbol of order $1$, to obtain a remainder $r$ which is of order $-2$,
we need to compute the first three terms in the symbolic expansion of~$\dns^*$.
To do this, note that there are some cancelations which follow directly from the
definition~\eqref{defiPsi}:
\begin{align*}
&\la\alpha\ra=1 \quad\Rightarrow \partial_{y}^\alpha (e^{i\Psi_{x}(y)\cdot\xi})\arrowvert
_{y=x} = 0,\\
&2\le \la\alpha\ra\le 3 \quad\Rightarrow \partial_{y}^\alpha (e^{i\Psi_{x}(y)\cdot\xi})\arrowvert
_{y=x} = i \partial_{x}^\alpha ( \chi(x)\cdot\xi),\\
&\la\alpha\ra=4 \quad\Rightarrow \partial_{y}^\alpha (e^{i\Psi_{x}(y)\cdot\xi})\arrowvert
_{y=x} = i \partial_{x}^\alpha (\chi(x)\cdot\xi) - \sum  \partial_{x}^\beta (\chi(x)\cdot\xi)
\partial_{x}^\gamma (\chi(x)\cdot\xi),
\end{align*}
where in the last line the sum is taken over all decompositions $\beta+\gamma=\alpha$ such that
$\la\beta\ra=2=\la\gamma\ra$.
Therefore, it follows from~\eqref{defi*} that
\begin{equation}
\dns^*(\chi(x),\xi) = \dns^1(x,{}^t\chi'(x) \xi) +b_{0}(x,\xi)+b_{-1}(x,\xi)+r(x,\xi),
\end{equation}%
where $r\in \Sigma^{-2}_{s-8}(\xT^2)$ and
\begin{align*}
b_{0}(x,\xi) &\defn \dns^0(x,{}^t\chi'(x) \xi)-i\sum_{\la\alpha\ra=2}  \frac{1}{\alpha !}
(\partial_{\xi}^{\alpha} \dns^1)(x,{}^t\chi'(x) \xi) \partial_{x}^\alpha \chi(x)\cdot\xi,\\
b_{-1}(x,\xi) &\defn  \dns^{-1}(x,{}^t\chi'(x) \xi)
+\sum_{\ell=0}^{1}\sum_{\la\alpha\ra=2+\ell}  \frac{1}{i^{\alpha}\alpha !}
(\partial_{\xi}^{\alpha} \dns^\ell)(x,{}^t\chi'(x) \xi) \partial_{x}^\alpha \chi(x)\cdot\xi\\
&\quad - \sum_{\la\beta\ra=2=\la\gamma\ra} (\partial_{\xi}^{\beta+\gamma} \dns^1)(x,{}^t\chi'(x) \xi)
(\partial_{x}^\beta \chi(x)\cdot\xi)
 (\partial_{x}^\gamma \chi(x)\cdot\xi).
\end{align*}

Recall that $\chi=(\chi^1,\chi^2)$ where $\chi^1$ is odd in $x_1$ and even in $x_2$, and $\chi^2$ is even in $x_1$ and
odd in $x_2$. Therefore, to prove the desired symmetry properties, it is sufficient to prove that $\dns^1,\dns^0,\dns^{-1}$ are invariant
by the two changes of variables
\begin{align*}
&(x_1,x_2,\xi_1,\xi_2)\mapsto (-x_1,x_2,-\xi_1,\xi_2),\\
&(x_1,x_2,\xi_1,\xi_2)\mapsto (x_1,-x_2,\xi_1,-\xi_2).
\end{align*}
We consider the first case only and use the notation
$$
f^\star(x_1,x_2)=f(-x_1,x_2).
$$
Observe that, since $\sigma^\star=\sigma$, it follows directly from the definition  of the Dirichlet to Neumann operator
(see \eqref{defi:dn}) that
$$
G(\sigma)f ^\star = \bigl[ G(\sigma) f\bigr]^\star.
$$
On the symbol level, this immediately gives the desired result:
$$
\dns(x^\star,\xi^\star)=\dns(x,\xi).
$$
Alternatively, one may use the
explicit definition of the symbols $\Slam_m$ given in the proof of Lemma~\ref{lemm:total}.

\smallbreak
Since, by notation,
\begin{equation*}
\mathcal{V}(x,\xi) \defn - \mathfrak{a}(x)^{-1}(V(x)\cdot \xi)^2 +i \cnx \left(\mathfrak{a}(x)^{-1}(V(x)\cdot \xi) V(x)\right),
\end{equation*}
the same reasoning implies that
\begin{align*}
\mathcal{V}^*(\chi(x),\xi)&=
- \mathfrak{a}^{-1}(V(x)\cdot {}^t\chi'(x) \xi)^2 +i \mathfrak{a}(x)^{-1}(V(x)\cdot {}^t\chi'(x) \xi) \cnx V(x)\\
&\quad +i \mathfrak{a}(x)^{-1}\left(V(x)\cdot \partialx V(x)\right)\cdot {}^t\chi'(x) \xi \\
&\quad +i\mathfrak{a}(x)^{-1}\sum_{1\le k,\ell \le 2}
V_k (x) V_\ell (x) \partial_{x_k}\partial_{x_\ell} (\chi(x)\cdot\xi).
\end{align*}
Gathering the last two terms in the right-hand side we thus obtain
\begin{align*}
\mathcal{V}^*(\chi(x),\xi)&=
- \mathfrak{a}^{-1}(V(x)\cdot {}^t\chi'(x)\xi)^2 +i \mathfrak{a}(x)^{-1}(V(x)\cdot {}^t\chi'(x) \xi) \cnx V(x)\\
&\quad 
+i\mathfrak{a}(x)^{-1}V (x)\cdot \partialx \left(V(x)\cdot {}^t\chi'(x) \xi\right)
\end{align*}
and hence, by construction of $\chi$, it follows that
$$
\mathcal{V}^*(\chi(x),\xi)=- \mathfrak{a}(x)^{-1}(V(x)\cdot {}^t\chi'(x)\xi)^2+ i a(x,\xi)\xi_1,
$$
where $a\in\Gamma^{0}_{s-5}(\xT^2)$ is such that 
$a\in \Gamma(o,e)$. 
This completes the proof of Lemma~\ref{prop:subsub}.
\end{proof}

We now are now in position to prove Proposition~\ref{C514}.
By using Theorem~\ref{theo:paracomp2}, it follows from Lemma~\ref{prop:subsub} and
Proposition~\ref{prop:r1}  that
$$
\Bigl(T_{\gamma} \bigl(\la D_x \ra  +\nu\partial_{x_1}^2\bigr) + T_{\alpha+a} \partial_{x_1}
+ T_b \Bigr) \bigl(\chi^{-1}\bigr)^{*} u \in H^{s+2}(\xT^2).
$$
We thus obtain~\eqref{eq:inter2} with
$A= (\alpha+a)/\gamma$ and $B= b/\gamma$.
\end{proof}

\subsection{Elliptic regularity far from the characteristic variety}

As usual, the analysis makes use of the division of the phase space into a region in which the inverse
of the symbol remains bounded by a fixed constant
and a region where the symbol is small and may vanish.
Here we consider the easy part and prove the following elliptic regularity result.

\begin{prop}\label{prop:ellipfar}
Let $\chi$ be the diffeomorphism determined
in Proposition~$\ref{prop:r1}$. Consider $\Theta=\Theta(\xi)$ homogenous of degree $0$ and such that
there exists a constant $K$ such that
$$
\la\xi_2\ra\geq K \la\xi_1\ra \Rightarrow \Theta(\xi_1,\xi_2)=0.
$$
Then,
$$
\Theta(D_x)( \bigl(\chi^{-1}\bigr)^{*} u) \in H^{s+2}(\xT^2).
$$
\end{prop}
\begin{rema}
Note that, on the characteristic variety, we have $\la\xi_2\ra\sim\nu\xi_1^2$.
On the other hand, on the spectrum of $\Theta(D_x)(\bigl(\chi^{-1}\bigr)^{*} u)$, we have $\la\xi_2\ra \le K\la\xi_1\ra$.
Therefore, the previous result establishes elliptic regularity very far from the characteristic variety.
\end{rema}
\begin{proof}Recall that $\bigl(\chi^{-1}\bigr)^{*} u$ satisfies \eqref{eq:inter2}. Set
$$
\wp(x,\xi)\defn\la\xi\ra-\nu\xi_1^2 + i A(x,\xi)\xi_1 +B(x,\xi).
$$
We have
\begin{equation}\label{ffcv1}
T_{\wp}(\bigl(\chi^{-1}\bigr)^{*} u) \in H^{s+2}(\xT^2).
\end{equation}
Note that
$\la\wp(x,\xi)\ra\geq c \la \xi\ra^2$ for some constant $c>0$ for all $(x,\xi)$ such that
$$
\la\xi_2\ra \le K\la\xi_1\ra \quad\text{and}\quad \la\xi\ra\geq M,
$$
for some large enough constant $M$ depending on $\sup_{x,\xi} \la A(x,\xi)\ra+\la B(x,\xi)\ra$.
Introduce a $C^\infty$ function $\tilde\Theta$ such that
\begin{equation*}
\begin{split}
&\tilde\Theta(\xi)=0
\quad\text{for}\quad \la\xi\ra\le M,\\
&\tilde\Theta(\xi)=\Theta(\xi) \quad\text{for }
\la\xi\ra\geq 2M.
\end{split}
\end{equation*}
Since $\Theta$ and $\tilde\Theta$ differ only on a bounded neighborhood of the origin, we have
$$
\Theta(D_x)(\bigl(\chi^{-1}\bigr)^{*} u)-\tilde\Theta(D_x)( \bigl(\chi^{-1}\bigr)^{*} u) \in C^{\infty}(\xT^2).
$$
Note that, since $\Theta$ is positively homogeneous of degree $0$,
$\tilde{\Theta}$ belongs to our symbol class ($\tilde{\Theta}\in \Gamma_{\rho}^{0}(\xT^2)$ for all $\rho\geq 0$). 
Set
$$
q=\frac{\tilde{\Theta}}{\wp}-\frac{1}{i}\partial_{\xi}\left( \frac{\tilde{\Theta}}{\wp}\right) \frac{\partial_x \wp}{\wp}\in
\Gamma^{-2}_{s-8}(\xT^2).
$$
According to Theorem~\ref{theo:sc} applied with $(m,m',\rho)=(2,-2,2)$, then,
$$
T_{q}T_{\wp} (\bigl(\chi^{-1}\bigr)^{*} u)- \tilde\Theta(D_x)(\bigl(\chi^{-1}\bigr)^{*} u)\in H^{s+2}(\xT^2).
$$
On the other hand, since $T_{q}$ is of order $\le 0$,
it follows from~\eqref{ffcv1} that
$$
T_{q}T_{\wp}(\bigl(\chi^{-1}\bigr)^{*} u)\in H^{s+2}(\xT^2).
$$
Which completes the proof.
\end{proof}

\subsection{The second reduction and the end of the proof of Theorem~\ref{theo:smooth}}
\label{srend}
We first set a few notations.
Introduce a $C^\infty$ function $\eta$ such that $0\le\eta\le 1$,
\begin{equation}\label{defieta}
\eta(t)=0\quad \text{for } t\in [-1/2,1/2],\qquad
\eta(t)=1\quad \text{for }\la t\ra\geq 1.
\end{equation}
Given $k\in\xN$, $\partial_{x_1}^{-k}$ denotes the Fourier multiplier
(defined on $\mathcal{S}'(\xR^d)$) with symbol
$\eta(\xi_1)(i\xi_1)^{-k}$. Note that, if $f$ is $2\pi$-periodic in $x_1$, then
$$
\partial_{x_1}^{0} f = f -\frac{1}{2\pi}\int_{0}^{2\pi} f(x_1,x_2)\, dx_2,
$$
and
$$
(\partial_{x_1}^{-1} f )(x_1,x_2)=\int_{0}^{x_1} \left( f (s,x_2)
-\frac{1}{2\pi}\int_{0}^{2\pi} f(x_1,x_2)\,dx_1 \right)\, ds.
$$
In particular, $\partial_{x_1}\partial_{x_1}^{-1}f=\partial_{x_1}^{0}f=f$ if and only if $f$ has zero mean value in~$x_1$.
We also have
$$
\partial_{x_1}^{-k-1} f = \partial_{x_1}^{-k}\partial_{x_1}^{-1}f.
$$

It will be
convenient to divide the frequency space into three pieces so that, in the two main parts,
$\xi_2$ is either positive or negative.
To do this, we need to use Fourier multipliers whose symbols belong to our symbol class,
which is necessary to apply symbolic calculus
in the forthcoming computations.
Here is one way to define such Fourier multipliers:
consider a $C^\infty$ function $J$ satisfying $0\le J\le 1$ and such that,
\begin{equation}\label{deltacomes}
J(s)=0 \quad\text{for } s\le 0.8,\qquad J(s)=1 \quad\text{for } s \geq 0.9,
\end{equation}
We define three $C^\infty$ functions $\jmath_0,\jmath_{-}$ and $\jmath_{+}$ by
\begin{equation*}
\jmath_{0}=1-\jmath_{-}-\jmath_{+},\quad
\jmath_{-}(\xi)=J\left(\frac{\la\xi\ra- \xi_2}{2\la\xi\ra}\right),
\quad\jmath_{+}(\xi)=J\left(\frac{\la\xi\ra+ \xi_2}{2\la\xi\ra}\right),
\end{equation*}
and then the Fourier multipliers
$$
\widehat{\jmath_{\eps}(D_x) f}(\xi)= \jmath_{\eps}(\xi)\widehat{f}(\xi) \qquad (\eps\in \{0,-,+\}).
$$
Note that there are constants $0<c_1<c_2$ such that
\begin{alignat}{6}
&\xi_2\le c_1\la\xi_1\ra  &\Rightarrow j_{+}(\xi)&=0, \qquad
&\xi_2\geq  c_2 \la\xi_1\ra &\Rightarrow j_{+}(\xi)&=1,\label{po1}\\
&\xi_2\geq  -c_1 \la\xi_1\ra  &\Rightarrow j_{-}(\xi)&=0,\qquad
&\xi_2\le -c_2\la\xi_1\ra  &\Rightarrow j_{-}(\xi)&=1.\label{po2}
\end{alignat}
Also, note that $\jmath_\pm$ is positively homogeneous of degree $0$ and hence
satisfies
$$
\la \partial_{\xi}^\alpha \jmath_\pm (\xi)\ra \le C_\alpha \la\xi\ra^{-\la\alpha\ra}.
$$
In view of \eqref{po1} and \eqref{po2}, Proposition~\ref{prop:ellipfar} implies that
\begin{equation}\label{j0}
\jmath_{0}(D_x)(\bigl(\chi^{-1}\bigr)^{*} u)\in H^{s+2}(\xT^2).
\end{equation}
As a result, it remains only to concentrate on the two other terms:
$$
\jmath_{\pm}(D_x)\left(\bigl(\chi^{-1}\bigr)^{*} u\right).
$$
Here is one other obvious observation that enable us to reduce
the analysis to the study of only one of these two terms:
Since $\bigl(\chi^{-1}\bigr)^{*} u$ is even in $x_2$, we have
\begin{equation*}
\widehat{\bigl(\chi^{-1}\bigr)^{*} u}(\xi_1,\xi_2)=\widehat{\bigl(\chi^{-1}\bigr)^{*} u}(\xi_1,-\xi_2),
\end{equation*}
Therefore,
\begin{equation}\label{jpm}
\jmath_{-}(D_x)\left(\bigl(\chi^{-1}\bigr)^{*} u\right)\text{ and }
 \jmath_{+}(D_x)\left(\bigl(\chi^{-1}\bigr)^{*} u\right) \text{ have the same regularity.}
\end{equation}
Consequently, it suffices to study one of these two terms. We chose to work with
$$
U\defn \jmath_{+}(D_x)\left(\bigl(\chi^{-1}\bigr)^{*} u\right).
$$

We shall prove that one can transform further the problem
to a linear equation with {\em constant} coefficients, using the method
of Iooss and Plotnikov~\cite{IP}.
The key to proving Theorem~\ref{toprove} is the following.

\begin{prop}\label{prop440}
There exist two constants $\kappa,\kappa'\in\xR$ and an operator
$$
Z_c = \sum_{0\le j\le 4} T_{c_j}\partial_{x_1}^{-j},
$$
where $c_{0},\ldots,c_{4}\in C^{1}(\xT^2)^{5}$ and $\la c_0\ra>0$, such that
\begin{equation}\label{canonique0}
\left( -i\partial_{x_2} + \nu\partial_{x_1}^2 +\kappa + \kappa' \partial_{x_1}^{-2} \right)
Z_{c} U\in H^{s+2}(\xT^2).
\end{equation}
\end{prop}
\begin{proof}
Proposition \ref{prop440} is proved in \S\ref{preparation} and \S\ref{psr}.
\end{proof}

We here explain how to conclude the proof of Theorem~\ref{toprove}.

\begin{proof}[Proof of Theorem~\ref{toprove} given Proposition~\ref{prop440}]
Since the symbol of the operator
$-i\partial_{x_2} + \nu\partial_{x_1}^2 +\kappa + \kappa' \partial_{x_1}^{-2}$
is $\xi_2 - \nu\xi_1^2 +\kappa - \kappa' \xi_1^{-2}$,
we set
$$
\nu(\var)=\periode \nu,\quad
\kappa_0(\var)=- \periode\kappa,
\quad
\kappa_1(\var)=\periode\kappa'.
$$
Assume that there exist $\delta\in [0,1[$ and $N\in \xN^*$ such that,
\begin{equation*}
\la \k_2 - \nu(\var) \k_1^2 -\kappa_0(\var)-\frac{\kappa_1(\var)}{\k_1^2}\ra\geq
\frac{1}{ \k_1^{2+\delta}},
\end{equation*}
for all $\k \in \xN^2$ with $\k_1$ sufficiently large.

Directly from the definitions of the coefficients, this assumption implies that
\begin{equation*}
\la \xi_2 - \nu \xi_1^2 +\kappa-\frac{\kappa'}{\xi_1^2}\ra\geq
\frac{1}{ \periode \xi_1^{2+\delta}},
\end{equation*}
for all $\xi=(\xi_1,\xi_2)\in \xN\times(\periode^{-1}\xN)$ with $\xi_1$ sufficiently large.
Since $\nu\ge 0$, the previous inequality holds for all
$\xi\in \xZ\times \periode^{-1}\xZ$ with $\la\xi_1\ra$ sufficiently large.

Now, since $\la\xi\ra\sim \nu\xi_1^2$ on the set where the above inequality
is not satisfied, this in turn implies that,
\begin{equation}\label{diophantine5}
\la \xi_2 - \nu \xi_1^2 +\kappa -\frac{ \kappa'}{\xi_1^2}\ra\geq
\frac{\nu}{\periode \la\xi\ra^{(2+\delta)/2}},
\end{equation}
for all $\xi\in \xZ\times(\periode^{-1}\xZ)$ with $\la\xi_1\ra$ sufficiently large.

Similarly, we obtain that
\begin{equation}\label{diophantine6}
\la \xi_2- \nu \xi_1^2 +\kappa - \frac{ \kappa'}{\xi_1^2}\ra\geq
 \frac{\sqrt{\nu}\la\xi_1\ra}{\periode\la\xi\ra^{(3+\delta)/2}},
\end{equation}
for $\la\xi_1\ra$ sufficiently large.

To use these inequalities, we take the Fourier transform of \eqref{canonique0}:
$$
\Bigl(\xi_2 - \nu \xi_1^2 +\kappa-\frac{\kappa' \eta(\xi_1)}{\xi_1^2}  \Bigr)
\widehat{Z_{c} U}(\xi) =:
\widehat{f}(\xi).
$$
A key point is that $Z_{c} U$ is doubly periodic.
Thus, if
$\xi$ belongs to the support of the Fourier transform of $Z_{c} U$, then
$\xi \in \xZ\times(\periode^{-1}\xZ)$. Therefore,
it follows from \eqref{diophantine5} and \eqref{diophantine6} that,
$$
Z_{c} U \in H^{s+1-\frac{\delta}{2}}(\xT^2)
 ,\quad \partial_{x_1} Z_{c} U
\in H^{s+\mez-\frac{\delta}{2}}(\xT^2).
$$
It follows that
$$
U \in H^{s+1-\frac{\delta}{2}}(\xT^2)
,\quad \partial_{x_1} U \in H^{s+\frac{1}{2}-\frac{\delta}{2}}(\xT^2).
$$
In view of \eqref{j0} and \eqref{jpm}, we end up with
\begin{align}
&\bigl(\chi^{-1}\bigr)^{*} u \in H^{s+1-\frac{\delta}{2}}(\xT^2), \label{u1}\\
&\partial_{x_1} (\bigl(\chi^{-1}\bigr)^{*} u)\in H^{s+\mez-\frac{\delta}{2}}(\xT^2).\label{u2}
\end{align}
Direclty from \eqref{u1}, by using Theorem~\ref{theo:paracomp1}, we obtain
\begin{equation}\label{preresult}
u\in H^{s+1-\frac{\delta}{2}}(\xT^2).
\end{equation}
We next claim that
\begin{equation}\label{claim5.8}
T_{V}\cdot\partialx  u \in H^{s+\mez-\frac{\delta}{2}}(\xT^2).
\end{equation}
To prove \eqref{claim5.8}, we first note that Theorem~\ref{theo:paracomp1} 
implies that it is enough to prove that
$$
\bigl(\chi^{-1}\bigr)^{*}\left(T_{V}\cdot\partialx  u\right) \in H^{s+\mez-\frac{\delta}{2}}(\xT^2).
$$
To prove this result, we shall apply our second theorem about paracomposition operators. Indeed, 
define the symbol $a$ by
$$
a^{*}(\chi(x),\eta)=V(x)\cdot \left( i {}^{t}\chi'(x)\eta\right).
$$
Theorem~\ref{theo:paracomp2} implies that
$$
\bigl(\chi^{-1}\bigr)^{*}T_{V}\cdot\partialx - T_{a^{*}}\bigl(\chi^{-1}\bigr)^{*}\quad\text{is of order }\le 0. 
$$
Consequently, in view of \eqref{preresult}, we obtain
\begin{equation}\label{preclaim5.8}
\bigl(\chi^{-1}\bigr)^{*}T_{V}\cdot\partialx u - T_{a^{*}}\bigl(\chi^{-1}\bigr)^{*} u \in H^{s+1-\frac{\delta}{2}}(\xT^2). 
\end{equation}
Now, directly from the definition of $\chi$ (see \S\ref{sec:3.1}), we compute that 
$$
a^{*}(\chi(x),\eta)= i \alpha(x) \xi_1 \quad\text{with}\quad \alpha (x)= V_{1}(x) \left( 1 +\partial_{x_1}\tilde{d} (\chi_1(x))\right).
$$
Therefore,  \eqref{u2} implies that
$$ 
T_{a^{*}}\bigl(\chi^{-1}\bigr)^{*} u \in H^{s+\mez-\frac{\delta}{2}}(\xT^2).
$$
Consequently, the claim~\eqref{claim5.8} follows from \eqref{preclaim5.8}

\bigbreak

Now recall that we have proved that (cf \eqref{sansa})
\begin{equation}\label{rpcsystem}
\sigma+T_{\mathfrak{a}^{-1}}T_{V}\cdot\partialx u \in H^{2s-4}(\xT^2).
\end{equation}
Since $T_{\mathfrak{a}^{-1}}$ is of order $\le 0$, by using \eqref{claim5.8} we obtain
$$
\sigma \in H^{2s-4}(\xT^2)+H^{s+\mez-\frac{\delta}{2}}(\xT^2),
$$
and hence
\begin{equation}\label{claim5.8b}
\sigma \in H^{s+\mez-\frac{\delta}{2}}(\xT^2).
\end{equation}
Since $T_\mathfrak{b}$ is of order $\le 0$, this yields 
\begin{equation}\label{preresult2}
T_{\mathfrak{b}}\sigma \in H^{s+\mez-\frac{\delta}{2}}(\xT^2).
\end{equation}
Writing $\psi=u+T_{\mathfrak{b}}\sigma$ (by definition of $u$), from \eqref{preresult} and \eqref{preresult2} we deduce that
$$
\psi\in H^{s+\mez-\frac{\delta}{2}}(\xT^2).
$$
This completes the proof of Theorem~\ref{toprove} and hence of Theorem~\ref{theo:smooth}.
\end{proof}

\subsection{Preparation}\label{preparation}
We have proved that there exist a change of variables $x\mapsto \chi(x)$
and two zero order symbols
$A=A(x,\xi)$ and $B(x,\xi)$ such that
\begin{equation}\label{eq:inter2bis}
\Bigl(\la D_x \ra  +\nu\partial_{x_1}^2 + T_A \partial_{x_1}  + T_B \Bigr) \bigl(\chi^{-1}\bigr)^{*} u =f\in H^{s+2}(\xT^2).
\end{equation}
We want to prove 
Proposition~\ref{prop440} which asserts that it is possible
to conjugate \eqref{eq:inter2bis} to a constant coefficient equation. 
Since the symbols $A$ and $B$
depend on the frequency variable, one more reduction is needed.

In this paragraph we shall prove the following preliminary result towards the proof of
Proposition~\ref{prop440}.
\begin{prop}\label{p:reduced}
There exist five functions
$$
a_j=a_j(x)\in C^{s-6-j}(\xT^2) \qquad (0\le j\le 4),
$$
where
\begin{equation*}
a_{j} \text{ is odd in }x_1
\text{ for }j\in \{0,2,4\},\quad
a_{j} \text{ is even in }x_1 
\text{ for }j\in \{1,3\},
\end{equation*}
such that
\begin{equation*}
\left( -i\partial_{x_2}  + \nu\partial_{x_1}^2 +
\sum_{0\le j\le 4}T_{a_j} \partial_{x_1}^{1-j} \right) U\in H^{s+2}(\xT^2),
\end{equation*}
where recall that $U= j_{+}(D_x)(\bigl(\chi^{-1}\bigr)^{*}u)$.
\end{prop}

To prove this result, we begin with the following localization lemma.
\begin{lemm}\label{li4}
Let $A=A(x,\xi)$ and $B=B(x,\xi)$ be as in \eqref{eq:inter2}.
Then,
\begin{equation}\label{eq:inter4}
\Bigl( \la D_{x}\ra  +\nu\partial_{x_1}^2 + T_{A} \partial_{x_1}  + T_B \Bigr) U\in H^{s+2}(\xT^2).
\end{equation}
\end{lemm}
\begin{proof}
This follows from Corollary~\ref{C514} and
Proposition~\ref{prop:ellipfar}. Indeed,
since $\jmath_{+}$ is positively homogeneous of degree $0$, it is
a zero-order symbol. According to Theorem~\ref{theo:sc} (applied with $a=\jmath_{+}(\xi)$,
$b= A(x,\xi)i\xi_1+B(x,\xi)$ and $(m,m',\rho)=(0,1,3)$), then,
\begin{multline*}
j_+(D_x)\Bigl(\la D_x \ra  +\nu\partial_{x_1}^2 + T_A \partial_{x_1}  + T_B \Bigr)(\bigl(\chi^{-1}\bigr)^{*}u)\\
=\Bigl(\la D_x \ra  +\nu\partial_{x_1}^2 + T_A \partial_{x_1}  + T_B \Bigr) U
+\sum_{1\le \la\alpha\ra\le 2} \frac{1}{i^{\alpha} \alpha!}
T_{\partial_{\xi}^\alpha \jmath_{+}(\xi) \partial_{x}^\alpha b (x,\xi)} (\bigl(\chi^{-1}\bigr)^{*}u)
+f,
\end{multline*}
with $f\in H^{s+2}(\xT^2)$. Corollary~\ref{C514} implies that the left hand side
belongs to $H^{s+2}(\xT^2)$. As regards the second term in the right hand side, observe that
$$
\la \xi_2\ra \geq \frac{3}{4} \la\xi_1\ra \Rightarrow
\partial_{\xi}\jmath_{+}(\xi) =0 .
$$
Moreover $\partial_{\xi}^\alpha \jmath_{+}(\xi) \partial_{x}^\alpha b (x,\xi)$ is of order $\le 0$ in $\xi$. 
Hence, by means of a simple symbolic calculus argument, it
follows from Proposition~\ref{prop:ellipfar} that
\begin{gather*}
T_{\partial_{\xi}^\alpha \jmath_{+}(\xi) \partial_{x}^\alpha b (x,\xi)} (\bigl(\chi^{-1}\bigr)^{*}u) \in H^{s+2}(\xT^2).
\end{gather*}
This proves the lemma.
\end{proof}

We are now in position to transform the equation. To clarify the articulation
of the proof, we proceed step by step. We first
prove that
\begin{enumerate}[i)]
\item
$\la D_x\ra U$ may be replaced by $-i\partial_{x_2}U$ ($-i\partial_{x_2}$ is the Fourier multiplier
with symbol $+\xi_2$). This point essentially follows from the fact that
$\la \xi\ra\sim \xi_2$ on the support of $\jmath_{+}(\xi)$.
\item One may replace the symbols $A$ and $B$ by a couple of symbols which are symmetric with respect to $\{\xi_2=0\}$
and vanish for $\la \xi_2\ra \le  \la \xi_1\ra/5$. The idea is that, since
$\widehat{U}(\xi)=0$ for $\xi_2\le \la\xi_1\ra/2$, for any $c<1/2$ one may freely change the values of $A(x,\xi)$
and $B(x,\xi)$ for $\xi_2\le c\la \xi_1\ra$.
\end{enumerate}

\begin{lemm}\label{ltAB}
There exists two symbols $\tilde{A}\in \Gamma^{0}_{s-6}(\xT^2),\tilde{B}\in \Sigma^{0}_{s-6}(\xT^2)$ such that
\begin{equation}\label{eq:inter2ter}
\Bigl( -i\partial_{x_2}  +\nu\partial_{x_1}^2 + T_{\tilde A} \partial_{x_1}  +
T_{\tilde B} \Bigr) U \in H^{s+2}(\xT^2),
\end{equation}
and such that
\begin{enumerate}[i)]
\item $\tilde A$ is homogeneous of degree $0$ in $\xi$;
$\tilde B=\tilde B_0+ \tilde B_{-1}$ where $\tilde B_{\ell}$ is homogeneous of degree $\ell$ in $\xi$;
\item $\tilde{A}(x^\star,\xi^\star)=-\tilde{A}(x,\xi)$ and
$\tilde{B}(x^\star,\xi^\star)=\tilde{B}(x,\xi)$;
\item $\tilde A (x,\xi_1,-\xi_2)=\tilde A(x,\xi_1,\xi_2)$
and $\tilde B (x,\xi_1,-\xi_2)=\tilde B(x,\xi_1,\xi_2)$;
\item $\tilde A (x,\xi)=0=\tilde B (x,\xi)$ for $\la\xi_2\ra\le \la \xi_1\ra/5$.
\end{enumerate}
\end{lemm}
\begin{proof}
The proof depends on Lemma~\ref{li4} and the fact that the Fourier multiplier $\jmath_{+}(D_x)$
is essentially a projection operator. Namely, we make use of a
 $C^\infty$ function $J'$ satisfying $0\le J'\le 1$ and such that,
\begin{equation*}
J'(s)=0 \quad\text{for } s\le 0.7,\qquad J'(s)=1 \quad\text{for } s \geq 0.8,
\end{equation*}
and set
$$
\jmath'_{\pm}(\xi)=J'\left(\frac{\la\xi\ra\pm \xi_2}{2\la\xi\ra}\right).
$$
Then $\jmath_{\pm}'(\xi)=0$ for $\la\xi_2\ra\le \la\xi_1\ra / 5$, and
\begin{equation}\label{apo}
\jmath'_{+}(\xi)\jmath_{+}(\xi)
=\jmath_{+}(\xi),
\quad
\jmath'_{-}(\xi)\jmath_{+}(\xi)=0.
\end{equation}
With $A=A(x,\xi)$ and $B=B(x,\xi)$ as in \eqref{eq:inter2}, set
\begin{align*}
\tilde{A}(x,\xi_1,\xi_2)&=
\jmath'_{+}(\xi)\left( A(x,\xi_1,\xi_2)
-\frac{i\xi_1}{\la \xi\ra+\la \xi_2\ra}\right)\\
&\quad +\jmath'_{-}(\xi)\left(A(x,\xi_1,-\xi_2)-\frac{i\xi_1}{\la \xi\ra+\la \xi_2\ra}\right) ,\\
\tilde{B}(x,\xi_1,\xi_2)&=\jmath'_{+}(\xi)B(x,\xi_1,\xi_2)+\jmath'_{-}(\xi)B(x,\xi_1,-\xi_2).
\end{align*}
Note that these symbols satisfy the desired properties.

On the symbol level, we have
$$
\la \xi\ra = \la\xi_2\ra  + \frac{\xi_1^2}{\la\xi\ra+\la\xi_2\ra}=
\la\xi_2\ra  - \frac{i\xi_1}{\la\xi\ra+\la\xi_2\ra}i\xi_1.
$$
On the other hand, by the very definition of paradifferential operators, for any couple of symbols $c_1=c_1(x,\xi)$ and
$c_2=c_2(\xi)$ depending only on $\xi$, we have
\begin{equation*}\label{c1c2}
T_{c_1}T_{c_2}=T_{c_1 c_2}.
\end{equation*}
Therefore, by means of \eqref{apo} we easily check that
\begin{equation*}
\left( -i\partial_{x_2} + \nu\partial_{x_1}^{2}+T_{\tilde A} \partial_{x_1}+ T_{\tilde{B}}\right)U=
\left(\la D_{x}\ra  +\nu\partial_{x_1}^2+ T_{A} \partial_{x_1}+T_{B}\right)U .
\end{equation*}
The desired result then follows from \eqref{eq:inter4}
\end{proof}

To prepare for the next transformation, we need a calculus lemma to handle
commutators of the form $[ T_{p}, \partial_{x_1}^{-j}]$.
Note that $\eta(\xi_1)(i\xi_1)^{-j}$
does not belong to our symbol classes. However, we have the following result.

\begin{prop}\label{scxi1}
Let $p\in \Gamma^{0}_{4}(\xT^2)$ and $v\in H^{-\infty}(\xT^2)$
be such that
$$
\partial_{x_1}^{-5}v\in H^{s+2}(\xT^2).
$$
If
\begin{equation}\label{assumeans}
\int_{-\pi}^{\pi} T_{p} v \, dx_1 = 0= \int_{-\pi}^{\pi} v \, dx_1,
\end{equation}
then
\begin{equation*}
\partial_{x_1}^{-1}T_{p} v
= T_{p}\partial_{x_1}^{-1}v-T_{\partial_{x_1}p}\partial_{x_1}^{-2}v
+ T_{\partial_{x_1}^{2}p}\partial_{x_1}^{-3}v-T_{\partial_{x_1}^{3}p}\partial_{x_1}^{-4}v+ f,
\end{equation*}
where $f\in H^{s+2}(\xT^2)$.
\end{prop}
\begin{proof}
We begin by noticing that~\eqref{assumeans} implies that
$$
\partial_{x_1}^0 T_{p} v = T_{p}v ,\quad \partial_{x_1}^{0} v= v,
$$
and hence
$$
\partial_{x_1}\left(\partial_{x_1}^{-1}T_{p}v-T_{p}\partial_{x_1}^{-1}v\right)
=
T_{p}v-T_{\partial_{x_1}p}\partial_{x_1}^{-1}v-T_{p}v=
-T_{\partial_{x_1}p}\partial_{x_1}^{-1}v.
$$
Since $u=\partial_{x_1}U$ implies $U=\partial_{x_1}^{-1}u$, this yields
\begin{equation}\label{ix1-1}
\partial_{x_1}^{-1}T_{p}v-T_{p}\partial_{x_1}^{-1}v
=-\partial_{x_1}^{-1}T_{\partial_{x_1}p}\partial_{x_1}^{-1}v.
\end{equation}
To repeat this argument we first note that, by definition of $\partial_{x_1}^{-1}$, we have
$$
\int_{-\pi}^{\pi} \partial_{x_1}^{-1}v \, dx_1=0.
$$
On the other hand,
\begin{align*}
\int_{-\pi}^{\pi} T_{\partial_{x_1}p}\partial_{x_1}^{-1}v \, dx_1&=
\int_{-\pi}^{\pi} \partial_{x_1}\left( T_{p}\partial_{x_1}^{-1}v \right) - T_{p}\partial_{x_1}^{0} v \, dx_1\\
&=
\int_{-\pi}^{\pi} \partial_{x_1}\left( T_{p}\partial_{x_1}^{-1}v \right) \, dx_1 - \int_{-\pi}^{\pi}T_{p} v \, dx_1\\
&=0,
\end{align*}
by periodicity in $x_1$ and \eqref{assumeans}.
We can thus apply \eqref{ix1-1} with $(p,v)$ replaced by $(\partial_{x_1}p,\partial_{x_1}^{-1}v)$ to obtain
\begin{equation}\label{ix1-2}
\partial_{x_1}^{-1}T_{\partial_{x_1} p}\partial_{x_1}^{-1}v-T_{\partial_{x_1}p}\partial_{x_1}^{-2}v
=-\partial_{x_1}^{-1}T_{\partial_{x_1}^2p}\partial_{x_1}^{-2}v.
\end{equation}
By inserting this result in \eqref{ix1-1} we obtain
$$
\partial_{x_1}^{-1}T_{p} v
= T_{p}\partial_{x_1}^{-1}v-T_{\partial_{x_1}p}\partial_{x_1}^{-2}v
+\partial_{x_1}^{-1}T_{\partial_{x_1}^2p}\partial_{x_1}^{-2}v.
$$
By repeating this reasoning  
we end up with
$$
\partial_{x_1}^{-1}T_{p} v
= T_{p}\partial_{x_1}^{-1}v-T_{\partial_{x_1}p}\partial_{x_1}^{-2}v
+ T_{\partial_{x_1}^{2}p}\partial_{x_1}^{-3}v-T_{\partial_{x_1}^{3}p}\partial_{x_1}^{-4}v+ f ,
$$
where
$$
f= T_{\partial_{x_1}^{4}p}\partial_{x_1}^{-5}v-\partial_{x_1}^{-1}T_{\partial_{x_1}^{4}p}\partial_{x_1}^{-5}v.
$$
By assumption $\partial_{x_1}^{-5} v\in H^{s+2}(\xT^2)$ and $\partial_{x_1}^{4}p\in \Gamma^{0}_{0}(\xT^2)$ so that
$T_{\partial_{x_1}^4p}$ is of order $\le 0$. Therefore we obtain $f\in H^{s+2}(\xT^2)$,
which concludes the proof.
\end{proof}

We have an analogous result for commutators $[T_p,\partial_{x_1}^{-j}]$ for $2\le j\le 4$.

\begin{coro}\label{scxi2}
Let $p\in \Gamma^{0}_{4}(\xT^2)$ and $v\in H^{-\infty}(\xT^2)$
be such that $\partial_{x_1}^{-5}v\in H^{s+2}(\xT^2)$. If
$$
\int_{-\pi}^{\pi} T_{p} v \, dx_1 = 0= \int_{-\pi}^{\pi} v \, dx_1,
$$
then
\begin{align*}
\partial_{x_1}^{-2}T_{p} v
&= T_{p}\partial_{x_1}^{-2}v-T_{\partial_{x_1}p}\partial_{x_1}^{-3}v
+ T_{\partial_{x_1}^{2}p}\partial_{x_1}^{-4}v+ f _2,\\
\partial_{x_1}^{-3}T_{p} v
&= T_{p}\partial_{x_1}^{-3}v-T_{\partial_{x_1}p}\partial_{x_1}^{-4}v
+ f _3,\\
\partial_{x_1}^{-4}T_{p} v
&= T_{p}\partial_{x_1}^{-4}v+ f _4,
\end{align*}
where $f_2,f_3,f_4\in H^{s+2}(\xT^2)$.
\end{coro}
\begin{proof}
This follows from the previous Proposition.
\end{proof}

An important remark for what follows is that $\partial_{x_1}$ and $\partial_{x_2}$
do not have the same weight.
Roughly speaking, the form of the equation \eqref{eq:inter2bis} indicates that
$$
\nu\partial_{x_1}^2 \sim \la D_x\ra \sim \la \partial_{x_2}\ra.
$$
In particular, we shall make extensive use of
$$
\nu^{-1}\partial_{x_1}^{-2} \sim \la D_x\ra^{-1}.
$$
The following lemma gives this statement a rigorous meaning.
\begin{lemm}There holds
$\partial_{x_1}^{-2}U\in H^{s+1}(\xT^2)$ and $\partial_{x_1}^{-4} U \in H^{s+2}(\xT^2)$.
\end{lemm}
\begin{proof}
Recall that
$A(x^\star,\xi^\star)=-A(x,\xi)$ and that $U$ is odd in $x_1$ and $2\pi$-periodic in $x_1$,
so we have
\begin{equation}\label{meansAUU}
\int_{-\pi}^{\pi} T_{A}\partial_{x_1}U\, dx_1=0, \quad
\int_{-\pi}^{\pi} \partial_{x_1}U\, dx_1=0, \quad
\int_{-\pi}^{\pi} U \, dx_1=0.
\end{equation}

Since
$$
\la D_x\ra U + \nu \partial_{x_1}^2 U+ T_A\partial_{x_1}U + T_B U \in H^{s+2}(\xT^2),
$$
by applying $\Lambda^{-1}\partial_{x_1}^{-2}$ with $\Lambda^{-1}=T_{\la\xi\ra^{-1}}$,
we have
\begin{align*}
\partial_{x_1}^{-2} U
= -\nu \Lambda^{-1}U - \Lambda^{-1} \partial_{x_1}^{-2}(T_A \partial_{x_1} U )
-\Lambda^{-1} \partial_{x_1}^{-2}(T_B  U ) +  F,
\end{align*}
where $F\in H^{s+2}(\xT^2)$.
The first term and the third term
in the right hand side obviously belong to $H^{s+1}(\xT^2)$.
Moving to the second term in the right-hand side, in view of~\eqref{meansAUU}, the argument establishing \eqref{ix1-1} also gives
$$
\partial_{x_1}^{-2}T_A \partial_{x_1}U=\partial_{x_1}^{-1} T_A U- \partial_{x_1}^{-1}T_{\partial_{x_1}A} U.
$$
Hence, using that $\partial_{x_1}^{-1},T_{A}$ and $T_{\partial_{x_1}A}$ are of order $\le 0$,
we have that
$$
\partial_{x_1}^{-2}T_A \partial_{x_1}U \in H^{s}(\xT^2),
$$
so that $\Lambda^{-1} \partial_{x_1}^{-2}(T_A \partial_{x_1} U ) \in H^{s+1}(\xT^2)$ and hence
\begin{equation}\label{L1d12}
\partial_{x_1}^{-2} U \in H^{s+1}(\xT^2).
\end{equation}

To study $\partial_{x_1}^{-4}U$ we start from
\begin{align*}
\partial_{x_1}^{-4} U
&= -\nu \Lambda^{-1}\partial_{x_1}^{-2}U - \Lambda^{-1} \partial_{x_1}^{-4}(T_A \partial_{x_1} U ) \\
&\quad -\Lambda^{-1} \partial_{x_1}^{-4}(T_B  U ) + \partial_{x_1}^{-2} F,
\end{align*}
We have just proved that the first term in the right hand side belongs to $H^{s+2}(\xT^2)$.
With regards to the second term
we use the third identity in
Corollary~\ref{scxi2} to obtain
$$
\partial_{x_1}^{-4}(T_A \partial_{x_1} U ) -
T_{A} \partial_{x_1}^{-3}U \in H^{s+2}(\xT^2).
$$
On the other hand,
by symbolic calculus, $ \Lambda^{-1} T_{A}-T_{A} \Lambda^{-1} $ is of order $\le -2$. Hence,
$$
 \Lambda^{-1} \partial_{x_1}^{-4}(T_A \partial_{x_1} U ) -
T_{A}  \Lambda^{-1} \partial_{x_1}^{-3}U \in H^{s+2}(\xT^2).
$$
In view of \eqref{L1d12} this yields
$$
\Lambda^{-1} \partial_{x_1}^{-4}(T_A \partial_{x_1} U ) \in H^{s+2}(\xT^2).
$$
Similarly we obtain that $\Lambda^{-1} \partial_{x_1}^{-4}(T_B U ) \in H^{s+2}(\xT^2)$. We thus end up with
$$
\partial_{x_1}^{-4} U \in H^{s+2}(\xT^2),
$$
which concludes the proof.
\end{proof}

The following definition is helpful for what follows.

\begin{defi} We say that an operator $R$ is of anisotropic order $\le -2$ if $R$ is of the form
$$
R=R_0\Lambda^{-2} + R_{1}\Lambda^{-1}\partial_{x_1}^{-2}
+ R_{2}\partial_{x_1}^{-4},
$$
where $R_0,R_1$ and $R_2$ are operators of order $\le 0$ and
$$
\Lambda^{-1} = T_{\la \xi\ra^{-1}}, \quad
\Lambda^{-2} = T_{\la \xi\ra^{-2}}.
$$
\end{defi}
It follows from the previous lemma that operators of anisotropic order $\le -2$ may be seen
as operators of order $\le -2$. Namely, the previous lemma implies the following result.
\begin{lemm}\label{anisotropic}
If $R$ is of anisotropic order $\le -2$, then
$R U \in H^{s+2}(\xT^2)$.
\end{lemm}

With these preliminaries established,
to prove Proposition~\ref{p:reduced}
the key point is a decomposition of zero-order symbols.
We want to decompose zero-order operators as sums of the form
$$
\sum_{0\le j\le 4}T_{a_j} \partial_{x_1}^{-j}+ R,
$$
where $T_{a_j}$ are paraproducts ($a_j=a_j(x)$ does not depend on $\xi$) 
and $R$ is of anisotropic
order $\le -2$ plus an admissible remainder. More generally, we consider below symbols of order $-2\le m\le 0$ and not only zero-order symbols.

\begin{lemm}\label{d0os}
Let $m\in \{0,1,2\}$ and $\rho\geq 0$. Let
$S\in \Gamma^{-m}_{\rho}(\xT^2)$
be an homogeneous symbol of degree $-m$ in $\xi$ such that
$$
S(x,\xi_1,-\xi_2)=S(x,\xi_1,\xi_2),
$$
and such that, for some positive constant $c$,
$$
\la\xi_2\ra\le c \la\xi_1\ra \Rightarrow S(x,\xi)=0.
$$
Then, for all $v$ whose spectrum is included in the semi-cone $\{\xi_2\ge c\la \xi_1\ra\}$,
\begin{equation}\label{JcSQ0}
T_{S(x,\xi)}\partial_{x_1}v
=
\sum_{j=2m}^{4}T_{S_{j} (x)}\partial_{x_1}^{1-j}v 
+ Q
(-i \partial_{x_2}+\nu\partial_{x_1}^{2})v
+ Rv,
\end{equation}
where $R$ is of anisotropic order $\le -2$,
\begin{equation}\label{defisj}
S_j(x)=\frac{1}{ i^j \nu^j j!}(\partial_{\xi_1}^j S)(x,0,1),
\end{equation}
and %
$$
Q=\sum_{k=2m+1}^{4}T_{q_{k}(x,\xi)}\partial_{x_1}^{1-k}, 
$$
where $q_k\in \Gamma^{-1}_{\rho}(\xT^2)$ for $1\le k\le 4$ is 
explicitly defined in \eqref{defiq} below
and satisfies
\begin{equation}\label{symq1}
q_k(x,\xi_1,-\xi_2)=q_k(x,\xi_1,\xi_2),
\end{equation}
and
\begin{equation}\label{symq2}
\la\xi_2\ra\le \frac{c}{2} \la\xi_1\ra \Rightarrow q_k(x,\xi)=0.
\end{equation}
\end{lemm}
\begin{rema}
This is a variant of the decomposition used in \cite{IP}. The main difference is that, having performed
the reduction to the case where $\xi_2\ge \la\xi_1\ra/2$, we do not need to consider the so-called elementary operators
in~\cite{IP}. Hence we obtain a decomposition where the $S_j$'s do not depend on~$\xi$.
\end{rema}
\begin{proof}
The proof is based on the following simple observation: $\xi\mapsto \xi_1$ is transverse
to the characteristic variety. We prove Lemma~\ref{d0os} for $m=0$.

Let $J_{c}$ be a real-valued function, homogeneous of degree $0$ such that 
$J_{c}(\xi_{1},\xi_{2})=J_{c}(\xi_{1},-\xi_{2})$ and
\begin{align*}
J_c(\xi) =0 \quad\text{for } \la \xi_2\ra \le \frac{c}{2}\la\xi_1\ra,
\quad 
J_c(\xi) =1 \quad\text{for } \la \xi_2 \ra\ge c\la\xi_1\ra
\end{align*}
(so that $J_c(\xi)S(x,\xi)=S(x,\xi)$ and $J_c(D_x)v=v$). We shall prove that
\begin{equation}\label{JcSQ}
T_{S(x,\xi)}\partial_{x_1}^{0}
=
\sum_{j=0}^{4}T_{S_{j} (x)}\partial_{x_1}^{-j}J_c 
+ \sum_{k=1}^{4}T_{q_{k}(x,\xi)}\partial_{x_1}^{-k} 
(\la\partial_{x_2}\ra+\nu\partial_{x_1}^{2})
+ R',
\end{equation}
where $S_j,q_k$ are as in the statement of the lemma and $R'$ is such that the operator $R'\partial_{x_1}$ is of 
anisotropic order $\le -2$. This yields the desired result \eqref{JcSQ0} since 
$J_c v=v$ and $\la\partial_{x_2}\ra v=-i\partial_{x_2}v$ by assumption on the spectrum of$v$.

For $m=0$, $S$ is an homogeneous symbol of degree $0$ in $\xi$. 
By the symmetry hypothesis $S(x,\xi_1,\xi_2)=S(x,\xi_1,-\xi_2)$, we can write $S(x,\xi)$ as
$$
S(x,\xi)=
S\left( x,\frac{\xi_1}{\la\xi_2\ra},1\right),
$$
(if $\xi_2=0$ then $S(x,\xi)=0$ by assumption) so we have
\begin{equation}\label{tausplit}
S(x,\xi)=
\sum_{j=0}^{4} \frac{1}{i^j j!}(\partial_{\xi_1}^j \tilde{S})(x,0,1)
\left(i \frac{\xi_1}{\la\xi_2\ra}\right)^j
+  r\left(x,\frac{\xi_1}{\la\xi_2\ra}\right)
\left(i \frac{\xi_1}{\la\xi_2\ra}\right)^5 ,
\end{equation}
where $r$ is given by Taylor's formula.

Next, by setting
$$
\mathcal{L}(\xi)\defn \la\xi_2\ra-\nu\xi_1^2,
$$
we claim that, for $1\le k\le 4$, there exists
a Fourier multiplier $Q_{k}(D_x)$ of order $\le -1$ such that,
$$
J_{c}(\xi) \,  i\xi_1 \left(i\frac{\xi_1}{\la\xi_2\ra}\right)^{k}= -\frac{J_{c}(\xi)}{i^{k-1} \nu^k \xi_1^{k-1}}
+  Q_{k}(\xi)\mathcal{L}(\xi).
$$
To see this, write
$$
\frac{\xi_1^2}{\la\xi_2\ra} = -\frac{\mathcal{L}(\xi)}{\nu\la\xi_2\ra} +\frac{1}{\nu},\quad
\frac{\xi_1^4}{\la\xi_2\ra^2} = \left(-\frac{\mathcal{L}(\xi)}{\nu\la\xi_2\ra} +\frac{1}{\nu}\right)^2,\quad
\frac{\xi_1^{6}}{\la\xi_2\ra^3} = \left(-\frac{\mathcal{L}(\xi)}{\nu\la\xi_2\ra} +\frac{1}{\nu}\right)^3,
$$
to obtain the desired identities with $Q_j(\xi)=J_{c}(\xi) \tilde Q_j (\xi)$ where
\begin{alignat*}{3}
\tilde Q_1(\xi)&=\frac{1}{\nu\la\xi_2\ra}, &
\quad \tilde Q_{3}(\xi)&=-\frac{\nu\la\xi_2\ra\xi_1^2+\la\xi_2\ra^2+\nu^2\xi_1^{4}}{\nu^3\xi_1^2\la\xi_2\ra^3},\\%
\tilde Q_{2}(\xi)&=-\frac{\la\xi_2\ra+\nu\xi_1^2}{i\nu^2\la\xi_2\ra^2 \xi_1},
& \quad  \tilde Q_4(\xi) &= -\frac{i\left(
c_1\la\xi_2\ra^3 +c_2 \la\xi\ra^2\xi_1^2+
c_3\la\xi_2\ra\xi_1^4+c_4 \xi_1^6\right)}{\nu^{4}\xi_1^3\la\xi_2\ra^4},
\end{alignat*}
where $c_1\defn 3-6\nu+4\nu^2$, $c_2\defn 3\nu+6\nu^2-8\nu^3$, $c_3\defn 3\nu^2+4\nu^4$ and $c_4\defn \nu^3 $.
Similarly, we have
$$
i\xi_1 \left(i\frac{\xi_1}{\la\xi_2\ra}\right)^{5}=
-\frac{1}{ \nu^3\la\xi_2\ra^2}  + \tilde Q_{5}(\xi)\mathcal{L}(\xi) \text{ with }
\tilde Q_{5}(\xi)\defn
\frac{\la\xi_2\ra^2+4\nu\la\xi_2\ra\xi_1^2+\nu^2\xi_1^{4}}{\nu^3 \la\xi_2\ra^5}.
$$

Our analysis of \eqref{tausplit} is complete; by inserting
the previous identities in \eqref{tausplit} premultiplied by $J_c(\xi)\eta(\xi_1)i\xi_1$ and next by 
dividing by $i\xi_1$, we obtain the desired decomposition~\eqref{JcSQ} with $q_k(x,\xi)=J_{c}(\xi) \tilde q_{k}(x,\xi)$ where
\begin{equation}\label{defiq}
\begin{aligned}
\tilde q_1(x,\xi)=
&\frac{S_1(x)}{\la\xi_2\ra}
-\frac{\nu S_3(x)}{\la\xi_2\ra^2}\\
&+i\xi_1\left( \frac{i\nu^2 S_3(x)\xi_1}{\la\xi_2\ra^3}
-\frac{\nu S_2(x)}{\la\xi_2\ra^2}-\frac{(c_3 \la \xi_2\ra+c_4 \xi_1^2)S_4(x)}{\la \xi_2\ra^4}\right)\\
&+ r\left(x,\frac{\xi_1}{\la\xi_2\ra}\right)\tilde Q_{5}(\xi), \\
\tilde q_2(x,\xi)=
&
\frac{-S_2(x)}{\nu^{2}\la\xi_2\ra}-i\frac{c_2 S_4}{\la\xi_2\ra^{2}},\\
\tilde q_3(x,\xi)=
&\frac{S_3(x)}{\nu^{3}\la\xi_2\ra},\\
\tilde q_4(x,\xi)=&\frac{-c_1 S_4(x)}{\nu^{4}\la\xi_2\ra}.
\end{aligned}
\end{equation}
Note that each term making up $q_k(x,\xi)$ is well-defined and $C^\infty$ for $\xi\neq 0$ and homogeneous of degree
$-1$ or $-2$ or $-3$ in $\xi$, so
$q_k\in \Gamma^{-1}_{\rho}(\xT^2)$.
\end{proof}

We are now prepared to conclude the proof of Proposition~\ref{p:reduced}.
We want to prove that there exist five functions
$$
a_j=a_j(x)\in C^{s-6-j}(\xT^2) \qquad (0\le j\le 4),
$$
where
\begin{equation*}
a_{j} \text{ is odd in }x_1 
\text{ for }j\in \{0,2,4\},\quad
a_{j} \text{ is even in }x_1 
\text{ for }j\in \{1,3\},
\end{equation*}
such that
\begin{equation}\label{dea}
\left( -i\partial_{x_2} + \nu\partial_{x_1}^2 +
\sum_{0\le j\le 4}T_{a_j} \partial_{x_1}^{1-j} \right) U\in H^{s+2}(\xT^2).
\end{equation}

To this end, since
$\tilde A$ and $\tilde{B}$ satisfy properties iii) and iv) in Lemma~\ref{ltAB}, 
since $U=\partial_{x_1}^0U$ 
and since the spectrum of $U$ is contained in the semi-cone $\{\xi_2\ge \frac{3}{4}\la \xi_1\ra\}$,  
we can use the above symbol-decomposition process to obtain
for $U$ an equation of the form
$$
\left(I+Q\right) \left( -i\partial_{x_2} U  + \nu\partial_{x_1}^2 U\right)
+\sum_{0\le j\le 4} T_{\alpha_j}\partial_{x_1}^{1-j} U =f
\in H^{s+2}(\xT^2),
$$
where $\alpha_j\in C^{s-6}(\xT^2)$ and 
$$
Q=\sum_{k=1}^{4}T_{q_{k}(x,\xi)}\partial_{x_1}^{1-k},
$$
where $q_k \in \Gamma^{-1}_{s-6}(\xT^2)$.  

Write
$$
\left(I+Q\right)\left(-i\partial_{x_2} U  +
\nu\partial_{x_1}^2 U+\sum_{0\le j\le 4} T_{\alpha_j}\partial_{x_1}^{1-j} U\right)
-Q \sum_{0\le j\le 4} T_{\alpha_j}\partial_{x_1}^{1-j} U =f.
$$
Now, by using Proposition~\ref{scxi1} and its corollary~\ref{scxi2}, 
notice that we can write the term 
$Q\sum_{0\le j\le 4} T_{\alpha_j}\partial_{x_1}^{1-j}U$ 
under the form
$$
\sum\sum\sum(-1)^\ell T_{q_k}T_{\partial_{x_1}^\ell\alpha_j}\partial_{x_1}^{2-j-k-\ell} U +F,
$$
where $F\in H^{s+2}(\xT^2)$ where the sum is taken over 
indices $k,j,\ell$ such that $1\le k\le 4$, $0\le j\le 4$ and 
$$
j+k+\ell \le 3.
$$
Indeed, for those indices such that $j+k+\ell \ge 4$ the operator 
$T_{q_k}T_{\partial_{x_1}^\ell\alpha_j}\partial_{x_1}^{2-j-k-\ell}$ is of anisotropic index $\le -2$ and hence 
$T_{q_k}T_{\partial_{x_1}^\ell\alpha_j}\partial_{x_1}^{2-j-k-\ell} U \in H^{s+2}(\xT^2)$, so that these terms  
can be handled as source terms. Next, noticing that  $T_{q_k}T_{\partial_{x_1}^\ell\alpha_j}-T_{q_k \partial_{x_1}^\ell\alpha_j}$ 
is order less than $ -3$, this implies that we have
$$
Q\sum_{0\le j\le 4} T_{\alpha_j}\partial_{x_1}^{1-j}U=\sum_{0\le p\le 4} T_{\mathcal{S}_p}\partial_{x_1}^{1-p}.
$$
for some symbols $\mathcal{S}_p$ of order $-1$ in $\xi$, with regularity $s-8$ in $x$ and satisfying the 
spectral assumptions in Lemma~\ref{d0os}. 
As a consequence, 
by applying Lemma~\ref{d0os} we obtain 
$$
Q \sum_{0\le j\le 4} T_{\alpha_j}\partial_{x_1}^{1-j} U= \sum_{2\le k\le 4} T_{\beta_k(x)}\partial_{x_1}^{1-k}U+ Q'
\left( -i\partial_{x_2} + \nu\partial_{x_1}^2\right)U+R'U,
$$
where $\beta_k\in C^{s-8}(\xT^2)$, $R'$ is of anisotropic order $\le -2$ and $Q'$ is of the form
$$
Q'=\sum_{k=3}^{4}T_{q'_{k}(x,\xi)}\partial_{x_1}^{1-k},
$$
where $q'_k(x,\xi)$ is of order $-1$ in $\xi$ (and regularity $C^{s-8}$ in $x$). 
Then
\begin{multline*}
\left(I+Q-Q'\right)\left(-i\partial_{x_2} + \nu\partial_{x_1}^2
+\sum_{0\le j\le 4} T_{\alpha_j}\partial_{x_1}^{1-j} -\sum_{2\le k\le 4} T_{\beta_k(x)}\partial_{x_1}^{1-k}\right)
U \\
+(Q-Q')\sum_{2\le k\le 4} T_{\beta_k(x)}\partial_{x_1}^{1-k}U+Q' \sum_{0\le j\le 4} T_{\alpha_j}\partial_{x_1}^{1-j}
\in H^{s+2}(\xT^2).
\end{multline*}
Again, one has
$$
(Q-Q')\sum_{2\le k\le 4} T_{\beta_k(x)}\partial_{x_1}^{1-k}U+Q' \sum_{0\le j\le 4} T_{\alpha_j}\partial_{x_1}^{1-j}
=T_{\gamma_4 (x)}\partial_{x_1}^{-3}+R''U.
$$
Now $(Q-Q')T_{\gamma_4 }\partial_{x_1}^{-3}$ is of anisotropic order $\le -2$. 
This yields
$$
\left(I+\mathcal{Q}\right)
\left(-i\partial_{x_2} U + \nu\partial_{x_1}^2 U+\sum_{0\le j\le 4} T_{a_j}\partial_{x_1}^{1-j} U
 \right)
\in H^{s+2}(\xT^2),
$$
where $\mathcal{Q}=Q-Q'$ and
$$
a_0=\alpha_0, ~a_1=\alpha_1,~a_2=\alpha_2-\beta_2,~a_3=\alpha_3-\beta_3,~a_4=\alpha_4-\beta_4+\gamma_4.
$$
Now we have an obvious left parametrix for $I+\mathcal{Q}$, in the following sense:
$$
\left(I-\mathcal{Q}+\mathcal{Q}^{2}-\mathcal{Q}^{3}\right)\left(I+\mathcal{Q}\right)
=I - \mathcal{Q}^{4},
$$
where $\mathcal{Q}^{4}$ is of order $\le -4$ so that
$$
\mathcal{Q}^{4}
\left( -i\partial_{x_2}U + \nu\partial_{x_1}^2 U
+\sum_{0\le j\le 4} T_{a_j}\partial_{x_1}^{1-j} U\right) \in H^{s+2}(\xT^2).
$$
This gives \eqref{dea}.

The symmetries of the coefficients $a_j$ can be checked on the
explicit expressions which are involved.
Indeed, it follows from \eqref{defisj}
that the function $s=(s_0,\ldots,s_4)$ given by Lemma~\ref{d0os}
satisfies the same symmetry as $S$ does: given $\eps\in \{-1,+1\}$ and $0\le j \le 4$,
we have
$$
S(x^\star,\xi^\star)=\eps S(x,\xi) ~\Rightarrow
S_{j}(x^\star)
=
\eps (-1)^{j}S_j(x).
$$
This concludes the proof of Proposition~\ref{p:reduced}.

\subsection{Proof of Proposition~\ref{prop440}}\label{psr}

Given Proposition~\ref{p:reduced}, the proof of Proposition~\ref{prop440}
reduces to establishing the following result.

\begin{nota}\label{529}
Given five complex-valued functions $a_0,\ldots,a_4$, we define
$$
Z_a = \sum_{0\le j\le 4}T_{a_j} \partial_{x_1}^{-j}.
$$
\end{nota}

\begin{prop}\label{prop44}
There exist two constant $\kappa,\kappa'\in\xR$ and five functions $c_0,\ldots,c_4$ with $c_j\in C^{6-j}(\xT^2)$ satisfying $\la c_0 \ra>0$ and
\begin{equation}\label{symc}
c_{k} \text{ is even in $x_1$ for } k\in \{0,2,4\},
\quad
c_{k} \text{ is odd in $x_1$ for } k\in \{1,3\},
\end{equation}
such that
\begin{equation}\label{canonique}
Z_c \Bigl(-i\partial_{x_2} +\nu \partial_{x_1}^2 +   Z_a\partial_{x_1}\Bigr)U =
\Bigl(-i\partial_{x_2}+\nu \partial_{x_1}^2  + \kappa +\kappa'\partial_{x_1}^{-2}\Bigr) Z_c U
+f,
\end{equation}
where $f\in H^{s+2}(\xT^2)$.
\end{prop}

\begin{proof}
The equation \eqref{canonique} is equivalent to a sequence of five transport equations
for the coefficients $c_j$ ($0\le j\le 4$), which
can be solved by induction. Indeed,
directly from the Leibniz rule we compute that
$$
\bigl(-i\partial_{x_2}u +\nu \partial_{x_1}^2+\kappa+\kappa'\partial_{x_1}^{-2}\bigr) Z_{c}u -
Z_{c} \bigl(-i\partial_{x_2}U +\nu \partial_{x_1}^2
\bigr)u = Z_{\delta}\partial_{x_1}U,
$$
where
\begin{align*}
\delta_{0}&=  2\nu\partial_{x_1}c_0 ,\\
\delta_{1}&= 2\nu \partial_{x_1}c_{1}-i\partial_{x_2}c_0+\nu \partial_{x_1}^{2} c_0
+c_0\kappa,\\
\delta_{2}&= 2\nu \partial_{x_1}c_2 -i \partial_{x_2}c_{1}
+\nu \partial_{x_1}^{2}c_{1}+c_1\kappa,\\
\delta_{3}&= 2\nu \partial_{x_1}c_3 -i \partial_{x_2}c_{2}
+\nu \partial_{x_1}^{2}c_{2}+c_2\kappa+c_0\kappa',\\
\delta_{4}&= 2\nu \partial_{x_1}c_4 -i \partial_{x_2}c_{3}
+\nu \partial_{x_1}^{2}c_{3}+c_3\kappa+c_1\kappa'.
\end{align*}
On the other hand, if \eqref{symc} is satisfied,
then we can apply Proposition~\ref{scxi1} and its corollary to obtain
$$
Z_c Z_a \partial_{x_1}U =Z_{\delta'}\partial_{x_1}U+ f,
$$
where $f\in H^{s+2}(\xT^2)$ and
$$
\delta'_k= \sum_{l+m+n=k}c_l (-\partial_{x_1})^{m}a_n
\qquad (0\le k,l,m,n\le 4).
$$
Hence, our purpose is to define $c=(c_0,\ldots,c_4)$ satisfying \eqref{symc} and two constants $\kappa$ and $\kappa'$ such that
$$
\delta=\delta'.
$$

\smallbreak
\noindent{{\sc Step} 1: Definition of $c_0$.}
We first define the principal symbol $c_0$ by solving the equation $\delta_{0}=\delta_{0}'$, which reads
\begin{equation*}
2\nu \partial_{x_1} c_{0} =  c_0 a_0.
\end{equation*}
We get a unique solution of this equation by imposing the initial condition $c_{0}(0,x_2)=C_{0}(x_2)$ on $x_1=0$,
where $C_0$ is to be determined.
That is, we set
$$
c_0(x)=C_0(x_2) e^{\gamma(x)/(2\nu)},
$$
where
$$
\gamma\defn \partial_{x_1}^{-1} a_0.
$$
Since $a_0$ is odd in $x_1$, we have $\int_{-\pi}^{\pi}a_0\, dx_1=0$ and hence
$$
\partial_{x_1}\gamma = a_0.
$$
Note that, directly from the definition, we have
$$
\gamma\in C^{s-6}(\xT^2),\quad
\gamma\in C(e,e), \quad \int_{-\pi}^{\pi} \gamma  \, dx_1=0,
\quad \gamma(x)\in i \xR.
$$

\smallbreak
\noindent{{\sc Step} 2: Definition of $c_1$, $C_0$ and $\kappa$.}
We next define $c_1$ by solving $\delta_{1}=\delta_{1}'$. This yields
\begin{equation*}
2\nu \partial_{x_1} c_{1}- a_0 c_1 =G_1
\end{equation*}
with
$$
G_1\defn i \partial_{x_2} c_{0}- \nu \partial_{x_1}^2 c_{0}
- \kappa c_0 +c_0 a_1 -c_{0} \partial_{x_1}a_0
$$
where $\kappa$ is determined later.
We impose the initial condition $c_1(0,x_2)=0$ on $x_1=0$, so that
$$
c_1(x_1,x_2)\defn \frac{1}{2\nu} e^{\gamma/(2\nu)} \int_{0}^{x_1} e^{-\gamma/(2\nu)} G_{1} \, ds.
$$

Note that $c_1$ is $2\pi$-periodic in $x_1$ if and only if
\begin{equation}\label{condG1}
\int_{0}^{2\pi} e^{-\gamma/(2\nu)} G_{1} \, dx_1 = 0.
\end{equation}
Directly from the definition of $G_{1}$, we compute that
\begin{equation*}
G_{1}= e^{\gamma/(2\nu)} \left[ i\partial_{x_2} C_0
+ C_0 \left( \frac{i}{2\nu}\partial_{x_2}\gamma -\kappa +a_1-\frac{3}{2}\partial_{x_1}a_0
-\frac{1}{4\nu} a_0 ^2\right)\right],
\end{equation*}
Using
$$
\int_{0}^{2\pi} \partial_{x_1}a_{0}\, dx_1=0 =\int_{0}^{2\pi}
 \gamma \, dx_1,
$$
this gives
\begin{equation*}
\int_{0}^{2\pi} e^{\gamma/(2\nu)} G_{1}\, dx_1
= 2i\pi  C_{0}'(x_2) +C_{0}(x_2)\int_{0}^{2\pi}
\left[-\kappa +a_1 -\frac{a_{0}^2}{4\nu} \right] \, dx_1.
\end{equation*}
Set
\begin{equation}\label{defi:beta}
\beta(x_2)\defn
-2\pi \kappa + \int_{0}^{2\pi}\left( a_1 -\frac{a_{0}^2}{4\nu}\right) \, dx_1.
\end{equation}
so that
\begin{equation*}
\int_{0}^{2\pi} e^{-\gamma/(2\nu)} G_{1} \, dx_1
= 2i\pi  C_0'(x_2)+\beta(x_2)C_0(x_2).
\end{equation*}

We thus define $\kappa$ by
\begin{equation}\label{defikappa}
\kappa = \frac{1}{\la\xT^2\ra}\iint_{\xT^2} \left(a_1  -\frac{a_{0}^2}{4\nu}\right)\, dx_1 dx_2.
\end{equation}
With this choice, we have
$$
\int_{0}^{2\pi\periode}\beta(s)\, ds=0.
$$
and hence
$$
C_0(x_2)\defn \exp\left(-\frac{1}{2i\pi}\int_0^{x_2}\beta(s)\, ds\right)
\text{ is $2\pi\periode$-periodic in $x_2$.}
$$
 With this particular choice of $C_0$,
the condition \eqref{condG1} is satisfied
and hence $c_1$ is bi-periodic.

Moreover, directly from these definitions, we have $C_0\in C^{6}$ and, by performing
an integration by parts to handle the term
$\int_{0}^{x_1}e^{\gamma/(2\nu)}(\partial^2 c_{1} / \partial {x_1}^2)\, ds
$, we obtain that $c_1\in C^{5}(\xT^2)$.

\smallbreak

\bigbreak
\noindent{\sc Step} 3: $\kappa\in\xR$.
It remains to prove that $\kappa\in \xR$.
To do this, we first observe that
$a_{0}(x)= A(x,0,1)$ where $A$ is given by Proposition~\ref{C514}.
In particular we easily check that $a_0(x)\in i \xR$ so that
$a_{0}(x)^2 \in \xR$. On the other hand, we claim that
\begin{equation}\label{a1o}
\IM a_{1}(x) \text{ is odd in $x_2$},
\end{equation}
so that
\begin{equation*}
\kappa = \frac{1}{\la\xT^2\ra}\iint_{\xT^2} \left(\RE a_1  -\frac{a_{0}^2}{4\nu}\right)\, dx_1 dx_2 \in \xR.
\end{equation*}
To prove \eqref{a1o}, still with the notations
of Proposition~\ref{C514} and Lemma~\ref{ltAB}, we first observe that the identity \eqref{defisj} and the definition of $a_1$ imply that
\begin{align*}
a_{1}(x)&= \frac{1}{i\nu} (\partial_{\xi_1} \tilde A)(x,\xi)
+\tilde B_0 (x,\xi)\bigg\arrowvert_{\xi=(0,1)}\\
&= \frac{1}{i\nu} \partial_{\xi_1}\left( A(x,\xi)
-\frac{i\xi_1}{\la \xi\ra+\la \xi_2\ra}\right)+B_0(x,\xi)\bigg\arrowvert_{\xi=(0,1)}\\
&=\frac{1}{i\nu} (\partial_{\xi_1} A)(x,0,1)  -\frac{1}{2\nu} + B_{0}(x,0,1),
\end{align*}
so that $\IM a_1(x)= \IM B_0 (x,0,1)$. 

Now, we claim that 
\begin{equation}\label{symb}
\IM B_0(x,-\xi)=-\IM B_{0}(x,\xi).
\end{equation}
Indeed, this follows from the definition of $B_0$ (cf Proposition~\ref{C514}) 
and the following symmetry of the symbol $\dns$ of the Dirichlet to Neumann operator
\begin{equation}\label{symdns}
\overline{\dns (x,\xi)}=\dns (x,-\xi).
\end{equation}
That \eqref{symdns} has to be true is clear since this symmetry means nothing more than the fact that $\DNS f$ is real-valued
for any real-valued function $f$.

Once \eqref{symb} is granted, 
using $B\in \Gamma(e,e)$ we obtain the desired result:
\begin{align*}
\IM a_1(x_1,-x_2)&=
\IM B_{0}(x_1,-x_2,0,1)\\
&=-\IM B_0(x_1,-x_2,0,-1)\\
&=-\IM B_0(x_1,x_2,0,1)\\
&= -\IM a_1 (x_1,x_2).
\end{align*}

\bigbreak
\noindent{\sc Step} 4: General formula. We can now give the scheme of the analysis.
For $k= 2,3,4$,
we shall define $c_{k}$ inductively by
$$
2\nu \partial_{x_1}\Big(e^{-\gamma/(2\nu)}c_{k}\Big) =
e^{-\gamma/(2\nu)} G_{k},
$$
where $G_{k}$ is to be determined by means of the equation
$\delta_k=\delta_k'+\delta_{k-1}''$.
That is, we set
$$
c_{k}(x_1,x_2)=\exp\left(\frac{\gamma(x_1,x_2)}{2\nu}\right) \bigl( C_{k}(x_2)+\Gamma_{k}(x_1,x_2)\bigr),
$$
where $C_k$ is to be determined and $\Gamma_{k}$ is given by
$$
\Gamma_{k}(x_1,x_2)=\frac{1}{2\nu}\int_{0}^{x_1}
\exp\left(\frac{-\gamma(s,x_2)}{2\nu}\right)G_{k}(s,x_2) \, ds.
$$

As in the previous step, we have to chose the initial data
$C_{k-1}(x_2)=c_{k-1}(0,x_2)$ such that
$\Gamma_k$ is $2\pi$-periodic in $x_1$. Now we note the following fact:
Starting from the fact that $a_0,a_2,a_4$ are odd in $x_1$ and $a_1,a_3$
are even in $x_1$, we successively check that:
$c_1$ is odd in $x_1$; $G_2$ is odd in $x_1$;
$c_2$ is even in $x_1$; $G_3$ is even in $x_1$; $c_3$ is odd in $x_1$; $G_4$ is odd in $x_1$.
As a result, we have
$$
\int_{-\pi}^{\pi} e^{-\gamma/(2\nu)} G_{2} \, dx_1 = 0=\int_{-\pi}^{\pi} e^{-\gamma/(2\nu)} G_{4} \, dx_1,
$$
which in turn implies that $\Gamma_2$ and $\Gamma_4$ are bi-periodic.
Consequently, one can impose
$$
C_1(x_2)=c_1(0,x_2)=0 \quad\text{and}\quad C_3(x_2)=c_3(0,x_2)=0.
$$
Moreover, we impose $C_4=0$ (there is no restriction on $C_4$ since we stop the expansion at this order).
Therefore, it remains only to prove that one can so define $C_2$ and $\kappa'$ that $\Gamma_3$ is bi-periodic.

\smallbreak
\noindent {\sc Step }5: Definition of $C_2$ and $\kappa'$.
We turn to the details and compute that
$$
G_3= i\partial_{x_2}c_2- \nu \partial_{x_1}^2 c_{2}
-(\kappa +\partial_{x_1}a_0  -a_1)c_2-\kappa' c_0 + c_0 a_3 +c_1 a_2 -c_1\partial_{x_1}a_1.
$$

The function $c_3$ is $2\pi$-periodic in $x_1$ if and only if
\begin{equation}\label{condG3}
\int_{0}^{2\pi} e^{-\gamma/(2\nu)} G_{3} \, dx_1=0.
\end{equation}
Directly from the definition of $G_3$, we have
\begin{align*}
&\int_{0}^{2\pi} e^{-\gamma/(2\nu)} G_{3}\, dx_1 \\
&\qquad\qquad=
2i\pi C_{2}'(x_2) +C_{2}(x_2)\int_{0}^{2\pi}
\left(-\kappa +a_1 -\frac{a_{0}^2}{4\nu} \right) \, dx_1 \\
&
\qquad\qquad\quad
+ \Gamma_{2}(x_1,x_2)\int_{0}^{2\pi}
\left( \frac{i}{2\nu}\partial_{x_2}\gamma -\kappa +a_1-\frac{3}{2}\partial_{x_1}a_0
-\frac{1}{4\nu} a_0 ^2\right)\, dx_1\\
&\qquad\qquad\quad + \int_{0}^{2\pi} \left( i\partial_{x_2}\Gamma_2 -\nu \partial_{x_1}^2 \Gamma_2
- \partial_{x_1}\gamma \partial_{x_1}\Gamma_2 \right)\, dx_1\\
&\qquad\qquad\quad +\int_{0}^{2\pi} (a_3 -\kappa')C_{0}\, dx_1.
\end{align*}
Now since $a_3$ is odd in $x_1$, the last term simplifies to $-2\pi\kappa' C_{0}(x_2)$. Note also that
$$
\int_{0}^{2\pi} \partial_{x_1}^2 \Gamma_2 \, dx_1=0 ,
$$
and
$$
\int_{0}^{2\pi} - \partial_{x_1}\gamma \partial_{x_1}\Gamma_2 \, dx_1
=\int_{0}^{2\pi} - a_0 \partial_{x_1}\Gamma_2 \, dx_1
=\int_{0}^{2\pi} \Gamma_2 \partial_{x_1}a_0  \, dx_1.
$$
By using the previous cancelations, we obtain that for
\eqref{condG3} to be true, a sufficient condition is that 
$C_2$ solves the equation
\begin{equation}\label{eqdefiC2}
2i\pi  C_{2}'(x_2)+\beta(x_2)
C_{2}(x_2)=F_{2}(x_2)-2\pi\kappa' C_{0}(x_2),
\end{equation}
with
$$
F_{2}(x_2)\defn\int_{0}^{2\pi}\left( \frac{i}{2\nu}\partial_{x_2}\gamma -\kappa +a_1-\frac{1}{2}\partial_{x_1}a_0
-\frac{1}{4\nu} a_0 ^2\right)\Gamma_{2} +i\partial_{x_2}\Gamma_2 \, dx_1.
$$
If we impose the initial condition $C_2=1$ on $x_2=0$, then equation~\eqref{eqdefiC2}
has a $2\pi\periode$-periodic solution if and only
if $\kappa'$ is given by
$$
\kappa'\defn \frac{1}{\la \xT^2\ra }
\int_{0}^{2\pi\periode} \exp \left(\frac{1}{2i\pi} \int_{0}^{x_2}\beta(s)\, ds\right)
F_{2}(x_2) \, dx_2.
$$
We are then in position to define a function $C_2$ such that $c_3$ is bi-periodic.

\smallbreak
\noindent {\sc Step }6: $\kappa'\in\xR$.
It remains to prove that $\kappa'\in \xR$. Firstly, observe that similar arguments to those
used in the third step establish that
\begin{equation}\label{symak}
\overline{a_k(x_1,x_2)}=a_k(x_1,-x_2)\quad \text{ for } 0\le k\le 4.
\end{equation}
Then, we successively check that
\begin{equation}\label{symck}
\begin{aligned}
\overline{c_0(x_1,x_2)}&=c_0(x_1,-x_2), \\
\overline{c_1(x_1,x_2)}&=c_1(x_1,-x_2), \\
\overline{c_2(x_1,x_2)}&=c_2(x_1,-x_2).
\end{aligned}
\end{equation}
To complete the proof we express $\kappa'$ as 
a function of these coefficients.

\begin{lemm}There holds
\begin{equation}\label{compk'}
\begin{split}
\kappa'&=\frac{1}{\la\xT\ra^2}\iint_{\xT^2}
\frac{c_2}{c_0}\left( \frac{i}{2\pi}\partial_{x_2}\gamma + a_1
-\frac{1}{2}\partial_{x_1}a_0 -\frac{1}{4\nu}a_0^2 - \kappa\right)\, dx_1dx_2\\
& +\frac{1}{\la\xT\ra^2}\iint_{\xT^2} \frac{c_1}{c_0}\left(a_2 - \partial_{x_1}a_1 +\partial_{x_1}^2 a_0\right)\, dx_1dx_2
\end{split}
\end{equation}
\end{lemm}
\begin{proof}
Introduce
$$
\gamma_1 = \frac{c_1}{c_0},\quad \gamma_2=\frac{c_2}{c_0},\quad \gamma_{3}=\frac{c_3}{c_0}.
$$
With this notation, directly from the equation $\delta_3=\delta_3'$ we have
\begin{align*}
\kappa'&=-2\nu \partial_{x_1}\gamma_3-\frac{\beta}{2\pi}\gamma_2+\frac{i}{2\nu}(\partial_{x_2}\gamma)\gamma_2 \\
&\quad -\nu \partial_{x_1}^2\gamma_2-\frac{1}{2}a_0 \partial_{x_1}\gamma_2-\frac{1}{4\nu}a_0^2 \gamma_2
+\frac{1}{2}(\partial_{x_1}a_0) \gamma_2-\kappa\gamma_2 \\
&\quad +a_3 -\partial_{x_1} a_2 +\partial_{x_1}^2 a_1 -\partial_{x_1}^3 a_0
+ \gamma_1 a_2 -\gamma_1 \partial_{x_1}a_1 + \gamma_1 \partial_{x_1}^2 a_0\\
&\quad + \gamma_2 a_1 -\gamma_2 \partial_{x_1}a_0 ,
\end{align*}
where we used various cancelations. By integrating over $\xT^2$, taking into accounts obvious cancelations of the form
$\int_{-\pi}^{\pi} \partial_{x_1} f \, dx_1=0$ and integrating by parts in $x_1$ the term $\iint a_0 \partial_{x_1}\gamma_2 \, dx_1 dx_2$, we obtain
the desired identity.
\end{proof}
Using \eqref{compk'}, \eqref{symak} and \eqref{symck}, we obtain $\overline{\kappa'}=\kappa'$ and hence $\kappa'\in \xR$.

\smallbreak
This completes the proof of the proposition, and hence of Theorem~\ref{theo:smooth}.
\end{proof}

\section{The small divisor condition for families of diamond waves}\label{SIP}

In this section, we prove Proposition~\ref{theo:uniform}, whose statement is recalled here. 

\begin{prop}
Let $\delta$ and $\delta'$ be such that
$$
1>\delta>\delta'>0.
$$
Let $\nu=\nu(\eps)$, $\kappa_0=\kappa_{0}(\eps)$ and
$\kappa_1=\kappa_1(\eps)$ be three real-valued functions defined on $[0,1]$ such that
\begin{equation}\label{nukappa02}
\begin{aligned}
&\nu(\eps)= \underline\nu+ \underline\nu' \eps^2 +\eps\varphi_1(\eps^2),\\
&\kappa_0(\eps)= \underline{\kappa_0}+ \varphi_2(\eps^2) ,\\
&\kappa_1(\eps)= \underline{\kappa_1}+\varphi_3(\eps^2),
\end{aligned}
\end{equation} 
for some constants $ \underline\nu, \underline\nu',\underline{\kappa_0},\underline{\kappa_1}'$ with
$$
 \underline\nu'\neq 0,
$$
and three Lipschitz functions $\varphi_j \colon [0,1]\rightarrow \xR$ satisfying $\varphi_j(0)=0$. 
Assume in addition that
there exists $n\geq 2$ such that
\begin{equation}\label{hd2}
\la \k_2 -\underline{\nu}\k_1^2 -\underline{\kappa_{0}}\ra \geq \frac{1}{\k_1^{1+\delta'}},
\end{equation}
for all $\k \in \xN^2$ with $\k_1\geq n$.
Then there exist $K>0$, $r_0>0$, $N_{0}\in\xN$ and a set $\mathcal{A}\subset [0,1]$ satisfying
$$
\forall r \in [0,r_0],\quad \frac{1}{r}\la  \mathcal{A} \cap [0,r] \ra \ge 1- K r^{\frac{\delta-\delta'}{3+\delta'}},
$$
such that, if $\eps^2\in\mathcal{A}$ and
$\k_1\geq N_{0}$ then
\begin{equation*}
\la \k_2 -\nu(\eps)\k_1^2 -\kappa_0(\eps)-\frac{\kappa_1(\eps)}{\k_1^2}\ra
\geq \frac{1}{\k_1^{2+\delta}},
\end{equation*}
for all $k_2\in \xN$.
\end{prop}
\begin{proof}
As in \cite{IP}, it is convenient to introduce
$$
\lambda =\eps^2,~
\tilde\nu (\lambda)=\nu (\sqrt{\lambda}),~
\tilde\kappa_{0}(\lambda)=\kappa_0(\sqrt{\lambda}),~
\tilde\kappa_1 (\lambda)=\kappa_1(\sqrt{\lambda}).
$$
Introduce also the notations
$$
d(k)=d(k_1,k_2)\defn\frac{k_2}{k_1^2}-\underline{\nu} -\frac{\underline{\kappa_0}}{k_1^2}-\frac{\underline{\kappa_1}}{k_1^4}
$$
and
$$
\omega(\lambda,k_1)
\defn \tilde\nu(\lambda)-\underline{\nu}+\frac{\tilde\kappa_0(\lambda)-\underline{\kappa_0}}{k_1^2}+\frac{\tilde\kappa_1(\lambda)-\underline{\kappa_1}}{k_1^4}.
$$

By assumption,
$$
\omega(\lambda,k_1) = \underline{\nu}'\lambda + \lambda^\mez \varphi_1(\lambda) + 
\frac{\varphi_2(\lambda)}{k_1^2}+\frac{\varphi_3(\lambda)}{k_1^4},
$$
where $\varphi_j(\lambda)$ ($j=1,2,3$) is a Lipschitz function vanishing at the origin. 
Therefore, for $k_1$ large enough and $\lambda$ small enough, we have
\begin{equation}\label{omega1}
\frac{\la\underline{\nu}'\ra}{2} \lambda \le \la\omega(\lambda,k_1)\ra\le  2\la\underline{\nu}'\ra \lambda.
\end{equation}
Similarly, for $k_1$ large enough and $\lambda,\lambda'$ small enough,
\begin{equation}\label{omega2}
\frac{\la\underline{\nu}'\ra}{2}  \la\lambda -\lambda'\ra  \le  \la \omega(\lambda,k_1)-\omega(\lambda',k_1)\ra.
\end{equation}

Given $r>0$ and $N\in \xN^*$, introduce the set
$$
\mathcal{A}(r,N)=\bigcap_{k_2\in\xN}\bigcap_{\xN\ni k_1\ge N}
\left\{\lambda \in [0,r]\, :\, 
\la d(k_1,k_2)-\omega(\lambda,k_1)\ra \ge \frac{1}{k_1^{4+\delta}}  \right\}.
$$
We have to prove that there exist $r_0>0$, $K>0$ and $N_0\in \xN^*$ such that
\begin{equation}\label{todoD}
\forall r\in ]0,r_0],\quad \la [0,r]\setminus \mathcal{A}(r,N_0)\ra \le K  r^{1+\frac{\delta-\delta'}{3+\delta'}}.
\end{equation}

Note that for $n\in \xN^*$,
\begin{equation*}
[0,r]\setminus \mathcal{A}(r,N)= \bigcup_{k_1\ge N}\bigcup_{k_2\in \xN}
E(r,k_1,k_2)
\end{equation*}
with
\begin{equation*}
E(r,k_1,k_2)=
 \left\{ \lambda \in [0,r]:
\la  d(k_1,k_2) -\omega(\lambda,k_1) \ra
 < \frac{1}{\k_1^{4+\delta}}  
 \right\}.
\end{equation*}
The proof then depends on the following lemma. 
\begin{lemm}\label{claim3eps}
There exist $r_0>0$, $N_{0}\in \xN$ and three positive constant $A_j$, $j=1,2,3$,such that, 
for all $r\in [0,r_0]$, for all $k_1\ge N_0$ and for all $k_2\in \xN$, we have\begin{enumerate}
\item If $E(r,k_1,k_2)\neq \emptyset$ then
\begin{equation}\label{claimkeps}
\k_1 \ge A_1  r^{-\frac{1}{3+\delta'}}.
\end{equation}
\item 
$\#\left\{ m\in \xN\,:\, E(r,k_1,m)\neq \emptyset \right\} \le  A_2 (r k_1^2+1)$.
\item 
$
\la E(r,k_1,k_2)\ra \le  A_3 k_1^{-4-\delta}.
$
\end{enumerate}
\end{lemm}

Assume this lemma for a moment and continue the proof. 
This lemma implies that, for 
$r\in [0,r_0]$ and $k_1\ge N_0$, we have
$$
\la \bigcup_{k_2\in \xN}
E(r,k_1,k_2)\ra
\le A_2 A_3 \frac{r k_1^2 +1}{k_1^{4+\delta}}.
$$
Consequently,
\begin{equation*}
\begin{aligned}
\la [0,r]\setminus \mathcal{A}(r,N_0)\ra &\le A_2 A_3 \sum_{\k_1\ge A_1 r^{-\frac{1}{3+\delta'}}}
\frac{r k_1^2 +1}{k_1^{4+\delta}}\\
&\le
A \left\{r \times r^\frac{1+\delta}{3+\delta'}+ r^\frac{3+\delta}{3+\delta'}\right\},
\end{aligned}
\end{equation*}
for some positive constant $A$ independent of $r$, 
which yields the desired bound \eqref{todoD}.

\smallbreak

It remains to prove Lemma~\ref{claim3eps}. Suppose that $E(r,k_1,k_2)\neq \emptyset$ 
and chose $\lambda\in E(r,k_1,k_2)$. 
To prove \eqref{claimkeps}, we first use assumption \eqref{hd2} which reads
\begin{equation*}
\frac{1}{\k_1^{3+\delta'}}
\le \la \frac{\k_2}{\k_1^2} -\underline{\nu} - \frac{\underline{\kappa_{0}}}{k_1^2}\ra,
\end{equation*}
for $k_1$ large enough. 
Hence, by definition of $d(k_1,k_2)$,
\begin{equation*}
\frac{1}{\k_1^{3+\delta'}}
\le \la d(k_1,k_2)\ra +  \frac{\la \underline{\kappa_{1}}\ra}{k_1^4}.
\end{equation*}
Now, by the triangle inequality, we have
\begin{align*}
\frac{1}{\k_1^{3+\delta'}}\le
\la d(k_1,k_2)-\omega(\lambda,k_1)\ra
+\la \omega(\lambda,k_1)\ra  +  \frac{\la \underline{\kappa_{1}}\ra}{k_1^4}.
\end{align*}
By definition, if $\lambda\in E(r,k_1,k_2)$ then
the first term in the right hand side is bounded by $k_1^{-4-\delta}$, so
\begin{equation*}
\frac{1}{\k_1^{3+\delta'}}
\le \frac{1}{\k_1^{4+\delta}} + \la \omega(\lambda,k_1)\ra  +  \frac{\la \underline{\kappa_{1}}\ra}{k_1^4}.
\end{equation*}
Hence, since $0<\delta$ and $\delta'<1$, we have
\begin{equation*}
\frac{1}{2}\frac{1}{k_1^{3+\delta'}}\le \la \omega(\lambda,k_1)\ra,
\end{equation*}
for $k_1$ large enough. 
Now we use the simple estimate \eqref{omega1} valid for $k_1$ large enough and $\lambda$ small enough, to obtain
\begin{equation*}
\frac{1}{k_1^{3+\delta'}}\le 4\la\underline{\nu}'\ra  \lambda\le 4\la\underline{\nu}'\ra r.
\end{equation*}
This establishes the first claim in Lemma~\ref{claim3eps}. 

To prove the second claim, note that, given $r>0$ and $k_1$, if 
$E(r,k_1,k_2)\neq \emptyset $ and  $E(r,k_1,k_2')\neq \emptyset$ then exist $\lambda$ and $\lambda'$ such that
$$
\la d(k_1,k_2) - \omega(\lambda,k_1) \ra
 \le \frac{1}{\k_1^{4+\delta}} ,\quad 
 \la  d(k_1,k_2') -\omega(\lambda',k_1)  \ra
 \le \frac{1}{\k_1^{4+\delta}} .
 $$
Now by the triangle inequality, using \eqref{omega1}, this yields
\begin{align*}
&\la \frac{k_2}{k_1^2}-\frac{k_2'}{k_1^2}\ra=\\
&\quad= \la d(k_1,k_2)- d(k_1,k_2') \ra \\
 &\quad\le 
 \la d(k_1,k_2)-\omega(\lambda,k_1)\ra+\la d(k_1,k_2')-\omega(\lambda',k_1) \ra 
   + \la \omega(\lambda,k_1)-\omega(\lambda',k_1) \ra 
\\
   &\quad\le \frac{2}{k_1^{4+\delta}} + 2\la\underline{\nu}'\ra(\lambda +\lambda')
   \le \frac{2}{k_1^{4+\delta}} + 4\la\underline{\nu}'\ra r.
\end{align*}
Multiplying both sides by $k_1^2$ we obtain the bound 
$$
\la k_2-k_2'\ra \le \frac{2}{k_1^{2+\delta}} + 4\la\underline{\nu}'\ra r k_1^2,
$$
and hence
$$
\#\left\{m\in \xN\,:\, E(r,k_1,m)\neq \emptyset \right\} \le  1 + \frac{2}{k_1^{2+\delta}} + 4\la\underline{\nu}'\ra r k_1^2.
$$

It remains to prove the last claim. 
Introduce
$$
\lambda_{min}(r,k_1,k_2)= \inf  E(r,k_1,k_2),\quad 
\lambda_{max}(r,k_1,k_2)= \sup  E(r,k_1,k_2).
$$
Then, using \eqref{omega2} we find
\begin{align*}
\la E(r,k_1,k_2)\ra &\le \lambda_{max}(r,k_1,k_2)-\lambda_{min}(r,k_1,k_2)\\
&\le \frac{2}{\la\underline{\nu}'\ra} \la \omega\big(\lambda_{max}(r,k_1,k_2),k_1\big)-
\omega\big(\lambda_{min}(r,k_1,k_2),k_1\big)\ra\\
&\le  \frac{2}{\la\underline{\nu}'\ra} \la \omega\big(\lambda_{max}(r,k_1,k_2),k_1\big)-d(k_1,k_2)
\ra\\
&\quad +\la d(k_1,k_2)-\omega\big(\lambda_{min}(r,k_1,k_2),k_1\big)\ra\\
&\le \frac{4}{\la\underline{\nu}'\ra}k_1^{-4-\delta}.
\end{align*}
This completes the proof of Lemma~\ref{claim3eps} 
and hence the proof of Proposition~\ref{theo:uniform}.
\end{proof}
\begin{rema}\label{numbertheory}
The previous proof follows essentially the analysis in \cite{IP}. The key difference is that,
in \cite{IP}, the authors need to prove that a diophantine condition of the form
\begin{equation}\label{dd0}
\la \k_2 -\nu(\eps)\k_1^2 -\kappa_0(\eps)\ra
\geq \frac{1}{\k_1^{2}},
\end{equation}
is satisfied for all $\eps^2\in \mathcal{A}$ for some set $\mathcal{A}$ satisfying
$$
\lim_{r\rightarrow 0} \frac{1}{r} \la \mathcal{A}\cap [0,r]\ra =1 .
$$
This corresponds to the case $\delta=0$ of the above theorem (which we precluded by assumption).
Now, observe that in this case the
above analysis only gives
$$
\la [0,r]\setminus \mathcal{A}(r,N_0)\ra \lesssim r^\frac{3}{3+\delta'}.
$$
Then to pass from this bound to $\la [0,r]\setminus \mathcal{A}(r,N_0)\ra = o(r)$, one has to gain extra decay in $r$.
To do this, Iooss and Plotnikov use an ergodic argument.
What makes the proof of the above Theorem simple is that we
proved only that a weaker diophantine condition is satisfied.
(Here ``weaker diophantine condition''
refers to the fact that, if \eqref{dd0} is satisfied then \eqref{df2} is satisfied
for any $\delta\ge 0$.)
In particular, this discussion clearly shows that it is simpler to prove that \eqref{diophantine}
is satisfied for some $\delta>0$ than for $\delta=0$.
This gives a precise meaning to what we claimed in the introduction:
our paradifferential strategy may be used to simplify
the analysis of the small divisors problems.
\end{rema}

\section{Two elliptic cases}

\subsection{When the Taylor condition is not satisfied}\label{Secgn}
Consider a classical $C^2$ solution $(\sigma,\phi)$ of the system
\begin{equation}\label{III}
\left\{
\begin{aligned}
&\DNS\psi =f_1\in C^\infty(\xT^2),\\[0.5ex]
&\fr \sigma+\frac{1}{2}\la \partialx \psi\ra^2
-\frac{1}{2}\frac{(\partialx  \sigma\cdot\partialx \psi)^2} {1+|\partialx  \sigma|^2}
=f_2\in C^{\infty}(\xT^2),
\end{aligned}
\right.
\end{equation}
where $x\in \xT^2$, and $f_1,f_2$ are given $C^\infty$ functions.

Our goal here is to show that the problem is much easier
in the case where the Taylor sign condition is not satisfied.
To make this more precise, set
$$
\mathfrak{a}\defn \fr+V\cdot\partialx \mathfrak{b} \quad\text{with}\quad
\mathfrak{b}\defn\frac{\partialx \sigma \cdot\partialx \psi}{1+|\partialx  \sigma|^2},
\quad
V\defn \partialx \psi -\mathfrak{b} \partialx\sigma.
$$
We prove a {\em local} hypoellipticity result near boundary points $(x,\sigma(x))$
where $\mathfrak{a}<0$.
We prove that, if $\sigma\in H^{s}$ and $\phi\in H^{s}$ for some $s>3$
near a boundary point $(x_0,\sigma(x_0))$ such that
$\mathfrak{a}(x_0)<0$, then
$\sigma,\phi$ are $C^\infty$ near $(x_0,\sigma(x_0))$.
(This can be improved; the result remains valid under the weaker
assumption that $\sigma,\phi\in C^s$ with $s>2$ for $x\in\xT^d$ with $d\geq 1$.
Yet, we will not address this issue.)

The main observation is that, in the case where $\mathfrak{a}<0$,
the boundary problem~\eqref{III} is weakly elliptic.
Consequently, any term which has the regularity of the unknowns
can be seen as an admissible remainder
for the paralinearization of the first boundary condition (that is why we can
localize the estimates). In addition, the fact that the problem is weakly elliptic implies that we can
obtain the desired sub-elliptic estimates by a simple integration by parts argument. 
To localize in Sobolev space, we use the following notation:
given an open subset $\omega\subset \xR^d$ and a distribution $u\in \mathcal{S}'(\xR^d)$, we say that
$u\in H^{s}(\omega)$ if $\chi u\in H^{s}(\xR^d)$ for every $\chi\in C^{\infty}_{0}(\omega)$.

\begin{theo}\label{theo:aneg} Let
$s>4$ and consider an open domain $\omega\subset\subset \xT^2$.
Suppose that $(\sigma,\psi) \in H^s(\omega)$ satisfies System~\ref{III} and $\mathfrak{a}(x)<0$ for all $x\in\omega$.
Then, for all $\omega' \subset\subset \omega$, there holds
$(\sigma,\psi) \in H^{s+1/2}(\omega')$.
\end{theo}
\begin{proof}
By using symbolic calculus, we begin by observing that we have
a localization property.
Consider two cutoff functions $\chi'\in C^\infty_0 (\omega')$ and
$\chi\in C^{\infty}_{0}(\xR^d)$
such that $\chi=1$ on $\omega$ and $\chi'=1$ on~$\omega'$.
Then $\tilde u = \chi' \psi - T_{\chi  \mathfrak{b}  } \chi' \sigma $ and $\tilde\sigma = \chi'\sigma$
satisfy
\begin{align}
&T_{\dns^1 } \tilde u -  T_V\cdot \partialx \tilde \sigma =\psi\in H^{s}(\omega'),\label{m2.7}
\\[0.5ex]
&T_\mathfrak{a} \tilde \sigma +  T_V \cdot\partialx \tilde u =\theta\in H^{2s - 3}(\omega'),
\label{m2.8}
\end{align}
where recall that
\begin{equation*}
\dns^1(x, \xi)= \sqrt{(1+|\partialx \sigma(x)|^2)| \xi |^2
-\left(\partialx \sigma(x) \cdot \xi\right)^2}.
\end{equation*}
The strategy of the proof is very simple:
We next form a second order equation from \eqref{m2.7}-\eqref{m2.8}.
The assumption $\mathfrak{a}(x_0)<0$ implies that the operator thus
obtained is quasi-homogeneous elliptic, which implies the desired
sub-elliptic regularity for System~\eqref{m2.7}--\eqref{m2.8}. Namely, we claim that
$\tilde{u} \in H^{\alpha + \mez}(\omega')$ and $\tilde{\sigma} \in H^{\alpha }(\omega')$
with
$$
\alpha \defn \min \big\{  s + \mez,  2 s - 3\big\}>s.
$$
To prove this claim, we set $\Lambda = (1 - \Deltax  )^{\mez}$ and use the G\aa rding's inequality for
\paradif operators, to obtain that
there are constants $C$ and $c > 0$ such that
\begin{gather*}
\Re  \big( T_V \cdot\partialx \tilde\sigma , \Lambda^{2\alpha} \tilde u \big)_{L^2} +
 \big( T_V \cdot\partialx \tilde u  , \Lambda^{2\alpha} \tilde\sigma \big)_{L^2}
 \le    C  \big\|  \tilde u   \big\|_{H^{\alpha}}   \big\|  \tilde\sigma    \big\|_{H^{\alpha}} , \\
 c  \big\|  \tilde u   \big\|^2_{H^{\alpha+ \mez }}   \le  \Re  \big( T_{\dns^1} \tilde u  ,
 \Lambda^{2\alpha} \tilde u \big)_{L^2}
+      C  \big\|  \tilde u   \big\|^2_{H^{\alpha}},\\
 c  \big\|  \tilde \sigma    \big\|^2_{H^{\alpha }}
 \le  \Re  \big( T_\mathfrak{a} \tilde \sigma   , -\Lambda^{2\alpha} \tilde \sigma \big)_{L^2}
+      C  \big\|  \tilde \sigma    \big\|^2_{H^{\alpha- \mez}}.
\end{gather*}
Therefore, taking the scalar product of the equations \eqref{m2.7} and
\eqref{m2.8} by $\Lambda^{2\alpha} \tilde u$ and $-\Lambda^{2\alpha} \tilde\sigma $
respectively, and adding the
real parts, implies that
\begin{align*}
&c\big\lVert  \tilde u   \big\rVert ^{2}_{H^{\alpha+ \mez }}
+c\big\lVert  \tilde\sigma \big\rVert^2_{H^{\alpha}}\\
&\quad \le  C  \big\lVert  \theta    \big\rVert_{H^{\alpha-\mez}}
\big\lVert \tilde u     \big\rVert_{H^{\alpha+\mez}}
+ C \big\lVert  \psi    \big\rVert_{H^{\alpha}}
\big\lVert  \tilde \sigma    \big\rVert_{H^{\alpha}}\\
&\quad\quad +  C\big\lVert  \tilde u   \big\rVert_{H^{\alpha}}
\big\lVert  \tilde\sigma    \big\rVert_{H^{\alpha}}
+    C\big\lVert  \tilde u   \big\rVert^2_{H^{\alpha}}
+ C\big\|  \tilde\sigma    \big\|^2_{H^{\alpha- \mez}},
\end{align*}
and the claim follows.

As a consequence we find that
$\sigma \in H^{\alpha}(\omega')$  and  $u\in H^{ \alpha + \mez}(\omega') $.
Going back to $\psi = u + T_{\mathfrak{b} }\sigma$, we obtain that
$\psi \in H^{\alpha}(\omega')$.
This finishes the proof of the Theorem~\ref{theo:aneg}.
\end{proof}

\subsection{Capillary gravity waves}\label{secst}
In this section, we prove {\em a priori} regularity for the
system obtained by adding surface tension:
\begin{equation*}
\left\{
\begin{aligned}
&\DNS\psi-c\cdot\partialx\sigma=0,\\
&\fr \sigma+c\cdot\partialx  \psi+ \frac{1}{2}\la\partialx \psi\ra^2  -\frac{1}{2}
\frac{\bigl(\partialx  \sigma\cdot\partialx \psi +c\cdot\partialx \sigma\bigr)^2}{1+|\partialx  \sigma|^2}
-H(\sigma)
= 0,
\end{aligned}
\right.
\end{equation*}
where $H(\sigma)$ denotes the mean curvature of the free surface $\{y=\sigma(x)\}$:
\begin{equation*}
H(\sigma)\defn \cnx \left(\frac{\partialx\sigma}{\sqrt{1+|\partialx\sigma|^2}}\right).
\end{equation*}
Recently there have been some results concerning {\em a priori} regularity for the
solutions of this system,
first for the two-dimensional case by Matei \cite{Matei},
and second for the general case $d\geq 2$ by
Craig and Matei~\cite{CM,CM2}. Independently, there is also
the paper by Koch, Leoni and Morini~\cite{KLM} which is motivated
by the study of the Munford--Shah functional.
Craig and Matei proved $C^\omega$ regularity for $C^{2+\alpha}$ solutions, and Koch, Leoni and Morini
proved this result for $C^2$  solutions (they also note that the result holds true for $C^1$ viscosity solutions).
Both proofs rely upon the hodograph and Legendre transforms introduced in this context by
Kinderlehrer, Nirenberg and Spruck in the well known papers \cite{KN,KNSI,KNSII}.
Here, as a corollary of Theorem~\ref{theo:paraDN},
we prove that $C^3$ solutions are $C^\infty$, without change of variables,
by using the hidden ellipticity given by surface tension. 
To emphasize this point, the following result is stated in a little more generality than is needed.
\begin{prop}\label{theo:H}
If $(\sigma,\psi)\in C^3(\xR^d)$ solves a system of the form
\begin{equation*}
\left\{
\begin{aligned}
&\DNS\psi = f_1 \in C^{\infty}(\xR^d),\\
&F(\psi,\partialx \psi,\sigma,\partialx\sigma)
+H(\sigma)=f_{2}\in C^{\infty}(\xR^d),
\end{aligned}
\right.
\end{equation*}
where $F$ is a smooth function of its arguments,
then $(\sigma,\psi)\in C^\infty(\xR^d)$.
\end{prop}
\begin{proof}
By using standard regularity results for quasi-linear elliptic PDE, we prove
that if $(\sigma,\psi)\in C^m$ for some $m\geq 2$, then
$(\sigma,\psi)\in C^{m+1-\eps}$ for any $\eps>0$.
For instance, it follows from Theorem 2.2.D in \cite{Taylor2} that,
$$
\sigma\in C^{m}, ~H(\sigma)\in C^{1} \Rightarrow \sigma\in C^{m+1-\delta},
$$
for any $\delta>0$.
As a result,
it follows from the paralinearization formula for the Dirichlet to Neumann operator
(cf Remark~\ref{rH} after Theorem~\ref{theo:paraDN}) that,
$$
T_{\dns^1}\bigl(\psi-T_{\mathfrak{b}}\sigma)\in
C^{m-\delta'}\quad\text{for any}\quad \delta'>0,
$$
where $\dns^1$ is the principal symbol of the Dirichlet to Neumann operator.
Since $\dns^1$ is a first-order elliptic symbol with regularity at least $C^1$ in $x$, this implies that
$\psi-T_{\mathfrak{b}}\sigma\in C^{m+1-\delta''}$ and hence
$\psi\in C^{m+1-\delta''}$ for any $\delta''>0$.
\end{proof}

\end{document}